\documentclass[a4paper, leqno, 11pt]{amsart}

\usepackage[applemac]{inputenc}
\usepackage{aeguill}
\usepackage[english,frenchb]{babel}

\swapnumbers

\numberwithin{equation}{section}

\newtheorem{theo}[equation]{Th\'eor\`eme}
\newtheorem{lem}[equation]{Lemme}
\newtheorem{cor}[equation]{Corollaire}
\newtheorem{prop}[equation]{Proposition}
\newtheorem{fact}[equation]{}

\theoremstyle{definition}
\newtheorem{definition}[equation]{Definition}
\newtheorem{remarque}[equation]{Remarque}

\newcommand{\NN}{{\mathbb N}}
\newcommand{\ZZ}{{\mathbb Z}}

\newcommand{\RR}{{\mathbb R}}
\newcommand{\CC}{{\mathbb C}}

\newcommand{\Reel}{\mathrm{Re}}
\newcommand{\Imag}{\mathrm{Im}}

\newcommand{\Un}{{\mathbf 1}}

\newcommand\cF{{\mathcal F}}

\newcommand\cS{{\mathcal S}}
\newcommand\cC{{\mathcal C}}

\newcommand\cD{{\mathcal D}}
\newcommand\cL{{\mathcal L}}
\newcommand\cN{{\mathcal N}}
\newcommand\cA{{\mathcal A}}
\newcommand\cB{{\mathcal B}}
\newcommand\cE{{\mathcal E}}

\let\wh=\widehat
\let\wt=\widetilde

\let\b\bra
\let\k\ket

\begin{document}
\noindent\vskip 1in

\begin{center}

\begin{Large}
\textbf{Entrelacement de co-Poisson}\\
\vskip 1cm
Jean-François Burnol\\
\end{Large}
\end{center}

\vskip 1cm

\begin{verbatim}
     0 Note importante (important notice)
     1 Introduction
       1.A Sommes
       1.B Propriété de support
       1.C Co-sommes
       1.D Mellin et dzêta
       1.E Fonctions entières et méromorphes
     2 Docteur Poisson et Mister Co
       2.A Des théorèmes de co-Poisson
       2.B Lemmes sur les sommes et les co-sommes
       2.C Preuve du théorème 2.4
       2.D Un théorème de Poisson presque sûr
       2.E Formule intégrale de co-Poisson
       2.F Sommes de Riemann
       2.G Un autre théorème de co-Poisson ponctuel
     3 Études sur une formule de Müntz
       3.A Dzêta et Mellin
       3.B Distributions tempérées et formule de Müntz
       3.C La transformation de Fourier de la fonction dzêta
       3.D Fonctions de carrés intégrables
     4 Entrelacement et fonctions méromorphes
       4.A Convolution multiplicative
       4.B Le théorème d'entrelacement
       4.C Transformation de Mellin
       4.D Propriété S et transformées de Mellin entières
       4.E Fonctions modérées et propriété S
       4.F Distributions homogènes et quasi-homogènes
       4.G Propriété S-étendue et fonctions méromorphes
       4.H Exemples
     Références
\end{verbatim}

\clearpage
  \thispagestyle{empty}
 
\ \vfill

\begin{otherlanguage}{english}

  \centerline{IMPORTANT NOTICE -- PLEASE READ}

   \medskip

   The author is thankful to the creator and maintainers of
   the ``arXiv'' repository system of research papers; he
   does want to point out that any amount of money provided
   by the (especially, French) institutions to help in the
   maintenance of the computers and networks, should
   definitely not be taken by the (especially, French)
   institutions as a method of claim to having contributed,
   even symbolically, and even in a minor way, to the
   research by this author, in the years past and present. 

   \medskip
  
   I shall not name the culprit French institution here, as I
   do not want it at this time to be associated with my work
   even in the most minor of ways. The past is never a past
   thing, especially when it deeply affects the distant
   future. Of course it is individuals who bear the real
   responsability for decisions, and they can trust that I do
   neither forget nor forgive. But I think I am acting
   responsibly in expecting justice and reparation to come
   from the institutional bodies which have that power. So
   far there has been neither justice nor reparation.
  
   \medskip

   It is fundamentally important to this author that the
   readers of his arXived papers become aware of this fact.
   The author would take it as a personal offense if he were
   to see some day in the future logos or other institutional
   stamps from this French institution electronically added
   to the front page or headers or footers of his
   manuscripts.

   \medskip
  
   The present statement is to be considered retroactively
   included into all the previous manuscripts the author has
   uploaded to the arXiv system.
\end{otherlanguage}
   \medskip

   ce 14 juillet 2004,

   \medskip
   Jean-François Burnol
 
   \vfill
 
\clearpage
\setcounter{page}1


\title{Entrelacement de co-Poisson}

\author{Jean-Fran\c cois Burnol}

\date{31 août 2004}

\address{Université Lille 1, UFR de Mathématiques, Cité
  Scientifique M2, F-59655 Villeneuve d'Ascq cedex, France}

\email{burnol@math.univ-lille1.fr}

\maketitle

\clearpage

\tableofcontents

\clearpage

\section{Introduction}

\subsection{Sommes}

Nous débuterons avec l'identité suivante:
\begin{equation}\label{eq:0}
\cF\left(\sum_{n\in\,\ZZ} \delta_n\right) = \sum_{m\in\,\ZZ}
\delta_m
\end{equation}
La transformation de Fourier $\cF$ est normalisée suivant:
\[ \cF(f)(y) = \wt f(y) = \int e^{2\pi i
  xy}f(x)dx = f^\vee(-y)\]
On a donc $(\wt{f}\;)^\vee = f = \wt{f^\vee}$, et $\wt{\wt
  f\;}(x) = f(-x)$.  Bien sûr, $\delta_n(x) = \delta(x-n)$
est la distribution de Dirac positionnée en $n$.

L'identité \eqref{eq:0} est une identité de distributions
tempérées. Suivant une idée de Jean-Pierre
Kahane\footnote{Lettre de J.-P. Kahane à l'auteur, 22 mars
  2002.}, nous allons la régulariser par convolution.
Notons $D(x)$ la \og distribution de Poisson\fg\ (aussi
appelée peigne de Dirac):
\[ D(x) = \sum_{n\in\,\ZZ} \delta_n(x) \]
On a l'implication:
\[ \cF(D(\cdot))(y) = D(y)\quad\Rightarrow\quad \cF(D(\cdot
- t))(y) = e^{2\pi i yt} D(y)\]
En moyennant avec une fonction $f(t)$, il vient:
\[ \cF\left(\int f(t)D(\cdot - t)dt\right)(y) = \int f(t)e^{2\pi i
  yt}dt D(y) = \wt f(y) D(y) \]
Soit encore:
\[ \sum_{n\in\,\ZZ} f(x-n)\quad \stackrel\cF\longrightarrow \quad
\wt f(y) D(y) = \sum_{m\in\,\ZZ} \wt f(m)\delta_m(y)\]
Nous régularisons maintenant le terme de droite:
\[ e^{-2\pi i
  ux}\sum_{n\in\,\ZZ} f(x-n)\quad
\stackrel\cF\longrightarrow \quad \sum_{m\in\,\ZZ} \wt
f(m)\delta_m(y-u)\]
\[ g^\vee(x)\sum_{n\in\,\ZZ} f(x-n)  \stackrel\cF\longrightarrow
\sum_{m\in\,\ZZ} \wt f(m)\int \delta_m(y-u)g(u)\,du =
\sum_{m\in\,\ZZ} \wt f(m)g(y-m)\]
Au final nous obtenons:
\begin{equation}\label{eq:phipsi}
\phi(x) =  g^\vee(x)\sum_{n\in\,\ZZ}
f(x+n)\stackrel\cF\longrightarrow \psi(y) = \sum_{m\in\,\ZZ}
\wt f(m)g(y-m)
\end{equation}
La méthode peut être justifiée par des théorèmes classiques
sur la convolution des distributions \cite{schwartz,
  hormander}, en faisant certaines hypothèses sur les
fonctions $f$ et $g$. Ici $f\in L^1(\RR,dx)$ et $g\in
L^1(\RR,dx)$ suffisent. Et plus généralement pour n'importe
quelle identité initiale $\cF(D_1) = D_2$ de distributions
tempérées l'hypothèse \og pour tout $N\in\NN$ $(1+x^2)^N
f\in L^1$ et $(1+x^2)^N g\in L^1$ \fg\ est suffisante. Si
$f$ et $g$ sont toutes deux dans la classe de Schwartz $\cS$
alors $\phi$ et $\psi$ le seront aussi. Il est important de
préciser cependant que l'identité \eqref{eq:phipsi} n'est a
priori qu'une identité de distributions tempérées, et que
les sommes définissant $\phi$ et $\psi$ sont prises en ce
sens.  Si $\phi$ et $\psi$ sont intégrables et continues
alors \eqref{eq:phipsi} se traduit par des identités
ponctuelles au sens classique.  Exprimons dans un tel cas la
formule d'inversion $\phi = \psi^\vee$, en supposant que la
permutation soit licite:
\[ \psi^\vee(x) = \int_\RR e^{-2\pi i
  xy} \left(\sum_{m\in\,\ZZ} \wt f(m)g(y-m)\right)\,dy =
\left(\sum_{m\in\,\ZZ} \wt f(m)e^{-2\pi i xm}\right)
g^\vee(x)\; ,\]
d'où en comparant avec \eqref{eq:phipsi}:
\begin{equation}\label{eq:poissonx}
 \sum_{n\in\,\ZZ} f(x+n) = \sum_{m\in\,\ZZ} \wt f(m)e^{-2\pi
  i xm} \; ,
\end{equation}
et la \emph{formule sommatoire de Poisson}:
\begin{equation}\label{eq:poisson}
 \sum_{n\in\,\ZZ} f(n) = \sum_{m\in\,\ZZ} \wt f(m)\;.
\end{equation}

Ces dernières étapes sont certainement justifiées si par
exemple $f$ est une fonction de la classe de Schwartz $\cS$,
mais alors la méthode suivie est follement indirecte,
puisqu'en un certain sens, l'équation \eqref{eq:0}, c'est
justement la validité de \eqref{eq:poisson} pour les
fonctions de $\cS$! Plus précisément \eqref{eq:0} se lit:
$\forall\ g\in\cS\ \sum_n \wt g(n) = \sum_m g(m)$, si l'on
suppose qu'une distribution (tempérée) est une fonctionnelle
linéaire sur $\cS$; elle l'est mais elle n'est pas que \c
ca: comme nous le verrons sous peu, il n'y a pas intérêt à
totalement identifier \eqref{eq:0} et \eqref{eq:poisson}.

Un théorème de Katznelson \cite{katznelson} nous avertit que
la formule de Poisson \eqref{eq:poisson} n'a pas une
validité universelle: il donne un exemple de fonctions
continues et intégrables $f$ et $g = \wt f$, telles que les
deux séries de \eqref{eq:poisson} soient absolument
convergentes avec deux sommes distinctes (en fait $f(0)=1$,
$f(n)= 0$ pour $n\in\ZZ, n\neq0$, $g(m) = 0$ pour
$m\in\ZZ$). Toutefois même si les fonctions $f$ et $\wt f$
ne sont pas dans $L^1(\RR)$ (avec $\wt f$ défini
ponctuellement par une intégrale impropre par exemple) il
est possible sous certaines hypothèses (par exemple, de
variation bornée, comme dans le Théorème 45 du livre de
Titchmarsh sur l'intégrale de Fourier \cite{titchfourier}),
d'établir la validité de \eqref{eq:poisson} avec des sommes
parfois seulement semi-convergentes (même après avoir
apparié $n$ et $-n$, $m$ et $-m$).  Un autre cas classique
de validité ponctuelle de \eqref{eq:poisson}, en fait de
\eqref{eq:poissonx}, avec cette fois des sommes absolument
convergentes, est, comme il est aisément établi, obtenu en
supposant que $f$ et $\wt f$ sont continues et toutes deux
bornées par une fonction paire intégrable $k$, décroissante
sur $[0,+\infty[$ (\textsl{cf.} la preuve de
\ref{prop:copoisson1}, page \pageref{prop:copoisson1}).

\subsection{Propriété de support}

Revenons à un aspect intéressant de l'équation
\eqref{eq:phipsi}, qui a disparu en la réduisant à la
formule sommatoire de Poisson \eqref{eq:poisson}. Supposons
que la fonction $f$ (que nous prendrons infiniment
dérivable) soit supportée dans l'intervalle $[-b,b]$ avec
$0<b<\frac12$. Et supposons que la fonction $g$, elle aussi
infiniment dérivable, soit également supportée dans
l'intervalle $[-b,b]$. Alors les deux fonctions $\phi$ et
$\psi$, de $\cS$, paire de Fourier, sont supportées dans
$[-b,b]+\ZZ$. Considérons maintenant la fonction $\phi_1(x)
= \exp(\pi i x)\phi(x-\frac12)$. On a $\wt \phi_1(y) =
i\exp(\pi i y)\psi(y+\frac12)$. Nous constatons que ces deux
fonctions de Schwartz $\phi_1$ et $\psi_1 = \wt \phi_1$ sont
simultanément nulles dans l'intervalle $]-a,+a[$ avec
$a=\frac12-b$. Si $f$ et $g$ ne sont pas toutes deux
identiquement nulles alors ce sera aussi le cas de $\phi_1$
et $\psi_1$.

Cette méthode de construction se heurte à la contrainte
$a<\frac12$. Suivant la lettre de Jean-Pierre Kahane, pour
obtenir des paires de Fourier avec la condition d'être
simultanément nulles dans $]-a,a[$, pour $a>0$ arbitraire,
on commence par modifier la distribution de Poisson $D(x)$
pour annuler un certain nombre de ses Dirac. Soit tout
d'abord $E(x) = i\exp(-\pi i x)D(x+\frac12)$, de sorte que
$\wt E(y) = - \exp(-\pi i y )\wt D(y-\frac12) = iE(y)$.
Remplaçons $E$ par $E_1(x) = \prod_{0\leq j<N} (x^2 -
(\frac12 + j)^2) E(x)$.  La transformée de Fourier de $E_1$
a son support en dehors de $]-\frac12,+\frac12[$, tandis que
$E_1$ elle-même a son support en dehors de
$]-N-\frac12,+N+\frac12[$. Si nous régularisons comme
ci-dessus (ce qui mènera à une formule plus compliquée que
\eqref{eq:phipsi} et que nous n'expliciterons pas), nous
obtenons des fonctions dans $\cS$ (non nulles généralement)
avec des conditions de support à peine moins bonnes. En
remplaçant $\phi(x)$ par $\phi(\sqrt N\;x)$ on obtiendra
l'objectif voulu, avec un $a>0$ arbitraire. Si $\phi$ n'est
pas nulle, alors soit sa partie paire, soit sa partie
impaire ne l'est pas. Si l'on dérive une fonction paire on
obtient une fonction impaire, et cela est compatible avec la
condition de support. Donc il existe pour tout $a>0$ des
fonctions dans $\cS$ non-nulles, paires ou impaires, nulles
et de Fourier nulles dans $]-a,a[$ (on pourra de plus les
prendre vecteurs propres de Fourier, pour l'une de ses
valeurs propres $1$, $i$, $-1$, $-i$, si l'on veut).

Nous ne savons pas si ce théorème a déjà été mentionné dans
la littérature. Un article de de~Branges de 1964
\cite{bra64} a comme corollaire l'existence de fonctions,
paires ou impaires, dans $L^2(\RR)$, avec la propriété de
support, pour $a>0$ quelconque. Pour une démonstration très
courte d'existence dans ce cadre $L^2$ il suffit de dire que
l'espace $L^2(-a,a;dx) + \cF\left(L^2(-a,a;dx)\right)$ est
fermé dans $L^2(\RR)$, et l'on trouvera une preuve simple de
cette dernière affirmation par exemple dans \cite{dymkean}.
Une fois connue une paire $L^2$, la régularisation par
convolution additive avec des fonctions test permet d'en
déduire une paire dans $\cS$ possédant une propriété de
support à peine moins bonne.  Cependant cela ne donne pas
d'exemples véritablement explicites, si le point de départ
lui-même n'est pas explicite.  La méthode de Jean-Pierre
Kahane pallie à cet inconvénient, car son point de départ
est explicite. Pour l'amusement du lecteur (de l'auteur
surtout), nous donnerons à la fin de l'article quelques
exemples anecdotiques explicites (obtenus par une autre
méthode) de fonctions dans $L^2(\RR)$ avec la condition de
support.

Dans \cite{bra64} un lien est fait par de~Branges, via la
transformation de Mellin, entre la condition de support pour
des paires de Fourier de carrés intégrables, plus
généralement pour des paires de transformées de Hankel, et
l'étude de certains espaces de Hilbert de fonctions
entières; ces \og espaces de Sonine\fg\ vérifient tous les
axiomes de la théorie générale présentée dans \cite{bra},
qui a comme point central un théorème d'existence abstrait
d'un développement isométrique. Il fut étudié de manière
concrète pour ces espaces par de~Branges \cite{bra64} et par
J. et V.~Rovnyak \cite{rov1,rov2}, mais principalement pour
les espaces de Sonine associés aux fonctions de Bessel
d'ordre entier.
Avec l'espace de Hardy d'un demi-plan et les espaces de
Paley-Wiener, les espaces de Sonine de de~Branges--Rovnyak
occupent une place remarquable à l'intersection de l'analyse
complexe et de l'analyse de Fourier.


\subsection{Co-sommes}

Reprenons maintenant cette idée de régularisation d'une
paire de Fourier de distributions tempérées telle que
\eqref{eq:0}. Mais moyennons de manière multiplicative, avec
$t\neq 0$ dans les calculs qui suivent:
\[ D(\frac{x}{t})  \stackrel\cF\longrightarrow |t| D(ty)
\]
\[
\int_{-\infty}^\infty f(t)D(\frac{x}{t}) \frac{dt}{|t|}
\stackrel\cF\longrightarrow \int_{-\infty}^\infty
f(t)D(ty)\,dt \]
Nous séparerons les Dirac en $0$ des autres. En effet,
puisque $\delta(\lambda x) = \frac1{|\lambda|}\delta(x)$,
pour $\lambda\neq 0$, on a les identités de distributions en
$x$ (resp. $y$):

\begin{align*}
  \int_{-\infty}^\infty f(t)\delta(\frac{x}{t} -
  n)\frac{dt}{|t|} &=
\begin{cases}
  \frac{f(x/n)}{|n|}&(n\neq0)\\ \left(\int_{-\infty}^\infty
    f(t)\, dt \right)\delta(x)&(n=0)
\end{cases}\quad ,\\
\int_{-\infty}^\infty f(t)\delta(ty - m)\,dt &=
\begin{cases}
  \frac{f(m/y)}{|y|}&(m\neq0)\\ \left(\int_{-\infty}^\infty
    f(t)\frac{dt}{|t|}\right)\delta(y)&(m=0)
\end{cases}\quad ,
\end{align*}
qui donnent:
\[ \left(\int_{-\infty}^\infty f(t)\, dt \right)\delta(x) +
\sum_{n\neq0} \frac{f(x/n)}{|n|} \stackrel\cF\longrightarrow
\left(\int_{-\infty}^\infty
  f(t)\frac{dt}{|t|}\right)\delta(y)+
\sum_{m\neq0}\frac{f(m/y)}{|y|}\;,\]
puis finalement:

\noindent\framebox[\textwidth][c]{%
  \parbox[c]{0.9\textwidth}{
\begin{equation}\label{eq:copoisson}
\sum_{n\neq0} \frac{f(x/n)}{|n|} - \int_{-\infty}^{+\infty}
f(t)\frac{dt}{|t|} \stackrel\cF\longrightarrow
\sum_{m\neq0}\frac{f(m/y)}{|y|} - \int_{-\infty}^{+\infty}
f(t)\, dt \; .
\end{equation}}}

\smallskip

La formule ne concernant en réalité que les fonctions
paires, nous pouvons donc la réécrire sous la forme
suivante, où il est entendu que $f$ est une fonction paire
(ou que $f$ est définie sur $]0,\infty[$ et que $\cF$ est la
transformation en cosinus $\cF(F)(y) = \int_0^\infty
2\cos(2\pi yx) F(x)\,dx$\;):
\[
\sum_{n\geq1} \frac{f(x/n)}{n} - \int_0^{+\infty}
f(t)\frac{dt}{t} \stackrel\cF\longrightarrow
\sum_{m\geq1}\frac{f(m/y)}{|y|} - \int_0^{+\infty} f(t)\, dt
\; .
\]

La méthode suivie est aisément justifiable telle quelle si
par exemple $f$ est une fonction intégrable, à support
compact, éloigné de l'origine, et le résultat est alors a
priori une identité de distributions tempérées, que nous
appellerons \textbf{entrelacement de co-Poisson}. On peut
prouver que \eqref{eq:copoisson} vaut au sens des
distributions sous la seule condition:
\[ \int_{\RR}
|f(t)|\,dt+\int_{\RR}\frac{|f(1/t)|}{|t|}\,dt<\infty\leqno{(C)}\]
Il s'agit là d'un théorème dû à Duffin et Weinberger (1997)
\cite{dw2} (sous une hypothèse additionnelle sur $f$; sous
cette seule hypothèse il apparaît dans \cite[Thm 4.2]{fz}).
Sous la seule hypothèse $(C)$, Duffin et Weinberger avaient
aussi donné antérieurement dans \cite{dw1} (1991) un énoncé
(voir plus loin \ref{theo:dw}$.(2)$, page \pageref{theo:dw})
qui exprime également que les deux termes de
\eqref{eq:copoisson} forment une paire de Fourier en un
certain sens généralisé, plus \og classique\fg. Citons dès
maintenant l'un des résultats que nous établirons ici dans
ce contexte:

\emph{ Soit $f$ une fonction mesurable telle que $
  \int_0^\infty (1+\frac1t)|f(t)|dt<\infty $. Soit $\xi\geq0$
  et soit $X>\xi$. On a, pour $\Lambda\to\infty$:
\[\begin{split}
  \int_0^\Lambda 2\cos(2\pi\xi x)\left(\sum_{n\geq1}
    \frac{f(n/x)}x - \int_0^\infty f(u)\,du\right)\,dx=\\
  \int_0^X \frac{\sin(2\pi\Lambda(t-\xi))}{\pi(t-\xi)}
  \sum_{n\geq1} \frac{f(t/n)}n\,dt - \int_0^\infty
  \frac{f(u)}u \,du + o(1)
\end{split}
\]
}

Ce théorème ramène donc, en toute généralité, la question de
la validité ponctuelle de la formule de co-Poisson à des
questions classiques portant sur le noyau de Dirichlet.

La formule suivante, cousine de \eqref{eq:copoisson}, fut
publiée par Duffin dès 1945! (\cite{duff}; l'apparence de la
formule dépend de la normalisation choisie pour la
transformation de Fourier)
\begin{equation}\label{eq:copoissonduffin}
\begin{split}
  \frac{f(2x)}{1/2} - \frac{f(2x/3)}{3/2} +
  \frac{f(2x/5)}{5/2}  - \dots \qquad\\
  \stackrel\cF\longrightarrow i\left( \frac{f(1/(2y))}{|y|}
    - \frac{f(3/(2y))}{|y|} +
    \frac{f(5/(2y))}{|y|}-\dots\right) \;.
\end{split}
\end{equation}
Elle vaut pour $f$ \emph{impaire}, et l'on pourrait la
déduire de l'identité de distributions impaires
\[ \wt E = iE\qquad E = \sum_{n\in\ZZ}
(-1)^{n-1}\delta_{n-\frac12}\;,\]
de la même façon que nous avons déduit ici
\eqref{eq:copoisson} de $\wt D = D$. Des généralisations de
la formule de Duffin \eqref{eq:copoissonduffin} furent
établies par Weinberger dans sa thèse de 1950 \cite{wein}.
La formule principale \eqref{eq:copoisson} n'apparaît
semble-t-il que plus tard en 1991 (\cite{dw1}).


Lorsque $f$ est choisie de classe $C^\infty$ à support
compact (éloigné de l'origine), les équations
\eqref{eq:copoisson} et \eqref{eq:copoissonduffin}
concernent des fonctions qui sont dans la classe de Schwartz
$\cS$ (comme on le voit aisément grâce à \dots\ la formule
de Poisson!). De cette façon on obtient des paires de
Fourier dans $\cS$ qui sont constantes (nulles, si l'on
veut) dans un intervalle $]-a,a[$ avec $0<a<1$ arbitraire
(avec \eqref{eq:copoissonduffin} on ne peut faire mieux que
$0<a<\frac12$).  Et si l'on modifie préalablement la
distribution de Poisson $D$ ou la variante $E$ (ou une autre
variante $\sum_{n\in\ZZ} c_n \delta_n$ avec
$(c_n)_{n\in\ZZ}$ périodique en $n$), en annulant certains
de ses Dirac comme dans la lettre de J.-P.~Kahane, on
obtient alors des exemples pour tout $a>0$.  Tacitement nous
avons utilisé la remarque simple suivante:

\medskip

\noindent\framebox[\textwidth][c]{%
\parbox[c]{0.9\textwidth}{%
  \emph{\textbf{Remarque simple}: si $f$ est nulle pour
    $0<t<a$, alors la co-somme $\sum_{n\geq1} f(t/n)/n$, et
    plus généralement les co-sommes $\sum_{n\geq1} c_n
    f(t/n)/n$, sont également nulles pour $0<t<a$.}%
}%
}

\medskip

Il est intéressant de remarquer que des sommes de ce type
apparaissent dès l'article de Riemann (plus précisément, des
sommes $\sum_{n\geq1} \frac{c_n}n F(x^{\frac1n})$). Le
calcul de leurs transformées de Fourier a été traité, pour
la première fois semble-t-il, par Duffin \cite{duff} et par
Weinberger \cite{wein} (pour $(c_n)_{n\in\ZZ}$ périodique,
de somme nulle sur une période). Ils appellent la formule
\eqref{eq:copoisson} la \emph{formule de Poisson dualisée}
\cite{dw1,dw2}; comme il y a de multiples liens entre la
formule de Poisson \eqref{eq:poisson} et la formule
\eqref{eq:copoisson}, nous optons plutôt pour la
terminologie \emph{formule de co-Poisson}.

La \og remarque simple\fg\ ci-dessus établit un lien entre
les formules de Duffin-\hbox{Weinberger} et les espaces de
de~Branges-Rovnyak. Ce lien, ainsi que certaines
constructions qui en résultent, furent mis en avant
(peut-être pour la première fois) dans la conférence
d'habilitation de l'auteur de décembre 2001 \cite{fz} (et,
plus implicitement, dans son article antérieur \cite{crashp}).

Nous avons déjà indiqué que la formule de Poisson
\eqref{eq:poisson} permet d'obtenir des informations sur les
sommes de Riemann dans la formule de co-Poisson
\eqref{eq:copoisson}. Donnons ici un exemple de théorème
portant sur la formule de Poisson que nous démontrons en
utilisant la formule de co-Poisson, et que l'on comparera au
contre-exemple de Katznelson \cite{katznelson}:

\emph{Soit $f$ et $\wt f$ une paire de Fourier de fonctions
  intégrables. La formule:
\[ 
\sum_{n\in\,\ZZ} \frac{f(n/x)}{|x|} = \sum_{m\in\,\ZZ} \wt
f(mx)\;,
\]
concerne des séries qui sont absolument convergentes pour
presque tout $x$ et elle est une identité pour presque tout
$x$.}

\medskip


\subsection{Mellin et dzêta}

Une formule de M\"untz (\cite{muntz};
\cite[II.11]{titchzeta}) rattache les sommes de Poisson à la
fonction dzêta de Riemann via la transformation de Mellin;
nous en donnerons une étude détaillée, en étendant fortement
ses conditions de validité, en examinant certains aspects
qui lui sont associés dans le contexte des distributions, et
en donnant une modification permettant de l'exprimer dans un
contexte de fonctions de carrés intégrables.  Nous mettrons
en évidence que l'équation fonctionnelle de la fonction
dzêta s'exprime via la formule de co-Poisson de manière tout
aussi naturelle que via la formule de Poisson; l'idée est de
distinguer entre les \emph{deux} formes de transformées de
Mellin $\int_0^\infty f(t)t^{s-1}\,dt$ et $\int_0^\infty
f(t) t^{-s}\,dt$, et donc d'espérer au moins {deux
  interprétations de l'équation fonctionnelle}: et
effectivement si l'une est Poisson, alors l'autre sera
co-Poisson.  Il est intéressant de remarquer que
les \emph{deux} transformées de Mellin sont présentes dans
l'article de Riemann.

\subsection{Fonctions entières et méromorphes}

Une dernière partie est consacrée aux distributions
tempérées vérifiant la condition de support initialement
considérée par de~Branges pour les fonctions. 

\emph{Soit $a>0$. Soit $D$ une distribution tempérée paire
  non nulle qui s'annule sur $]-a,a[$. On peut donner un
  sens à $\wh{D}(s) = \int_0^\infty D(t)t^{-s}dt$ comme
  fonction analytique de $s$ pour $\Reel(s)\gg 1$. Si
  $\cF(D)$ est à nouveau nulle dans $]-a,a[$ alors
  $\wh{D}(s)$ est une fonction entière de $s$, qui a des
  zéros triviaux en $-2n$, $n\in\NN$ et vérifie l'équation
  fonctionnelle
\[\pi^{-\frac s2}\Gamma(\frac s2)\wh{D}(s) =
\pi^{-\frac{1-s}2}\Gamma(\frac{1-s}2)\wh{\cF(D)}(1-s)\;.\]
Elle a dans le secteur $|\arg(s-\frac12) -
\frac\pi2|<\epsilon$ ($0<\epsilon<\frac\pi2$) un nombre de
zéros de modules au plus $T$ asymptotiquement équivalent à
$\frac{T}{2\pi}\log(T)$, et dans $|\arg(s)|<\frac\pi2
-\epsilon$ les zéros de modules au plus $T$ sont en nombre
$o(T)$.}


\emph{Les fonctions entières $F(s)$ de la forme $\wh D(s)$
  où $D$ est une distribution tempérée nulle et de Fourier
  nulle dans un voisinage de l'origine sont les fonctions
  $F(s)$ qui sont $O(A^{\Reel(s)}(1+|s|)^N)$ dans tout
  demi-plan $\Reel(s)\geq \sigma$, (pour $A>0$ et $N\in\NN$
  dépendant de $\sigma$), et telles que $\chi(s)F(1-s)$ soit
  aussi entière et avec la même propriété dans ces
  demi-plans ($\chi(s) = \zeta(s)/\zeta(1-s)$).}

Une condition de support plus générale est explicitée qui
correspond à la possibilité pour $\wh D(s)$ d'avoir un
nombre fini de pôles, et nous décrivons alors la partie
polaire de $\pi^{-\frac s2}\Gamma(\frac s2)\wh{D}(s)$ en
fonction de la distribution $D$.

Ces outils permettraient de suivre dans le détail tous les
calculs qui dans \cite{crassonine} ont abouti à la
définition de certaines distributions paires (et leurs
analogues impaires) $A_a$ (invariante par $\cF$) et $B_a$
(anti-invariante par $\cF$) qui ont la propriété de support
pour l'intervalle $]-a,a[$.
Nous les évoquerons brièvement, mais une discussion plus
approfondie nécessiterait un travail spécialement dédié à
cet effet, trop indirectement lié à la formule de co-Poisson
pour avoir trouvé sa place ici (voir aussi
\cite{crasdirac}).  Nous conclurons donc plutôt avec
quelques exemples explicites amusants de fonctions de carrés
intégrables avec la propriété de support.

\section{Docteur Poisson et Mister Co}

\subsection{Des théorèmes de co-Poisson}

Soit $f(x)$ une fonction paire mesurable. Posons:
\begin{subequations}
\begin{align}\label{eq:F}
  F(x) &= \sum_{n\geq1} \frac{f(n/x)}{|x|} - \int_0^\infty
  f(u)\,du\\ \label{eq:K} K(x) &= \sum_{n\geq1}
  \frac{f(x/n)}{n} - \int_0^\infty \frac{f(1/u)}{u}\,du
\end{align}
\end{subequations}
On pose aussi:
\[F(0) = - \int_0^\infty f(u)\,du\qquad K(0) = -\int_0^\infty
\frac{f(1/u)}{u}\,du\]
Nous verrons que la somme dans $F$ (resp. $K$) est presque
partout absolument convergente dès que $\int_0^\infty
|f(u)|\,du<\infty$ (resp. $\int_0^\infty
\frac{|f(1/u)|}{u}\,du<\infty$).  Remplacer $f$ par
$f(1/x)/|x|$ revient à interchanger $F$ et $K$.
La formation de $K$ commute aux changements d'échelle,
contrairement à la formation de $F$ qui renverse les
changements d'échelle. La fonction $F$ représente l'écart
entre une somme de Riemann et l'intégrale de $f$, et tend
donc vers zéro pour $x\to\infty$ par exemple lorsque $f$ est
R-intégrable et à support compact.  Citons dans ce contexte
un résultat élémentaire ($\cS$ désigne la classe de
Schwartz):
\begin{lem}
  Si $f\in\cS$ et $f(0)=0$ alors $F\in\cS$.
\end{lem}
En effet, on voit d'abord directement sur \eqref{eq:F} que
$F$ est $C^\infty$ y-compris en $x=0$. De plus la formule de
Poisson \eqref{eq:poisson} permet d'écrire $F(x) =
\sum_{m\geq1} \wt f(mx)$ et donc d'obtenir la décroissance
rapide de $F$ pour $x\to\infty$. Si l'on veut que la somme
$\sum_{n\geq1} f(n/x)/|x|$ ait une chance d'être en règle
générale de carré intégrable sur $[0,+\infty[$, il
\emph{faut} donc la modifier en $F$.

On utilisera souvent la condition suivante sur la fonction
paire $f$:
\[
\int_0^\infty |f(x)|\left( 1 + \frac1x\right)\,dx <
\infty\leqno{(C)}
\]

\begin{theo}[Duffin-Weinberger]\label{theo:dw}
  Soit $f$ une fonction paire mesurable vérifiant $(C)$.  On
  a:
\begin{enumerate}
\item \cite[Th. 1]{dw1} Les sommes dans $F$ et $K$ sont
  presque partout absolument convergentes et elles sont
  convergentes dans $L^1(0,\Lambda;dx)$ pour tout
  $\Lambda<\infty$.
\item \cite[Th. 1, éq. 3.8]{dw1} Les fonctions $F$ et $K$
  sont une paire de Fourier en le sens généralisé:
  \begin{align*}
    \int_0^\Lambda F(x)\,dx &= \int_0^{\to\infty}
    \frac{\sin(2\pi\Lambda t)}{\pi t} K(t)\,dt \\
    \int_0^T K(t)\,dt &= \int_0^{\to\infty} \frac{\sin(2\pi
      Tx)}{\pi x} F(x)\,dx
  \end{align*}
\item \cite[Lem. 2]{dw2} Si $f$ est de variation bornée et
  si $\int_0^\infty (1+x)|df|(x)<\infty$ alors $F$ et $K$
  sont une paire de Fourier au sens des distributions.
\item \cite[Th. 1]{dw2} Si $f$ est $C^1$, si $f'$ est abs.
  cont.  et si $\int_0^\infty |f''(x)|\,dx<\infty$ alors $F$
  est $L^1$, $K$ est partout définie et est continue, et
  \begin{equation*}
    \forall\xi\quad\int_0^\infty 2\cos(2\pi\xi x) F(x)\,dx = K(\xi)
  \end{equation*}
\end{enumerate}
\end{theo}

Nous allons montrer:

\begin{theo}\label{theo:copoisson}
  Soit $f$ une fonction paire mesurable vérifiant $(C)$. On
  a:
  \begin{enumerate}
  \item $F$ et $K$ sont tempérées en tant que distributions
    et sont une paire de Fourier en ce sens.
  \item On a: $\int_0^{\to\infty} F(x)\,dx = \frac12 K(0)$.
  \item Soit $\Lambda\geq0$ et soit $\xi\geq0$. On a:
    \begin{multline*}
      \qquad\int_0^\Lambda 2\cos(2\pi \xi x) F(x)\,dx
      =\\
      \int_0^{\to\infty}
      \left(\frac{\sin(2\pi\Lambda(t-\xi))}{\pi(t - \xi)} +
        \frac{\sin(2\pi\Lambda(t+\xi))}{\pi(t +\xi)}\right)
      K(t)\,dt
    \end{multline*}
  \item Soit $\xi\geq0$ et soit $X>\xi$.  On a, pour
    $\Lambda\to\infty$:
    \begin{multline*}
      \qquad\int_0^\Lambda 2\cos(2\pi\xi x)F(x)\,dx=\\
      \int_0^X \frac{\sin(2\pi\Lambda(t-\xi))}{\pi(t-\xi)}
      \sum_{n\geq1} \frac{f(t/n)}n\,dt - \int_0^\infty
      \frac{f(u)}u \,du + o(1)
    \end{multline*}
\end{enumerate}
\end{theo}

Le point \ref{theo:copoisson}$.(1)$ est énoncé par
Duffin-Weinberger dans \cite{dw2} sous une hypothèse plus
forte sur $f$ (\textsl{cf.}  \ref{theo:dw}$.(3)$); il est
prouvé sous la seule hypothèse $(C)$ dans \cite{fz}
(théorème $4.2.$). Le point \ref{theo:copoisson}$.(2)$ est
établi à la fin de la démonstration du théorème 4.6 dans
\cite{fz}, avec un argument qui n'utilise en fait que $(C)$.
Le point \ref{theo:copoisson}$.(3)$ est une extension du
théorème principal \cite[Th. 1]{dw1} de Duffin-Weinberger
\textsl{cf.}  \ref{theo:dw}$.(2)$. Finalement
\ref{theo:copoisson}$.(4)$ résout la question de la validité
ponctuelle en un sens classique de la formule de co-Poisson,
puisqu'il a comme corollaire immédiat:

\begin{cor}\label{cor:1}
  Soit $f$ une fonction paire mesurable vérifiant $(C)$.
  L'intégrale impropre $\int_0^{\to\infty}2\cos(2\pi\xi
  x)F(x)\,dx$ existe si et seulement si $\xi$ est un point
  de Dirichlet de la fonction (loc. intégrable) définie
  presque partout $K$, et sa valeur est alors la valeur de
  Dirichlet de $K$ en $\xi$.
\end{cor}

En effet si $\phi$ est localement intégrable on dit que
$x_0$ est un point de Dirichlet de $\phi$ si $L =
\lim_{\Lambda\to\infty} \int_{x_0-\delta}^{x_0+\delta}
\frac{\sin(2\pi\Lambda(x-x_0))}{\pi(x-x_0)}\phi(x)\,dx$
existe ($\delta>0$ fixé petit quelconque), et la valeur de
Dirichlet de $\phi$ en $x_0$ est par définition $L$.
Remarquons que par \ref{theo:copoisson} le point $\xi=0$ est
toujours un point de Dirichlet de $K$ (il s'agit là en fait
d'un lemme que nous établirons préalablement à la preuve de
\ref{theo:copoisson}). Comme corollaire, on obtient un
critère, non restreint au cadre $L^1$ et à des fonctions
continues, qui est suffisant pour la validité ponctuelle de
la formule de co-Poisson sur tout un intervalle:

\begin{theo}\label{cor:2}
  Soit $f$ une fonction paire mesurable vérifiant $(C)$ et
  telle que, pour un certain $X>0$, l'on ait:
  \begin{enumerate}
  \item $f$ est de variation bornée sur $[0,X]$,
  \item aux discontinuités dans $]0,X[$, $f$ a la valeur
    moyenne.
  \end{enumerate}
  Alors:
  \[\forall \xi\in\; ]0,X[\qquad
  \int_0^{\to\infty} 2\cos(2\pi\xi x) F(x)\,dx = K(\xi)
  \]
\end{theo}

Prouvons ce corollaire. Comme $f$ a une variation totale
finie sur $[0,X]$, elle a une limite en $0^+$, et cette
limite doit être nulle à cause de $(C)$. On peut donc aussi
bien poser $f(0) =0$. Soit $\mu$ la mesure complexe $df$ de
sorte que $f(t) = \mu([0,t[) + \frac12\mu(\{t\})$ pour
$t<X$. Soit $|\mu|$ la mesure des variations de $\mu$. Soit
$0<t_1<t_2< X$ fixés et soit $V_m$, pour $m\geq1$, la
variation totale de $f(t/m)/m$ sur l'intervalle $[t_1,
t_2]$. On a $V_m \leq |\mu|([t_1/m, t_2/m])/m$. Nous
obtenons alors:
\[ \sum_{m\geq1} V_m \leq \int_{[0,t_2]} \sum_{t_1/u\leq
  m\leq t_2/u} \frac1m d|\mu|(u)\;, \]
qui est
\[\leq C+ \int_{[0, t_1/2]} \log\frac{t_2/u}{t_1/u - 1}
d|\mu|(u) \leq C + \log(2t_2/t_1)|\mu|([0, t_1/2]) <\infty\]
Nous savons par \ref{theo:dw}$.(1)$ (voir lemme \ref{lem:k})
que la série $\sum_m f(t/m)/m$ converge presque partout donc
elle converge en au moins un point $t_0$ de $[t_1,t_2]$.  La
série $\sum_m (f(t/m)/m - f(t_0/m)/m)$ est uniformément de
Cauchy donc elle converge ponctuellement partout et
uniformément sur $[t_1,t_2]$ (qui est arbitraire). La
fonction $\sum_m f(t/m)/m$ a sa variation totale sur
$[t_1,t_2]$ bornée par $\sum_{m\geq1} V_m$ qui est fini.
Comme la série converge uniformément la règle de la valeur
moyenne aux discontinuités est vérifiée. Par le critère de
Jordan si une fonction $\phi$ est de variation bornée dans
un voisinage d'un point $\xi$, alors $\xi$ est un point de
Dirichlet de $\phi$ et la valeur de Dirichlet $\phi_D(\xi)$
est $\frac12(\phi(\xi+) + \phi(\xi-))$. Le critère de Jordan
s'applique donc et le corollaire est démontré.

Un autre théorème assurant une validité ponctuelle de la
formule de co-Poisson est donné dans la section
\kern-1em\ref{subsec:pointwise}.

\subsection{Lemmes sur les sommes et les co-sommes}

Ces lemmes serviront pour la démonstration de
\ref{theo:copoisson}.

\begin{lem}\label{lem:k}
  Soit $f\in L^1(0,\infty;dx)$. La série $\sum_{n\geq1}
  \frac{f(n/x)}x$ est presque partout absolument convergente
  et est absolument convergente dans $L^1(0,\Lambda;dx)$
  pour tout $\Lambda>0$. Soit $k\in L^1(0,\infty;dx)$
  décroissante (donc à valeurs positives ou nulles). On a
  \[ \int_0^\infty \sum_{n\geq1}  \frac{|f(n/x)|}x k(x)\,dx
  \leq \int_0^\infty |f(x)|\,dx \int_0^\infty k(x)\,dx\;.\]
\end{lem}

La convergence ponctuelle et en norme $L^1$ est un lemme de
Duffin et Weinberger dans \cite{dw1}, que l'on peut voir
comme une conséquence de l'inégalité intégrale, puisque l'on
pourra y prendre pour $k$ la fonction indicatrice de
l'intervalle $]0,\Lambda[$ par exemple. On utilisera aussi
plus loin $k(x) = 1/(1+x^2)$. Par le théorème de la
convergence monotone on a
\[  \int_0^\infty \sum_{n\geq1}  \frac{|f(n/x)|}x k(x)\,dx
= \int_0^\infty \sum_{n\geq1} |f(x)|\frac{k(n/x)}x \,dx
\;.\]
Or $\forall x>0\ \sum_{n\geq1} \frac{k(n/x)}x\leq
\int_0^\infty k(x)\,dx$. Le lemme est prouvé.

\begin{lem}\label{lem:g}
  Soit $g\in L^1(0,\infty;dx)$. La série $\sum_{m\geq1}
  g(mx)$ est absolument convergente dans
  $L^1(0,\infty;\frac{dx}{1+(\log x)^{2}})$, donc aussi
  presque partout et dans $L^1(a,A;dx)$ pour tout
  $0<a<A<\infty$.
\end{lem}

On a en effet:
\[ \int_0^\infty \sum_{m\geq1}  |g(mx)|
\frac{dx}{1+(\log x)^{2}} = \int_0^\infty \sum_{m\geq1}
\frac{|g(m/x)|}x \frac{dx}{x(1+(\log x)^{2})}\;,\] et l'on
applique le lemme précédent (la fonction $x^{-1}(1+(\log
x)^2)^{-1}$ n'est pas décroissante sur tout $]0,\infty[$
mais on peut certainement la majorer par une fonction
décroissante et intégrable).

\begin{lem}\label{lem:dirichlet}
  Soit $X>0$ et soit $f\in L^1(0,X;\frac{dt}t)$. On a:
  \[ \lim_{\Lambda\to\infty}  \int_0^X \frac{\sin(2\pi\Lambda
    t)}{\pi t}\sum_{m\geq1} \frac{f(t/m)}m\,dt = 0\]
\end{lem}

On peut poser $f\equiv0$ pour $t>X$. On a $f(t/m)/m =
f_1(m/t)/t$ avec $f_1(t)=f(1/t)/t$, donc $f_1\in
L^1(0,\infty;dt)$. Par \ref{lem:k} $\sum_{m\geq1}
\frac{f(t/m)}m$ converge absolument presque partout et dans
$L^1(0,X;dt)$.  Le lemme est certainement vrai si $f$ est
identiquement nulle sur $]0,\epsilon[$, $\epsilon>0$, par le
lemme de Riemann-Lebesgue; de plus il existe une constante
absolue telle que:
\[ \left| \int_0^X \frac{\sin(2\pi\Lambda t)}{\pi
    t}\sum_{m\geq1} \frac{f(t/m)}m\,dt \right| \leq C
\int_0^X \frac{|f(u)|}u\,du \]
En effet
\begin{align*}
  &\int_0^X \frac{\sin(2\pi\Lambda t)}{\pi
    t}\sum_{m\geq1} \frac{f(t/m)}m\,dt\\
  &= \sum_{m\geq1} \int_0^X \frac{\sin(2\pi\Lambda t)}{\pi
    t}\frac{f(t/m)}m\,dt \\
  &= \sum_{m\geq1} \int_0^{X/m} \frac{\sin(2\pi\Lambda
    mu)}{\pi
    m}f(u)\frac{du}u \\
  &= \int_0^X \sum_{m\geq1} \Un_{[0,X]}(mu)
  \frac{\sin(2\pi\Lambda mu)}{\pi m}f(u)\frac{du}u
\end{align*}
On a utilisé le théorème de la convergence dominée et le
fait bien connu que les sommes partielles $\sum_{1\leq m\leq
  M}\sin(2\pi mv)/\pi m$ sont bornées uniformément en $M\in\NN$,
$v\in\RR$. Cette borne absolue $C$ est celle qui apparaît
alors dans l'inégalité plus haut qui est ainsi prouvée. En
combinant ces deux observations on obtient:
\[ \forall\epsilon>0\quad \limsup_{\Lambda\to\infty} \left|
  \int_0^X \frac{\sin(2\pi\Lambda t)}{\pi t}\sum_{m\geq1}
  \frac{f(t/m)}m\,dt \right| \leq C \int_0^\epsilon
\frac{|f(u)|}u\,du \] Comme $\epsilon>0$ est arbitraire le
lemme est démontré.

\subsection{Preuve du théorème 2.4}

\begin{prop} 
  \label{prop:copoisson1}
  Soient $\phi$ et $\psi$ deux fonctions paires continues et
  intégrables sur $\RR$ et qui sont transformées de Fourier
  l'une de l'autre. On suppose de plus qu'il existe une
  fonction $k:[0,+\infty[\to\RR$ décroissante et intégrable
  qui majore à la fois $\phi$ et $\psi$. Soit $f$ une
  fonction paire mesurable vérifiant $(C)$. Alors on a
  l'égalité intégrale de co-Poisson:
  \[ \int_0^\infty \phi(x)F(x)\,dx = \int_0^\infty
  \psi(y)K(y)\,dy \;,\]
  où $F$ et $K$ sont définies selon \eqref{eq:F} et
  \eqref{eq:K}.
\end{prop}

Les deux intégrales sont absolument convergentes par le
lemme \ref{lem:k}. Le même lemme justifie:
\[ \int_0^\infty \phi(x)\sum_{n\geq1} \frac{f(n/x)}x \,dx
= \int_0^\infty \sum_{n\geq1} \frac{\phi(n/x)}x f(x)\,dx \]
La formule de Poisson
\[ \phi(0) + 2\sum_{n\geq1} \phi(n) = 
\psi(0) + 2\sum_{m\geq1} \psi(m) \; , \]
vaut. En effet $\sum_{n\geq0} \phi(n+\alpha)$ est
uniformément absolument convergente pour $0\leq\alpha\leq1$
puisque dominée par $\sum_{n\geq0} k(n)$. Il en est de même
de $\sum_{n<0} \phi(n+\alpha)$. Donc $A(\alpha) =
\sum_{n\in\ZZ} \phi(n+\alpha)$ est une fonction continue
$1$-périodique.  Ses coefficients de Fourier sont les
$\psi(-m)$, $m\in\ZZ$, qui forment une suite de $l^1$. La
fonction continue $B(\alpha) = \sum_{m\in\ZZ}
\psi(-m)\exp(2\pi i m\alpha)$ est donc partout égale à
$A(\alpha)$, en particulier en $\alpha=0$. Le même
raisonnement vaut pour $\phi(\cdot/x)/x$ et $\psi(x\cdot)$,
$x>0$, et permet donc d'écrire:
\[ \int_0^\infty \sum_{n\geq1} \frac{\phi(n/x)}x f(x)\,dx
= \int_0^\infty (\sum_{m\geq1}\psi(mx) - \frac{\phi(0)}{2x}
+ \frac12\psi(0))f(x)\,dx\;,\] puis en utilisant $\phi(0)=
2\int_0^\infty \psi(u)\,du$, $\psi(0) = 2\int_0^\infty
\phi(t)\,dt$, et un changement de variable $x=1/y$:
\begin{align*}
  &\int_0^\infty \phi(x)\left(\sum_{n\geq1} \frac{f(n/x)}x -
    \int_0^\infty f(t)\,dt\right)\,dx\\
  =& \int_0^\infty (\sum_{n\geq1} \frac{\phi(n/x)}x -
  \int_0^\infty \phi(t)\,dt) f(x)\,dx\\
  =& \int_0^\infty (\sum_{m\geq1}\psi(mx) -
  \frac{\phi(0)}{2x}) f(x)\,dx\\
  =& \int_0^\infty (\sum_{m\geq1}\frac{\psi(m/y)}y -
  \int_0^\infty \psi(u)\,du) \frac{f(1/y)}y \,dy\;.
\end{align*}
La dernière ligne donne alors par les mêmes calculs sur
$\psi$ et $f(1/y)/y$ que précédemment pour $\phi$ et $f$:
\[ = \int_0^\infty \psi(y)\left(\sum_{m\geq1}\frac{f(y/m)}m -
  \int_0^\infty \frac{f(1/u)}u\,du\right)\,dy\;.\]
Ceci conclut la preuve de la Proposition. Les calculs sont
exactement semblables à ceux de Duffin-Weinberger dans
\cite{dw2} (ou aux calculs menant à \cite[Th. 4.2]{fz}),
seules les hypothèses (et notations) diffèrent.

En particulier, cela établit le point $(1)$ du théorème
\ref{theo:copoisson} (\textsl{cf.} aussi le Lemma $4.1$ de
\cite{fz}). Nous établissons maintenant le point $(4)$.
Soit $\xi\geq0$ fixé et soit $\Lambda>0$.  Considérons la
fonction continue
\[ \phi(x)=
\begin{cases}
  2\cos(2\pi\xi x) - 2\cos(2\pi\Lambda\xi)&|x|\leq
  \Lambda\;,\\
  0&|x|>\Lambda\;.
\end{cases}\]
Sa transformée de Fourier $\psi(y)$ est donnée par
\begin{equation}\label{eq:psi}
\psi(y) = \frac{\sin(2\pi\Lambda(y-\xi))}{\pi(y - \xi)} +
\frac{\sin(2\pi\Lambda(y+\xi))}{\pi(y +\xi)} -
2\cos(2\pi\Lambda\xi)\frac{\sin(2\pi\Lambda y)}{\pi y} 
\end{equation}
On peut la réécrire sous la forme
\begin{multline*}
  \psi(y) = \sin(2\pi\Lambda
  y)\cos(2\pi\Lambda\xi)\left\{\frac1{\pi(y-\xi)}+\frac1{\pi(y+\xi)}
    - \frac2{\pi y}\right\}\\+ \cos(2\pi\Lambda
  y)\sin(2\pi\Lambda\xi)\left\{\frac1{\pi(y+\xi)}-\frac1{\pi(y-\xi)}\right\}\;,
\end{multline*}
qui montre qu'elle est elle aussi $L^1$, car en fait majorée
par un multiple (dépendant de $\xi$ et de $\Lambda$) de
$1/(1+y^2)$. Les hypothèses de la proposition
\ref{prop:copoisson1} s'appliquent.  Ainsi, avec $X>\xi$:
\begin{multline}\label{eq:Jxi}
  J_\xi(\Lambda) := \int_0^\Lambda 2\cos(2\pi\xi x) F(x)\,dx
  = 2\cos(2\pi \Lambda\xi)\int_0^\Lambda F(x)\,dx \\+
  \int_0^X \psi(y)K(y)\,dy + \int_X^\infty
  \psi(y)K(y)\,dy\;.
\end{multline}
Le dernier terme est $o(1)$ lorsque $\Lambda\to\infty$
puisque $\int_0^\infty |K(y)|\frac{dy}{1+y^2}<\infty$, et
que l'on applique le lemme de Riemann-Lebesgue, compte tenu
de l'expression ci-dessus pour $\psi(y)$. On a de plus:
\[ \int_0^X
\psi(y)K(y)\,dy = \int_0^X \psi(y)\sum_{m\geq1}
\frac{f(y/m)}m\,dy - \int_0^X \psi(y)\,dy \int_0^\infty
\frac{f(1/u)}u\,du\]
Compte tenu de \eqref{eq:psi}, on a:
\begin{multline*}
  \int_0^X \psi(y)\,dy = \int_{-X}^X
  \frac{\sin(2\pi\Lambda(y-\xi))}{\pi(y - \xi)} \,dy -
  \cos(2\pi\Lambda\xi)\int_{-X}^X
  \frac{\sin(2\pi\Lambda y)}{\pi y} \,dy\\
  = 1 - \cos(2\pi\Lambda\xi) + o(1)\;,
\end{multline*}
et ainsi:
\begin{multline*}
  J_\xi(\Lambda) = 2\cos(2\pi\Lambda\xi)\left(\int_0^\Lambda
    F(x)\,dx + \frac12\int_0^\infty
    \frac{f(1/u)}u\,du\right) - \int_0^\infty
  \frac{f(1/u)}u\,du\\
  + \int_0^X \psi(y)\sum_{m\geq1} \frac{f(y/m)}m\,dy +
  o(1)\;.
\end{multline*}
La dernière intégrale est, par l'expression \eqref{eq:psi}
de $\psi(y)$, le lemme de Riemann-Lebesgue et le lemme
\ref{lem:dirichlet}:
\[ \int_0^X \psi(y)\sum_{m\geq1}
\frac{f(y/m)}m\,dy = \int_0^X
\frac{\sin(2\pi\Lambda(y-\xi))}{\pi(y - \xi)}\sum_{m\geq1}
\frac{f(y/m)}m\,dy + o(1) \] En combinant tous ces éléments
on obtient:
\begin{multline*}
  J_\xi(\Lambda) = 2\cos(2\pi\xi
  \Lambda)\left(\int_0^\Lambda F(x)\,dx +
    \frac12\int_0^\infty \frac{f(1/u)}u\,du\right) -
  \int_0^\infty
  \frac{f(1/u)}u\,du\\
  + \int_0^X \frac{\sin(2\pi\Lambda(y-\xi))}{\pi(y -
    \xi)}\sum_{m\geq1} \frac{f(y/m)}m\,dy + o(1)
\end{multline*}
Il ne reste plus qu'à évoquer
\[
\lim_{\Lambda\to\infty} \int_0^\Lambda F(x)\,dx =
-\frac12\int_0^\infty \frac{f(1/u)}u\,du\;,\leqno{(E)}
\]
pour que la preuve de \ref{theo:copoisson}$.(4)$ soit alors
complète.

L'équation $(E)$, c'est-à-dire \ref{theo:copoisson}$.(2)$,
est établie dans \cite{fz}, à la fin de la démonstration du
théorème $4.6$. Pour la commodité du lecteur nous
reproduisons ici la preuve qui n'utilise que la seule
condition $(C)$ sur $f$. Notons $g(t) = f(1/t)/t$.

Tout d'abord
\[\int_0^\Lambda \sum_{n\geq1} \frac{g(t/n)}n dt = \sum_{n\geq1}
\int_0^{\Lambda/n} g(t)dt = \int_0^\infty \left[\frac\Lambda
  t\right]g(t)dt\;,\]
donc
\[ \int_0^\Lambda F(t)\,dt = \int_0^\Lambda
\left(\sum_{n\geq1} \frac{g(t/n)}n - \int_0^\infty
  \frac{g(u)}u\,du\right)\,dt\]
\[ =  -\int_0^\infty\left\{\frac\Lambda
  t\right\}g(t)dt = - \int_0^\infty \{\Lambda v\}h(v)dv\;,
\]
avec $h(v) = g(1/v)/v^2 = f(v)/v$ appartenant à
$L^1(0,+\infty;dv)$.  Il est clair que
\[0\leq A\leq B\quad\Rightarrow\quad \lim_{\Lambda\to\infty}
\int_0^\infty \{\Lambda v\}\Un_{A\leq v\leq B}(v)dv =
\frac{B-A}2\] donc par l'argument usuel de densité dans la
preuve standard du lemme de Riemann-Lebesgue on prouve
\[\lim_{\Lambda\to\infty} \int_0^\infty \{\Lambda v\}h(v)dv =
\frac12\int_0^\infty h(v)dv\]
ce qui complète la preuve de \ref{theo:copoisson}$.(2)$ et
donc de \ref{theo:copoisson}$.(4)$.

Venons-en à \ref{theo:copoisson}$.(3)$. Nous avons par
\eqref{eq:Jxi}:
\[ \int_0^\Lambda 2\cos(2\pi \xi x) F(x)\,dx  =
2\cos(2\pi\Lambda\xi)\int_0^\Lambda F(x)\,dx + \int_0^\infty
\psi(y) K(y)\,dy \;.\]
Par le Théorème de Duffin et Weinberger \ref{theo:dw}$.(2)$
(\cite{dw1}), on a:
\[ \int_0^\Lambda
F(x)\,dx = \int_0^{\to\infty} \frac{\sin(2\pi\Lambda y)}{\pi
  y} K(y)\,dy \;,\]
donc:
\[ \int_0^\Lambda 2\cos(2\pi \xi x) F(x)\,dx  =
\int_0^{\to\infty} \left( \psi(y) + 2\cos(2\pi\Lambda\xi)
  \frac{\sin(2\pi\Lambda y)}{\pi y}\right) K(y)\,dy\;, \]
et finalement par \eqref{eq:psi}:
\[ \int_0^\Lambda 2\cos(2\pi \xi x) F(x)\,dx =
\int_0^{\to\infty}
\left(\frac{\sin(2\pi\Lambda(y-\xi))}{\pi(y - \xi)} +
  \frac{\sin(2\pi\Lambda(y+\xi))}{\pi(y +\xi)}\right)
K(y)\,dy\;.\]

\begin{remarque}
  Dans la preuve de \ref{theo:copoisson}$.(4)$ nous nous
  sommes arrangé par la définition spéciale de $\phi(x)$
  pour n'avoir que des intégrales absolument convergentes ce
  qui nous a permis de faire appel à la proposition générale
  \ref{prop:copoisson1}. On peut procéder plus directement
  en définissant $\phi_1(x) = 2\cos(2\pi\xi x)$ pour
  $|x|<\Lambda$, $\phi_1(\pm \Lambda) =
  \cos(2\pi\Lambda\xi)$, $\phi_1(x) = 0$ pour $|x|>\Lambda$.
  Soit alors $\psi_1(y)$ sa transformée de Fourier. La
  fonction $\phi_1(\cdot/x)/x$ ($x>0$) est intégrable et de
  variation bornée et vérifie la règle de la valeur moyenne
  aux discontinuités, on peut donc lui appliquer le Théorème
  de Poisson $45$ de Titchmarsh \cite{titchfourier}. Pour
  poursuivre les calculs et justifier les interversions de
  sommes et d'intégrales on aura besoin d'établir que les
  sommes partielles $\sum_{1\leq m\leq M} y\psi_1(my)$ sont
  bornées uniformément en $M$ et en $y$ (pour $\Lambda$ et
  $\xi$ fixés). On se ramène ensuite aux mêmes
  considérations que dans la preuve exposée ici. Les détails
  sont laissés au lecteur intéressé.
\end{remarque}

\subsection{Un théorème de Poisson presque sûr}

Une conséquence intéressante de la formule de co-Poisson est
le théorème de Poisson presque sûr:
\begin{theo}\label{theo:poissonpp}
  Soient $\phi$ et $\psi$ deux fonctions paires continues et
  intégrables qui sont transformées de Fourier l'une de
  l'autre. Alors, pour presque tout $x>0$:
\[
\sum_{n\in\,\ZZ} \frac{\phi(n/x)}{|x|} = \sum_{m\in\,\ZZ}
\psi(mx)\;,
\]
avec des séries absolument convergentes.
\end{theo}

Nous savons déjà par les lemmes \ref{lem:k} et \ref{lem:g}
que les séries sont presque partout absolument convergentes
et aussi dans $L^1(a,A;dx)$ pour tout $0<a<A<\infty$. Soit
$f$ une fonction paire infiniment dérivable avec son support
(pour $x>0$) dans $[a,A]$, $0<a<A<\infty$. On a
\begin{align*}
  \int_0^\infty |\phi(x)|\sum_{n\geq1} \frac{|f(n/x)|}x\,dx &<\infty\\
  \int_0^\infty |\psi(x)|\sum_{m\geq1} \frac{|f(x/m)|}m\,dx
  &<\infty\;,
\end{align*}
puisque les sommes impliquant $f$ sont bornées. Ainsi le
théorème de la convergence monotone implique
\begin{align*}
  \int_0^\infty \phi(x)\sum_{n\geq1} \frac{f(n/x)}x\,dx &=
  \int_0^\infty \sum_{n\geq1} \frac{\phi(n/x)}x f(x)\,dx
  \\
  \int_0^\infty \psi(x)\sum_{m\geq1} \frac{f(x/m)}m\,dx &=
  \int_0^\infty \sum_{m\geq1}\psi(mx) f(x)\,dx\;.
\end{align*}
Soient $F$ et $K$ définies selon \eqref{eq:F} et
\eqref{eq:K}. Ce sont des fonctions de Schwartz qui sont
transformées de Fourier l'une de l'autre, par la formule de
co-Poisson. On a donc:
\[\int_0^\infty \phi(x)F(x)\,dx = \int_0^\infty
\psi(x)K(x)\,dx \;,\]
que l'on peut écrire compte tenu des équations précédentes
selon:
\begin{align*}
  &\int_0^\infty \sum_{n\geq1} \frac{\phi(n/x)}x f(x)\,dx
  - \int_0^\infty \phi(x)\,dx\int_0^\infty f(x)\,dx\\
  = &\int_0^\infty \sum_{m\geq1}\psi(mx) f(x)\,dx -
  \int_0^\infty \psi(x)\,dx\int_0^\infty \frac{f(x)}x\,dx\;,
\end{align*}
puis finalement, puisque $\int_0^\infty \phi(x)\,dx =
\frac{\psi(0)}2$, $\int_0^\infty \psi(x)\,dx =
\frac{\phi(0)}2$:
\[
\int_a^A \left(\sum_{n\geq1}
  \frac{\phi(n/x)}x+\frac{\phi(0)}{2x} -
  \sum_{m\geq1}\psi(mx) -\frac12\psi(0)\right)f(x)\,dx =
0\;.\] Cela conclut la preuve du théorème.

\subsection{Formule intégrale de co-Poisson}

Nous allons maintenant utiliser ce Théorème de Poisson
(déduit de la formule de co-Poisson) pour prouver\dots\ la
formule de co-Poisson!

\begin{theo}
  Soit $f$ une fonction paire mesurable vérifiant:
  \[
  \int_0^\infty |f(x)|\left( 1 + \frac1x\right)\,dx <
  \infty\leqno{(C)}
  \]
  Soient $\phi$ et $\psi$ deux fonctions intégrables paires
  qui sont transformées de Fourier l'une de l'autre. Si:
  \begin{align*}
    &&\int_0^\infty |\phi(x)|\sum_{n\geq1}
    \frac{|f(n/x)|}x\,dx <\infty \quad \text{ou}\quad
    &\int_0^\infty \sum_{n\geq1} \frac{|\phi(n/x)|}x
    |f(x)|\,dx<\infty
    \\
    &\text{et si}\quad &\int_0^\infty |\psi(x)|\sum_{n\geq1}
    \frac{|f(x/n)|}n\,dx <\infty \quad \text{ou}\quad
    &\int_0^\infty \sum_{n\geq1}|\psi(nx)|\,|f(x)|\,dx
    <\infty\;,
  \end{align*}
  alors la formule intégrale de co-Poisson
  \[\begin{split}
    &\int_0^\infty \phi(x)\left(\sum_{n\geq1}
      \frac{f(n/x)}{x} - \int_0^\infty f(u)\,du\right)\,dx
    \\
    = &\int_0^\infty \psi(x)\left(\sum_{n\geq1}
      \frac{f(x/n)}{n} - \int_0^\infty
      \frac{f(1/u)}{u}\,du\right)\,dx \;,
  \end{split}\]
  est valable (avec des intégrales absolument convergentes).
\end{theo}

\begin{remarque}
  Par le théorème de la convergence monotone les conditions
  ci-dessus séparées par des \og ou\fg\ sont équivalentes et
  impliquent la convergence absolue des intégrales de la
  formule intégrale de co-Poisson. Elles seront vérifiées
  si, par exemple, soit $f(x)$ et $f(1/x)/x$, soit $\phi(x)$
  et $\psi(x)$ sont toutes deux majorées par une fonction
  $k$ décroissante et intégrable sur $]0,\infty[$ (ce
  dernier cas est déjà traité par \ref{prop:copoisson1}).
\end{remarque}

En ce qui concerne la démonstration du théorème, on peut
supposer que $\phi$ et $\psi$ sont choisies continues dans
leur classe d'équivalence. Il suffira ensuite de suivre à
l'identique la preuve de la proposition
\ref{prop:copoisson1}, et d'invoquer au moment opportun le
théorème de Poisson presque sûr \ref{theo:poissonpp}.

\subsection{Sommes de Riemann}

Le théorème de Poisson presque sûr permet aussi d'obtenir
certains résultats sur le comportement de sommes de Riemann
pour des fonctions intégrables.

\begin{prop}\label{prop:A}
  Soit $\phi$ une fonction paire intégrable continue, dont
  la transformée de Fourier $\psi$ vérifie
  \[ \int_2^\infty \log(y)|\psi(y)|\,dy < \infty\]
  Soit:
  \[ A(x) = \sum_{n\geq1} \frac{\phi(n/x)}x - \int_0^\infty
  \phi(t)\,dt\;,\]
  qui est définie pour presque tout $x>0$. On a alors:
  \[ \int^\infty \left| A(x) +
    \frac12\frac{\phi(0)}x\right|\,dx < \infty \]
\end{prop}

En effet, posons $B(x) = \sum_{m\geq1} \psi(mx)$. On vérifie
aisément
\[ \int_2^\infty \log(x)|\psi(x)|\,dx <\infty
\quad\Rightarrow \quad \int_2^\infty |B(x)|\,dx <\infty\;,\]
or, par \ref{theo:poissonpp}, on sait que pour presque tout
$x>0$ on a
\[ A(x) +\frac12\frac{\phi(0)}x = B(x)\]
ce qui complète la preuve.

Sous une hypothèse plus faible, nous obtenons tout de même à
nouveau un résultat qui montre que $A(x)$ tend vers zéro en
un sens moyen lorsque $x$ tend vers plus l'infini:

\begin{prop}\label{prop:A2}
  Soit $\phi$ une fonction paire continue intégrable dont la
  transformée de Fourier $\psi$ est aussi intégrable. Soit:
\[ A(x) = \sum_{n\geq1} \frac{\phi(n/x)}x - \int_0^\infty
\phi(t)dt\;,\] qui est défini pour presque tout $x>0$. On a:
\[ \forall \epsilon>0\qquad \int^\infty |A(x)|\frac{dx}{\log
  x(\log\log x)^{1+\epsilon}} < \infty \;,\] ce qui
implique, pour $\lambda\to\infty$:
\[ \int_\lambda^{\lambda^2} |A(x)| dx =
o(\log\lambda(\log\log\lambda)^{1+\epsilon})\;,\] et aussi
certainement:
\[ \lim_{\lambda\to\infty} \frac1\lambda
\int_\lambda^{2\lambda} |A(x)|\,dx = 0 \; . \]
\end{prop}

Nous avons à nouveau presque partout, par la formule de
Poisson presque sûr:
\[ A(x) +\frac12\frac{\phi(0)}x = \sum_{m\geq1} \psi(mx)\]
et il suffira donc de s'intéresser à $C(x) = \sum_{m\geq1}
|\psi(mx)|$. Définissons
\[ k(x) = \begin{cases}\frac1x\frac1{|\log x|}{(\frac1{\log|\log
      x|})^{1+\epsilon}}&(0<x\leq x_0)\\
  0&(x>x_0)
\end{cases}\;,\]
où $x_0>0$ est choisi suffisamment petit de sorte que $k$
est décroissante. Alors:
\begin{multline*}
  \int_0^\infty \frac{C(1/x)}x k(x)\,dx = \int_0^\infty
  |\psi(x)|\sum_{m\geq1} \frac{k(m/x)}x\, dx \\
  \leq \int_0^\infty |\psi(x)|\,dx\int_0^\infty k(x)\,dx <
  \infty
\end{multline*}
Ceci donne
\[ \int^\infty C(x) \frac{dx}{\log
  x(\log\log x)^{1+\epsilon}} < \infty \;,\] qui établit
l'affirmation principale, et les deux suivantes en sont de
directs corollaires.

\subsection{Un autre théorème de co-Poisson ponctuel}
\label{subsec:pointwise}

Le théorème suivant a le théorème \ref{theo:dw}$.(4)$ de
Duffin-Weinberger comme corollaire.

\begin{theo}\label{theo:casepointwise}
  Soit $f$ une fonction mesurable vérifiant
\[ \int_0^\infty (1+\frac1x)|f(x)|dx<\infty \leqno{(C)}\]
On suppose aussi que sa transformée en cosinus
\[ \widetilde f(y) = \int_0^\infty 2\cos(2\pi xy)f(x)\,dx \]
a la propriété
\[ \int^\infty \log(y)|\widetilde f(y)|\,dy <\infty\;. \]
En particulier $\wt f$ est $L^1$ et $f$ est essentiellement
continue, donc on la supposera continue. Sous ces
hypothèses, la fonction
\[ F(x) = \sum_{n\geq1} \frac{f(n/x)}x - \int_0^\infty
f(x)\,dx \;, \]
qui est définie pour presque tout $x>0$, a la propriété
\[ F \in L^1(0,\infty;dx) \;, \]
la série $\sum_{m\geq1} \frac{f(\xi/m)}m$ converge
uniformément sur tout intervalle $[0,X]$, et l'on a la
formule de co-Poisson ponctuelle:
\[
\forall \xi\geq0\quad \int_0^\infty 2\cos(2\pi\xi
x)F(x)\,dx= \sum_{m\geq1} \frac{f(\xi/m)}m - \int_0^\infty
\frac{f(u)}u \,du\;.
\]
\end{theo}

On notera que par $(C)$ la fonction continue $f$ vérifie
$f(0) = 0$ et que l'intégrabilité de $F$ résulte alors de la
proposition \ref{prop:A}. Nous savons que la formule de
co-Poisson vaut au sens des distributions. Donc il suffira
de prouver que la série $\sum_{m\geq1} \frac{f(\xi/m)}m$
converge uniformément sur $[0,X]$, car sa somme sera une
fonction continue égale presque partout, donc partout, à la
transformée en cosinus de $F$. Compte tenu de $f(0) = 0$, on
peut écrire:
\[\frac{f(\xi/m)}m = \int_0^\infty 2\cos(2\pi \xi y)\wt
f(my)\,dy = \int_0^\infty 2(\cos(2\pi \xi y) - 1)\wt
f(my)\,dy\]
De l'hypothèse $\int^\infty \log(y)|\wt f(y)|\,dy <\infty $,
il résulte:
\[\sum_{m\geq1} \int_1^\infty |\wt f(my)|\,dy < \infty \]
et il suffit donc d'examiner la série
\[ \sum_{m\geq1} \int_0^1 2(1 - \cos(2\pi \xi y))\,|\wt
f(my)|\,dy\;,\]
qui est majorée par
\[ \sum_{m\geq1} \int_0^1 4\pi\xi \,y\,|\wt f(my)|\,dy\;.\]
Or par le lemme \ref{lem:g} on a certainement
\[ \wt f\in L^1\quad\Rightarrow\quad \int_0^\infty
\sum_{m\geq1}|\wt f(my)|\frac{dy}{y+\frac1y} <\infty \;,\]
ce qui implique la convergence de la série précédente et
complète la démonstration du théorème
\ref{theo:casepointwise}.

\section{Études sur une formule de Müntz}

\subsection{Dzêta et Mellin}

Dans son article, Riemann multiplie $\zeta(s)$ par
$\Gamma(s)$ ou par $\pi^{-s/2}\Gamma(s/2)$ pour obtenir des
expressions intégrales permettant d'établir le prolongement
analytique et l'équation fonctionnelle. On peut adopter la
perspective suivant laquelle c'est $\zeta(s)$ qui multiplie
$\Gamma(s)$ ou $\pi^{-s/2}\Gamma(s/2)$, elles-mêmes des
intégrales, et que c'est cette multiplication par $\zeta(s)$
qui a un effet sur les quantités intégrées. Ce point de vue
est donc un point de vue opératoriel sur l'action de
$\zeta(s)$. Dans les deux cas cités les intégrales sont des
transformées de Mellin \textit{gauches} $\int_0^\infty
f(x)x^{s-1}\,dx$. On trouve alors que l'action de $\zeta(s)$
est celle d'une somme $\sum f(nx)$.  Mais dans d'autres cas,
les intégrales sont plus naturellement exprimées par des
transformées de Mellin \textit{droites} $\int_0^\infty
f(x)x^{-s}\,dx$.  Alors, l'action de $\zeta(s)$ est celle
d'une co-somme $\sum f(x/n)/n$. Il y a donc la possibilité,
au moins, de \emph{deux} traductions de l'équation
fonctionnelle via l'interprétation opératorielle de la
multiplication par $\zeta(s)$. Et effectivement, pour la
transformée de Mellin gauche, on trouve que cette traduction
est la formule de Poisson, tandis que pour la transformée de
Mellin droite on trouve la formule de co-Poisson.

Pour correctement expliquer ce qui précède il faut préciser
un point important. La formule
\begin{equation}\label{eq:premuntz}
\int_0^\infty \sum_{n\geq1}
f(nx)\,x^{s-1}\,dx = \zeta(s)\int_0^\infty f(x)x^{s-1}\,dx\;,
\end{equation}
pour $\Reel(s)=\sigma$, $\int_0^\infty
|f(x)|x^{\sigma-1}\,dx<\infty$, ne vaut (en général) que
pour $\sigma>1$ (et est justifiée par le théorème de
Fubini-Tonelli). Si l'on veut interpréter opératoriellement
l'équation fonctionnelle il est indispensable que
l'opération $s\to1-s$ soit licite.  Il nous faut donc une
version de \eqref{eq:premuntz} qui soit valable pour
$0<\sigma<1$ (ou pour $\sigma=\frac12$). Il semble d'après
Titchmarsh \cite[II.11]{titchzeta} que c'est Müntz
\cite{muntz} qui le premier a donné la formule. Voici le
Théorème de Müntz, tel qu'on le trouve dans le traité de
Titchmarsh:

\begin{theo}[M\"untz {\cite{muntz} \cite[II.11]{titchzeta}}]
\label{theo:muntz}
Soit $f$ une fonction de classe $C^1$ sur $[0,\infty[$ qui
est $O(x^{-\alpha})$ pour $x\to\infty$ et telle que aussi
$f'(x) = O(x^{-\beta})$, avec $\alpha>1$, $\beta>1$. On a
alors:
  \begin{equation}\label{eq:muntz}
    \int_0^\infty \left(\sum_{n\geq1} f(nx) -
      \frac{\int_0^\infty f(t)\,dt}{x}\right) x^{s-1}\,dx =
    \zeta(s)\int_0^\infty f(x)x^{s-1}\,dx\;,
  \end{equation}
  pour $0<\Reel(s)<1$.
\end{theo}

Nous appellerons
\begin{equation}\label{eq:A}
  A_f(x) = \sum_{n\geq1} f(nx) - \frac{\int_0^\infty f(t)\,dt}{x}
\end{equation}
la sommation de Poisson (Müntz)-modifiée de $f$.
Plusieurs démonstrations de l'équation fonctionnelle sont
obtenues dans \cite{titchzeta} à partir de la formule de
Müntz, en choisissant $f$ de sorte que $A_f$ ait des
propriétés spéciales, comme par exemple d'être invariante
sous la transformée en sinus, ou sous la transformation
$A_f(x)\to A_f(1/x)/x$.

Avant toute chose, faisons la remarque suivante: le
changement de variable $x\to1/x$ transforme \eqref{eq:muntz}
en la formule équivalente:
  \begin{equation*}
    \int_0^\infty \left(\sum_{n\geq1} \frac{f(n/x)}x -
      \int_0^\infty f(t)\,dt\right) x^{-s}\,dx =
    \zeta(s)\int_0^\infty f(x)x^{s-1}\,dx\;,
  \end{equation*}
  puis en remplaçant $f(x)$ par $f(1/x)/x$ on obtient la
  formule de co-Müntz:
  \begin{equation}\label{eq:comuntz}
    \int_0^\infty \left(\sum_{n\geq1} \frac{f(x/n)}n -
      \int_0^\infty \frac{f(1/t)}t\,dt\right) x^{-s}\,dx =
    \zeta(s)\int_0^\infty f(x)x^{-s}\,dx\;.
  \end{equation}
  Ceci met en évidence que si l'on utilise la transformée de
  Mellin \emph{droite}, la multiplication par $\zeta(s)$ est
  associée aux \emph{co-sommes}.
  
  Que ce soit pour Mellin-droit ou pour Mellin-gauche, la
  droite $\Reel(s)=\frac12$ est spéciale du point de vue
  Hilbertien: les applications $f(x)\mapsto M_g(f)(s) =
  \int_0^\infty f(x)x^{s-1}\,dx$ et $f(x)\mapsto M_d(f)(s) =
  \int_0^\infty f(x)x^{-s}\,dx$ sont isométriques de
  $L^2(0,\infty;dx)$ sur
  $L^2(s=\frac12+i\tau;\frac{d\tau}{2\pi})$. Supposons que
  la fonction paire $f$ de carré intégrable soit sa propre
  transformée de Fourier. Comment cela se traduit-il pour
  $M_d(f)$ et $M_g(f)$?
  
  Soit $\Gamma = \cF I$, avec ici $\cF$ la transformation en
  cosinus, et $I:f(x)\mapsto f(1/x)/x$. L'opérateur $\Gamma$
  est unitaire sur $L^2(0,\infty;dx)$ et il commute aux
  changements d'échelle. Il est donc diagonalisé par $M_d$
  (ou $M_g$), qui n'est pas autre chose que la
  transformation de Fourier multiplicative (et l'on sait
  bien que les opérateurs sur $L^2(-\infty,+\infty;du)$ qui
  commutent aux translations sont diagonalisés par la
  transformation de Fourier additive). Donc il existe une
  fonction mesurable $\chi(s)$ sur la droite critique, de
  module $1$ presque partout, et telle que:
\[\forall  f\in L^2(0,\infty;dx) \quad
\Reel(s)=\frac12\Rightarrow
M_d(\Gamma(f))(s)=\chi(s)M_d(f)(s)\ \text{p.p.}\;,\]
soit encore, en remplaçant $f$ (considérée paire) par $If$:
\[\forall  f\in L^2(0,\infty;dx) \quad
\Reel(s)=\frac12\Rightarrow M_d(\wt
f)(s)=\chi(s)M_d(f)(1-s)\ \text{p.p.}\] On identifie $\chi$
en prenant par exemple $f(x) = \exp(-\pi x^2)$. On trouve
ainsi:
\begin{equation}\label{eq:eqfct}
\forall f\qquad\pi^{-\frac s2}\Gamma(\frac
s2)M_d(\wt
f)(s)=\pi^{-\frac{1-s}2}\Gamma(\frac{1-s}2)M_d(f)(1-s)\;.
\end{equation}
Donc:

\medskip

\noindent\framebox[\textwidth][c]{%
  \parbox[c]{0.9\textwidth}{ \emph{On a l'égalité $f = \wt
      f$ si et seulement si $\int_0^\infty f(x)x^{-s}\,dx$
      vérifie l'équation
      fonctionnelle de $\zeta(s)$, \\
      \centerline{ou encore:} On a l'égalité $f = \wt f$ si
      et seulement si $\int_0^\infty f(x)x^{s-1}\,dx$
      vérifie l'équation fonctionnelle de $\zeta(1-s)$.}}}

\medskip

Ces constatations, certes élémentaires, n'en sont pas moins
utiles pour dissiper une impression qui est laissée par le
traitement devenu habituel de l'équation fonctionnelle de
dzêta, où sont combinés le caractère autoréciproque de
$x\mapsto \exp(-\pi x^2)$, la formule de Poisson, et la
transformation de Mellin gauche $\int_0^\infty
f(x)\,x^{s-1}\,dx$.

Supposons en effet que l'on désire construire des fonctions
$f$ auto-ré\-ci\-pro\-ques sous la transformation en
cosinus. Il suffira que $\int_0^\infty f(x)x^{-s}\,dx =
\zeta(s)\int_0^\infty k(x)x^{-s}\,dx$ avec $k$ vérifiant
$k(x) = k(1/x)/x$, qui correspond à $\int_0^\infty
k(x)x^{s-1}\,dx = \int_0^\infty k(x)x^{-s}\,dx$, puisque
cela garantit que $\int_0^\infty f(x)x^{-s}\,dx$ satisfasse
à la même équation fonctionnelle que dzêta. Par l'équation
de co-Müntz \eqref{eq:comuntz} cela signifie que $f$ est la
co-somme de Poisson associée à $k$:
\[ f(x) = \sum_{n\geq1} \frac{k(x/n)}n -
\int_0^\infty \frac{k(1/t)}t\,dt \qquad k(x) =
\frac{k(1/x)}x\Rightarrow f=\wt f\]
Plus généralement on voit donc que l'équation fonctionnelle
de la fonction dzêta est, via \eqref{eq:eqfct},
essentiellement équivalente à la formule de co-Poisson.

Comme il est bien connu elle est aussi équivalente à la
formule de Poisson. Pour le voir dans le même style,
reformulons \eqref{eq:eqfct} ainsi:
\begin{equation}\label{eq:eqfct2}
 \zeta(1-s)\int_0^\infty \wt f(x) x^{-s}\,dx
=\zeta(s)\int_0^\infty f(x) x^{s-1}\,dx\;.
\end{equation}
Par la formule de Müntz le terme de droite est la
transformée de Mellin (gauche) de $A_f$. Aussi le terme de
gauche est la transformée de Mellin (gauche) de $I(A_{\wt
  f})$. Il vient donc ($f$ paire, $x>0$; l'identité vaut
peut-être presque partout seulement):
\[\sum_{n\geq1} \frac{\wt f(n/x)}x - \int_0^\infty \wt f(t)\,dt
= \sum_{n\geq1} f(nx) - \frac{\int_0^\infty
  f(t)\,dt}{x}\;,\] c'est-à-dire, la formule de Poisson!

Dans ce qui précède nous n'avons pas cherché à rendre
précises les conditions sur $f$ qui assurent la validité de
nos manipulations, car notre but était principalement de
faire passer le message que co-Poisson est aussi intimement
lié à l'équation fonctionnelle que Poisson. Nous proposons
maintenant une étude précise des conditions de validité de
\eqref{eq:muntz}.

\begin{remarque}
  Dans tout le reste de ce chapitre, la notation $\wh f(s)$
  désignera la transformée de Mellin gauche:
  \[ \wh f(s) = \int_0^\infty f(x) x^{s-1}\, dx \]
\end{remarque}

\subsection{Distributions tempérées et formule de Müntz}

Nous allons étudier de manière précise les conditions de
validité de la formule de Müntz. Rappelons nos notations: on
a $f\in L^1(0,\infty;\,dx)$, et l'on note
\[ A_f(x) = \sum_{n\geq1} f(nx) -
\frac{\int_0^\infty f(t)\,dt}{x} \;, \]
la somme modifiée, qui par le lemme \ref{lem:g} est presque
partout absolument convergente et est intégrable contre $(1
+ \log^2 x)^{-1}\,dx$. On notera aussi parfois $\wh f(s) =
\int_0^\infty f(t)t^{s-1}\,dt$ lorsque cette intégrale
existe.

Le théorème principal exprime la validité de la formule de
Müntz au sens des distributions:

\begin{theo}
  \label{theo:muntzdistrib}
  Soit $f\in L^1(0,\infty;\,dx)$. Soit $0<\sigma<1$. On
  suppose de plus
  \[ \int_0^\infty |f(x)|x^{\sigma-1}\,dx<\infty\;. \]
  La fonction $x^\sigma A_f(x)$ est alors une distribution
  tempérée en $\log(x)\in\RR$. L'identité de Müntz:
  \[ \int_0^\infty x^\sigma A_f(x) x^{i\tau}\frac{dx}x 
  = \zeta(\sigma+i\tau)\wh f(\sigma+i\tau)\]
  vaut comme une identité de distributions tempérées en
  $\tau\in\RR$, l'intégrale étant vue comme représentant une
  transformation de Fourier au sens des distributions
  tempérées.
\end{theo}

\begin{remarque}
  \label{rem:muntzdistrib}
  La situation sur la droite $\Reel(s)=1$ est différente et
  plus délicate. Les énoncés relatifs à ce cas seront donnés
  plus loin.
\end{remarque}

\begin{remarque}
  Dans l'expression:
\[ x^\sigma A_f(x) = \sum_{n\geq1} x^\sigma f(nx) -
x^{\sigma-1}\int_0^\infty f(t)\,dt\;,\] ni la somme ni le
terme avec l'intégrale ne sont des distributions tempérées
en $\log x$ (sauf si l'intégrale s'annule). C'est leur
différence qui est tempérée comme distribution en $\log x$.
\end{remarque}

\begin{cor}
  Soit $f\in L^1(0,\infty;\,dx)$. Soit $0<\sigma<1$. Si
  \begin{eqnarray*}
    \int_0^\infty |f(x)|x^{\sigma-1}\,dx<\infty\;,\qquad \int_0^\infty
  |g(x)|x^{\sigma-1}\,dx<\infty\;,\\
    \text{et}\quad \Reel(s) =\sigma\Rightarrow\ 
    \int_0^\infty
    g(x)x^{s-1}\,dx = \zeta(s)\int_0^\infty f(x)x^{s-1}\,dx\;,
  \end{eqnarray*}
  alors
  \[ g(x) = \sum_{n\geq1} f(nx) -
  \frac{\int_0^\infty f(t)\,dt}{x} \]
  pour presque tout $x>0$.
\end{cor}

Il s'agit effectivement d'un simple corollaire du théorème
\ref{theo:muntzdistrib}, puisque les deux fonctions
$x^\sigma g(x)$ et $x^\sigma A_f(x)$, en tant que
distributions tempérées en $\log(x)$, ont la même
transformée de Fourier.

\begin{theo}\label{theo:muntz1}
  Soit $f\in L^1(0,\infty;\,dx)$. Soit $0<\sigma<1$. Si
\[\begin{matrix}
  \hfill(1)&\int_0^1 \left|\sum_{n\geq1} f(nx) -
    \frac{\int_0^\infty f(t)\,dt}{x}\right|
  x^{\sigma-1}\,dx&<\infty\hfill\\
  \text{et si\ }(2)&\hfill\int_0^1
  |f(x)|x^{\sigma-1}\,dx&<\infty\;,
\end{matrix}\]
alors la formule de Müntz \ref{eq:muntz} concerne des
intégrales absolument convergentes et est valable sur la
droite $\Reel(s)=\sigma$.
\end{theo}

\begin{remarque}\label{rem:muntz1}
  On sait déjà par le lemme \ref{lem:g} que la somme
  $\sum_{n\geq1} f(nx)$ est presque partout absolument
  convergente et que les intégrales prises de $1$ à $\infty$
  sont absolument convergentes. De plus si les hypothèses
  sont vérifiées pour $\sigma$ elles le sont aussi pour
  $1>\sigma'>\sigma$. La formule de Müntz vaut donc alors au
  moins dans la bande $\sigma\leq \Reel(s) <1$.
\end{remarque}

Ce théorème est lui aussi un corollaire immédiat de
\ref{theo:muntzdistrib}.  Cependant nous en donnerons une
démonstration directe (plus longue\dots) qui ne fait pas
appel à la transformée de Fourier-Schwartz des
distributions, pour illustrer d'autres techniques plus \og
classiques\fg.

\begin{prop}\label{prop:muntz1}
  Soit $f$ une fonction intégrable sur $]0,\infty[$
  vérifiant au choix l'une des deux conditions suivantes:
  \begin{enumerate}
  \item elle est de variation totale bornée sur
    $]0,\infty[$,
  \item ou sa transformée en cosinus $\wt f$ est intégrable.
  \end{enumerate}
  La formule de Müntz \eqref{eq:muntz} est alors valable
  pour $0<\Reel(s)<1$ et concerne des intégrales absolument
  convergentes.
\end{prop}

Prouvons-le comme corollaire de \ref{theo:muntz1}. Si $f$
est de variation totale bornée, elle est bornée, donc
certainement la condition $(2)$ de \ref{theo:muntz1} est
vérifiée. Soit $df$ la mesure des variations de $f$. Quitte
à modifier $f$ en un nombre dénombrable de points, on peut
supposer $f(x) = -\int_{t>x} df(t)$. Soit $x>0$ fixé. Pour
$n\geq1$:
\[ xf(nx) - \int_{(n-1)x}^{nx} f(t)\,dt = \int_{(n-1)x<t\leq
  nx} (t - (n-1)x)df(t)\;,\]
donc:
\[\left|f(nx) - \frac1x\int_{(n-1)x}^{nx} f(t)\,dt\right| \leq
\int_{\big](n-1)x,nx\big]} \left\{\frac
  {t}x\right\}^*|df|(t)\;.\] On a noté $\{t\}^* = t-[t]$ si
$t\notin\ZZ$, $=1$ pour $t\in\ZZ$. Ainsi pour tout $x>0$:
\[ \left| \sum_{n\geq1} f(nx) - \frac{\int_0^\infty f(t)dt}x \right|
\leq \int_{]0,\infty[} \left\{\frac tx\right\}^*|df|(t) \leq
\int_{]0,\infty[} |df|(t) \;.\]
La sommation modifiée $A_f(x)$ est donc bornée, et la
condition $(1)$ de \ref{theo:muntz1} est satisfaite.

Supposons maintenant $\wt f\in L^1$. En particulier $f$
(considérée comme une fonction paire) est essentiellement
continue sur $\RR$, et la condition $(2)$ de
\ref{theo:muntz1} est donc vérifiée. Par le théorème presque
sûr de Poisson \ref{theo:poissonpp}, on a pour presque tout
$x>0$:
\[\sum_{n\geq1} f(nx) - \frac{\int_0^\infty f(t)\,dt}{x}
= \sum_{n\geq1} \frac{\wt f(n/x)}x - \int_0^\infty \wt
f(t)\,dt \;,\]
Pour vérifier la condition $(1)$ il suffit de s'intéresser à
$B(x)=\sum_{n\geq1} \frac{|\wt f(n/x)|}x$.  Par le lemme
\ref{lem:k}, $B$ est intégrable contre la fonction
décroissante et intégrable $k(x) = x^{\sigma-1}$, $x\leq1$,
$k(x) = 0$, $x>1$. Ainsi $(1)$ vaut aussi.

Un autre corollaire du théorème \ref{theo:muntz1} est donné
par:

\begin{prop}
  Soit $0<\sigma<1$.  Soit $f$ une fonction localement de
  variation bornée telle que
  \[ \int_0^\infty |f(x)|\,dx + \int_{]0,1]} x^\sigma\,|df|(x) +
  \int_{]1,\infty[} |df|(x) < \infty \;.\]
  La formule de Müntz \eqref{eq:muntz} est alors valable
  comme identité d'intégrales absolument convergentes dans
  la bande $\sigma<\Reel(s)<1$.
\end{prop}

Quitte à modifier $f$ (donc $A_f$) en un nombre dénombrable
de points, on peut supposer $f(x) = -\int_{t>x} df(t)$.
Remarquons ensuite pour $x<1$: $\int_{]x,1]}
t^\sigma\,|df|(t)\geq x^\sigma \int_{]x,1]} |df|(x)$, et
donc $f(x) = O(x^{-\sigma})$ pour $x\to0$. Ainsi
$\int_0^\infty f(x)x^{s-1}\,dx$ converge absolument pour
$\sigma<\Reel(s)\leq 1$. L'inégalité établie pour $A_f(x)$
dans la preuve de \ref{prop:muntz1}:
\[ \left| \sum_{n\geq1} f(nx) - \frac{\int_0^\infty f(t)dt}x \right|
\leq \int_{]0,\infty[} \left\{\frac tx\right\}^*|df|(t)\;,\]
vaut aussi sous la nouvelle hypothèse (on a $\int_{]0,1]}
t|df|(t)<\infty$). Comme $\left\{\frac tx\right\}^*\leq
(\left\{\frac tx\right\}^*)^\sigma\leq (\frac tx)^\sigma$,
on en déduit:
\[ |A_f(x)| \leq x^{-\sigma}\int_{]0,1]} t^\sigma|df|(t) +
\int_{]1,\infty[}|df|(t)\]
Donc $A_f$ est aussi $O(x^{-\sigma})$ pour $x\to0$. Ainsi
$\int_0^\infty A_f(x)x^{s-1}\,dx$ converge absolument pour
$\sigma<\Reel(s)<1$, la convergence sur $]1,\infty[$ par le
lemme \ref{lem:g}. La proposition est donc prouvée, comme
corollaire du théorème \ref{theo:muntz1}.

\medskip

Nous en venons maintenant aux preuves de
\ref{theo:muntzdistrib} et de \ref{theo:muntz1}, en
commençant par ce dernier. Bien qu'il puisse être vu comme
un simple corollaire de \ref{theo:muntzdistrib}, nous en
donnerons une preuve directe, dans un style plus classique.
Nous débutons par un lemme:

\begin{lem}\label{lem:muntz}
  Soit $f\in L^1(0,\infty;\,dx)$ et soit $0<\sigma<1$. Si
  $f$ est de classe $C^1$ sur $]0,\infty[$, si $f'$ est
  intégrable sur $[1,\infty[$, et si $f'(x)$ est
  $O(x^{-1-\sigma})$ pour $x\to0$ alors la formule de Müntz
  \ref{eq:muntz} vaut dans la bande $\sigma<\Reel(s)<1$ avec
  des intégrales absolument convergentes.
\end{lem}

Tout d'abord $f(x)$ est $O(x^{-\sigma})$ pour $x\to0$, donc
$\wh f(s) = \int_0^\infty f(x) x^{s-1}\,dx$ est absolument
convergent et analytique pour $\sigma<\Reel(s)\leq 1$. On
note que $\int_0^1 t|f'(t)|\,dt$ est fini. Par
\[ xf(nx) - \int_{(n-1)x}^{nx} f(t)\,dt = \int_{(n-1)x}^{nx}
(t - (n-1)x)f'(t)\,dt\;,\]
on obtient la convergence absolue pour tout $x>0$ de
$\sum_{n\geq1} f(nx)$ et, pour
\[ A(x) = \sum_{n\geq1} f(nx) - \frac{\int_0^\infty
  f(t)\,dt}{x} = \int_0^\infty f^\prime(t) \left\{\frac
  tx\right\}\,dt\;,\]
l'inégalité:
\[|A(x)| \leq \int_0^\infty |f^\prime(t)|\left\{\frac tx\right\}\,dt
\leq \int_0^x |f^\prime(t)|\frac tx\,dt + \int_x^\infty
|f^\prime(t)|\,dt\;,\]
qui est $O(x^{-\sigma})$ pour $x\to0$. Donc $\int_0^1 A(x)
x^{s-1}\,dx$ est absolument convergent et analytique pour
$\Reel(s)>\sigma$.  Par le lemme \ref{lem:g} on a
$\int^\infty |A(x)|\frac{dx}{\log^2(x)}<\infty$ et donc
certainement $\sigma<1\Rightarrow\ \int_1^\infty
|A(x)|x^{\sigma-1}\,dx<\infty$. Ainsi $\wh A(s) =
\int_0^\infty A(x)x^{s-1}\,dx$ est absolument convergent
pour $\sigma<\Reel(s)<1$.

Soit $\epsilon>0$ et soit $f_\epsilon(x) = f(x)e^{-\epsilon
  x}$. On remarque que $f_\epsilon$ vérifie les trois mêmes
hypothèses que $f$. On a par convergence dominée:
\[ \sigma<\Reel(s)\leq 1\Rightarrow
\quad \wh f(s) = \lim_{\epsilon\to0} \wh{f_\epsilon}(s)
\;.\]
Soit $A_\epsilon(x)$ la sommation de Poisson modifiée de
$f_\epsilon$. À nouveau par convergence dominée on a
$\forall x>0\ \lim_{\epsilon\to0} A_\epsilon(x) = A(x)$. De
plus les intégrales $\int_1^\infty \sum_{n\geq1}
f_\epsilon(nx) x^{s-1}\,dx$ sont dominées par $\int_1^\infty
\sum_{n\geq1}|f(nx)|x^{\Reel(s)-1}\,dx$. Donc
\[ \sigma<\Reel(s)< 1\Rightarrow
\quad \int_1^\infty \sum_{n\geq1} f(nx) x^{s-1}\,dx=
\lim_{\epsilon\to0} \int_1^\infty \sum_{n\geq1}
f_\epsilon(nx) x^{s-1}\,dx\;.\]
En ce qui concerne l'intervalle $]0,1]$ on a les expressions
\[A_\epsilon(x) = \int_0^\infty f^\prime_\epsilon(t)
\left\{\frac tx\right\}\,dt\;,\]
puis les majorations ($0<\epsilon<1$)
\[ |A_\epsilon(x)| \leq \int_0^\infty |f^\prime(t)| \left\{\frac
  tx\right\}\,dt + \int_0^\infty |f(t)|\,dt \;, \]
qui sont indépendantes de $\epsilon$ et intégrables contre
$|x^{s-1}|\,dx$ sur $]0,1]$ pour $s$ dans la bande
considérée. Nous pouvons donc affirmer:
\[ \sigma<\Reel(s)<1\Rightarrow
\quad \wh A(s) = \lim_{\epsilon\to0} \wh{A_\epsilon}(s)
\;.\]

Il suffira donc d'établir la formule de Müntz pour
$f_\epsilon$, $\epsilon>0$ fixé. Écrivons pour
$\sigma<\Reel(s)<1$:
\[ \int_0^\infty A_\epsilon(x) x^{s-1}\,dx  = \int_0^1 A_\epsilon(x)
x^{s-1}\,dx + \int_1^\infty \sum_{n=1}^\infty f_\epsilon(nx)
x^{s-1}\,dx + \frac{\int_0^\infty f_\epsilon(y)dy}{s-1} \]
La première intégrale est analytique pour $\Reel(s)>\sigma$,
et la deuxième intégrale est une fonction entière, puisque
\[ x\geq1 \Rightarrow\quad\left|\sum_{n=1}^\infty
  f_\epsilon(nx)\right| = O(1/(\exp(\epsilon x) - 1)) \;.\]
Ainsi $\wh{A_\epsilon}(s)$ est méromorphe dans
$\Reel(s)>\sigma$, et on peut la réécrire pour $\Reel(s)>1$
selon:
\[\wh{A_\epsilon}(s) = \int_0^1 \Big( A_\epsilon(x) +  \frac{\int_0^\infty
  f_\epsilon(y)dy}x \Big)x^{s-1}\,dx + \int_1^\infty
\sum_{n=1}^\infty f_\epsilon(nx)x^{s-1}\,dx \]
\[=  \int_0^\infty \sum_{n=1}^\infty
f_\epsilon(nx) x^{s-1}\,dx = \zeta(s)\int_0^\infty
f_\epsilon(x) x^{s-1}\,dx\]
La preuve du lemme est complète, puisque le principe du
prolongement analytique donne $\wh{A_\epsilon}(s) =
\zeta(s)\wh{f_\epsilon}(s)$ aussi dans la bande
$\sigma<\Reel(s)<1$.

\medskip

Prouvons maintenant le théorème \ref{theo:muntz1}.
Rappelons-en les hypothèses: on a $f\in L^1(0,\infty;\,dx)$,
on note $A(x)$ la somme de Poisson (Müntz)-modifiée de $f$
qui est définie presque partout et l'on suppose que $f$ et
$A$ sont toutes deux dans $L^1(0,\infty;\,x^{\sigma-1}dx)$
pour un certain $0<\sigma<1$. Il s'agit de prouver la
formule de Müntz. Comme nous l'avons déjà indiqué dans la
remarque \ref{rem:muntz1}, les fonctions $f$ et $A$ seront
aussi dans $L^1(0,\infty;\,x^{\sigma'-1}dx)$ pour tout
$1>\sigma'>\sigma$, et les fonctions $\wh f$ et $\wh A$ sont
continues dans la bande $\sigma\leq \Reel(s)<1$. Il suffit
donc d'établir la validité de la formule de Müntz sur la
bande $\sigma<\Reel(s)<1$.

Soit $g(t)$ une fonction quelconque non identiquement nulle
et de classe $C^1$ et de support dans $[a,b]$, $0<a<b$. Soit
$k$ quelconque localement intégrable sur $]0,\infty[$. On
note alors
\[ (g*k)(t) = \int_0^\infty k(u)g(\frac
tu)\frac{du}u = \int_{t/b}^{t/a} k(u)g(\frac tu)\frac{du}u
\;,\] la convolution multiplicative de $g$ et de $k$. Pour
tout nombre complexe $s$ tel que $\int_0^\infty
k(t)t^{s-1}\,dt$ soit absolument convergent, le théorème de
Fubini prouve que $g*k$ a la même propriété et que
\[ \int_0^\infty (g*k)(t)t^{s-1}\,dt = \int_0^\infty
g(t)t^{s-1}\,dt \int_0^\infty k(t) t^{s-1}\,dt\;.\]
En particulier comme $f$ est intégrable, il en est de même
de $F = g*f$. De plus de l'expression $F(t) = \int_0^\infty
f(u)g(\frac tu)\frac{du}u$ il ressort que $F$ est de classe
$C^1$ sur $]0,\infty[$, avec $F'(t) = \frac1t \int_0^\infty
f(u)\frac tu g^\prime(\frac tu)\frac{du}u$. Cela montre en
particulier que $t\,F'(t)$ est dans $L^1(0,\infty;\,dx)$. Et
comme $g$ est $C^1$ avec son support compact éloigné de
l'origine, on a certainement $\forall x>0\ |xg'(x)|\leq
Cx^{-\sigma}$ pour une certaine constante $C$. Il en résulte
que $F'$ est $O(t^{-\sigma-1})$ pour $t\to0$. Donc $F$
vérifie toutes les hypothèses du lemme \ref{lem:muntz} et la
formule de Müntz $\wh A_F(s) = \zeta(s) \wh F(s)$ vaut pour
$F$ dans la bande $\sigma<\Reel(s)<1$. Nous savons déjà $\wh
F(s) = \wh g(s) \wh f(s)$ et il suffit donc de montrer $\wh
A_F(s) = \wh g(s) \wh A_f(s)$, ce qui résultera de $A_F =
g*A_f$. Nous avons pour tout $x>0$ fixé:
\[\sum_{n\geq1} \int_0^\infty |f(u)||g(\frac{nx}{u})|\frac{du}u =
\int_{x/b}^{x/a} \sum_{n\geq1}
|f(nv)||g(\frac{x}{v})|\frac{dv}v < \infty\;,\]
où nous avons utilisé l'intégrabilité locale de
$\sum_{n\geq1} |f(nv)|$. Ainsi
\[ \sum_{n\geq1}  (g*f)(nx) =
\int_{x/b}^{x/a} \sum_{n\geq1} f(nv)g(\frac{x}{v})\frac{dv}v
= g*\Big(\sum_{n\geq1} f(n\,\cdot)\Big)(x)\;, \]
On vérifie aussi l'identité $(g*\frac1t)(x) =
\frac1x\int_0^\infty g(t)dt$. Il en résulte $A_F = g*A_f$ et
ceci termine la preuve du théorème \ref{theo:muntz1}.

\medskip

Venons-en au théorème \ref{theo:muntzdistrib}. On y fait les
hypothèses $0<\sigma<1$, $f\in L^1(0,\infty;\,dx)$,
$\int_0^\infty |f(x)|x^{\sigma-1}\,dx<\infty$, et l'on note
  \[ A_f(x) = \sum_{n\geq1} f(nx) -
  \frac{\int_0^\infty f(t)\,dt}{x} \]
  qui est presque partout définie et est comme fonction de
  $\log(x)\in\RR$ localement intégrable.
  
  Soit $\epsilon>0$ fixé, et soit $f_\epsilon(x) =
  f(x)\exp(-\epsilon x)$, $A_\epsilon = A_{f_\epsilon}$.
  Pour toute fonction $\phi(x)$ de classe $C^\infty$,
  supportée dans $0<a<x<b$, on a par convergence dominée:
\[ \lim_{\epsilon\to0} \int_0^\infty
A_\epsilon(t)\phi(t)\,dt = \lim_{\epsilon\to0} \int_a^b
A_\epsilon(t)\phi(t)\,dt = \int_0^\infty
A_f(t)\phi(t)\,dt\;.\]
La fonction entière $\alpha(s) = \int_a^b \phi(t)t^{-s}\,dt$
est de décroissance rapide dans toute bande $\sigma_1\leq
\Reel(s)\leq\sigma_2$. Par la formule d'inversion on a:
\[ \forall t>0\ \forall c\in\RR \quad
\phi(t) = \int_{c-i\infty}^{c+i\infty}
\alpha(s)t^{s-1}\frac{|ds|}{2\pi}\]
Prenons tout d'abord $c>1$. Comme $\int_0^\infty
|f_\epsilon(t)|t^{c-1}\,dt<\infty$ on a
\[ c>1\Rightarrow\quad \int_0^\infty \left|\sum_{n\geq1}
  f_\epsilon(nt)\right|t^{c-1}\,dt \leq \zeta(c)
\int_0^\infty |f_\epsilon(t)|t^{c-1}\,dt < \infty\;.\]
Ainsi, pour $c>1$ et $0<\Lambda<\infty$ on a
\[ \int_0^\Lambda |A_\epsilon(t)|t^{c-1}\,dt \leq 
\zeta(c) \int_0^\infty |f_\epsilon(t)|t^{c-1}\,dt +
\int_0^\infty|f_\epsilon(t)|\,dt\;\frac{\Lambda^{c-1}}{c-1}
<\infty\;.\]
On peut donc écrire, avec $c=2$ et pour tout $\Lambda>b$:
\[ \int_0^\infty A_\epsilon(t)\phi(t)\,dt = \int_0^\Lambda
A_\epsilon(t)\phi(t)\,dt = \int_{2-i\infty}^{2+i\infty}
\alpha(s)\left(\int_0^\Lambda
  A_\epsilon(t)t^{s-1}\,dt\right)\frac{|ds|}{2\pi}\]
\[ =  \int_{2-i\infty}^{2+i\infty}
\alpha(s)\left(\int_0^\Lambda \sum_{n\geq1} f_\epsilon(nt)
  t^{s-1}\,dt\right)\frac{|ds|}{2\pi} - \int_0^\infty
f_\epsilon(t)\,dt\int_{2-i\infty}^{2+i\infty}
\alpha(s)\frac{\Lambda^{s-1}}{s-1}\frac{|ds|}{2\pi} \]
On a de plus:
\[\int_{2-i\infty}^{2+i\infty} 
\alpha(s)\frac{\Lambda^{s-1}}{s-1}\frac{|ds|}{2\pi} =
\int_0^\Lambda \int_{2-i\infty}^{2+i\infty}
\alpha(s)t^{s-2}\frac{|ds|}{2\pi}\,dt = \int_0^\Lambda
\frac{\phi(t)}t\,dt\;.\]
On fait tendre $\Lambda$ vers $+\infty$, et on obtient:
\[
\int_0^\infty A_\epsilon(t)\phi(t)\,dt =
\int_{2-i\infty}^{2+i\infty} \alpha(s)\zeta(s)\int_0^\infty
f_\epsilon(t) t^{s-1}\,dt\; \frac{|ds|}{2\pi} - \alpha(1)
\int_0^\infty f_\epsilon(t)\,dt\;.\]
La décroissance rapide de $\alpha(s)$ pour
$|\Imag(s)|\to\infty$ permet de décaler l'intégrale complexe
de la droite $\Reel(s)=2$ à la droite $\Reel(s)=\sigma$, ce
qui donne un résidu au passage en $s=1$ qui vaut
$\alpha(1)\int_0^\infty f_\epsilon(t)\,dt$, et qui compense
exactement le dernier terme. Ainsi:
\[ \int_0^\infty A_\epsilon(t)\phi(t)\,dt = \int_{s=\sigma+i\tau}
\alpha(s)\zeta(s)\left(\int_0^\infty f_\epsilon(t)
  t^{s-1}\,dt\right) \frac{d\tau}{2\pi}\;.\] On a sur la
droite $\Reel(s) = \sigma$ convergence pour $\epsilon\to0$
de $\int_0^\infty f_\epsilon(t) t^{s-1}\,dt$ vers
$\int_0^\infty f(t)t^{s-1}\,dt$, avec $\int_0^\infty
|f(t)|t^{\sigma-1}\,dt$ donnant une borne supérieure fixe.
On peut donc affirmer
\[ \int_0^\infty
A_f(t)\phi(t)\,dt = \int_{s=\sigma+i\tau}
\left(\int_0^\infty \phi(t)t^{-s}\,dt\right)\zeta(s)\wh f(s)
\frac{d\tau}{2\pi}\;.\]
Posons maintenant $\phi(t) = t^{\sigma-1}\psi(t)$, et
réécrivons l'équation précédente:
\[ \int_0^\infty
t^\sigma A_f(t)\psi(t)\frac{dt}t = \int_\RR
\left(\int_0^\infty
  \psi(t)t^{-i\tau}\frac{dt}t\right)\zeta(\sigma+i\tau)\wh
f(\sigma+i\tau) \frac{d\tau}{2\pi}\;.\]
Cela exprime exactement que la distribution tempérée
$\zeta(\sigma+i\tau)\wh f(\sigma+i\tau)=D(\tau)$ a comme
transformation de Fourier inverse au sens des distributions
tempérées $\int_\RR e^{-i\tau u}D(\tau)\frac{d\tau}{2\pi}$
la fonction localement intégrable $e^{\sigma u} A_f(e^u)$.
Donc $t^\sigma A_f(t)$ est une distribution tempérée en $u =
\log(t)\in\RR$ et l'identité de Müntz:
\[ \int_0^\infty t^\sigma A_f(t)\,t^{i\tau}\frac{dt}t = 
\zeta(\sigma+i\tau)\wh f(\sigma+i\tau)\;,\]
vaut au sens des distributions tempérées. Ainsi le théorème
\ref{theo:muntzdistrib} est établi.

\medskip

La situation sur la droite $\Reel(s)=1$ est un peu
différente.

\begin{theo}
   \label{theo:muntzdistrib1a}
   Soit $f\in L^1(0,\infty;\,dx)$.
   \begin{enumerate}
   \item La fonction $D(u) = \sum_{n\geq1} x f(nx)$,
     $x=\exp(u)$ est définie ponctuellement pour presque
     tout $u\in\RR$, est localement intégrable, et est
     tempérée comme distribution.
   \item On a les identités de distributions tempérées:
   \begin{align*}
     \int_\RR D'(u) e^{i\tau u}\,du
     &= -i\tau \zeta(1+i\tau)\wh f(1+i\tau)\\
     \text{sur }\RR\setminus\{0\}\quad \int_\RR D(u)
     e^{i\tau u}\,du&= \zeta(1+i\tau)\wh f(1+i\tau)
    \end{align*}
  \item Si $f\in L^1(0,\infty;\,(1+|\log x|)dx)$, alors:
   \[ \int_0^\infty \left(\sum_{n\geq1} xf(nx) -
     \frac12\int_0^\infty f(t)\,dt\right)
   x^{i\tau}\;\frac{dx}x = \mbox{\rm
     V.P.}\;\zeta(1+i\tau)\wh f(1+i\tau)\]
   vaut comme une identité de distributions tempérées en
   $\tau$, le symbole $\mbox{\rm V.P.}$ représentant la
   valeur principale au sens de Cauchy.
   \end{enumerate}
\end{theo}

\begin{remarque}
  On notera le $\frac12$ dans
  \ref{theo:muntzdistrib1a}$.(3)$.
\end{remarque}

Sous la seule hypothèse $f\in L^1(0,\infty;\,dx)$, le lemme
\ref{lem:g} permet d'affirmer que $\sum_{n\geq1} f(nx)$ est
intégrable sur $]0,\infty[$ contre $dx/(1+\log^2 x)$. Ainsi
$\sum_{n\geq1} xf(nx)$ comme fonction de $u=\log(x)$ est
intégrable sur $\RR$ contre $du/(1+u^2)$.  Donc
$\sum_{n\geq1} xf(nx)$ est une distribution tempérée en
$u=\log(x)$. Ceci prouve $(1)$.

Pour la preuve du point $(2)$ nous reprenons la technique de
la preuve de \ref{theo:muntzdistrib}.  Soit $\epsilon>0$ et
$f_\epsilon(x) = e^{-\epsilon x}f(x)$; soit $\phi(t)$ de
classe $C^\infty$ à support dans $[a,b]$, $0<a<b$ et soit
$\alpha(s) = \int_a^b \phi(t)t^{-s}\,dt$. On justifie comme
précédemment:
\[
\int_0^\infty \sum_{n\geq1} f_\epsilon(nt)\phi(t)\,dt =
\int_{2-i\infty}^{2+i\infty} \alpha(s)\zeta(s)\int_0^\infty
f_\epsilon(t) t^{s-1}\,dt\; \frac{|ds|}{2\pi} \;.\]
Faisons l'hypothèse $\alpha(1) = \int_0^\infty \phi(t)
\frac{dt}t = 0$.  Pour $\epsilon>0$ fixé, la fonction
$\int_0^\infty f_\epsilon(t) t^{s-1}\,dt$ est analytique
pour $\Reel(s)>1$, continue pour $\Reel(s)\geq1$, bornée
pour $1\leq \Reel(s)\leq2$. On peut donc décaler l'intégrale
complexe sur la droite $\Reel(s) = 1$. Sur cette droite les
quantités $\int_0^\infty f_\epsilon(t) t^{s-1}\,dt$ sont
bornées indépendamment de $\epsilon$ et de $s$. On peut donc
faire tendre $\epsilon$ vers $0$, et ainsi, ($\alpha(1) =
0$):
\[
\int_0^\infty \sum_{n\geq1} tf(nt)\phi(t)\frac{dt}t =
\int_\RR \alpha(1+i\tau)\zeta(1+i\tau)\wh
f(1+i\tau)\frac{d\tau}{2\pi}\;.
\]
Choisissons $\phi(t)$ de la forme $t\frac{d}{dt}\psi(t)$
avec $\psi\in C^\infty([a,b])$. On a alors
\[ \alpha(s) = \int_0^\infty t\psi'(t) t^{-s}\,dt
= (s-1)\int_0^\infty \psi(t) t^{-s}\,dt = (s-1)\beta(s)\;.\]
On obtient donc pour toute fonction $\psi$ de $t>0$, $C^\infty$,
supportée dans $[a,b]$:
\[
\int_0^\infty (\sum_{n\geq1} tf(nt))t\psi'(t)\frac{dt}t =
\int_\RR i\tau\beta(1+i\tau)\zeta(1+i\tau)\wh
f(1+i\tau)\frac{d\tau}{2\pi}\;.
\]
En passant à la variable $u=\log(t)$ et en utilisant la
distribution tempérée $D(u) = \sum_{n\geq1} e^u f(ne^u)$
cela s'écrit
\[
- \int_\RR D'(u)\theta(u)\,du = \int_\RR i\tau\left(\int_\RR
  \theta(u)e^{-i\tau u}\,du\right) \zeta(1+i\tau)\wh
f(1+i\tau)\frac{d\tau}{2\pi}\;,\]
pour toute fonction $\theta(u)$ de classe $C^\infty$ et à
support compact. Cela signifie exactement que la
distribution tempérée
\[ E(u) = \int_\RR i\tau 
\zeta(1+i\tau)\wh f(1+i\tau)e^{-i\tau
  u}\frac{d\tau}{2\pi}\;,\]
vérifie $E(u) = -D'(u)$. Donc $\int_\RR e^{i\tau u}E(u)\,du$
comme distribution tempérée en $\tau$ est d'une part
$+i\tau\int_\RR e^{i\tau u}D(u)\,du$, d'autre part $i\tau
\zeta(1+i\tau)\wh f(1+i\tau)$. Ceci prouve que la
distribution $\int_\RR e^{i\tau u}D(u)\,du$ a sa restriction
à l'ouvert $\{\tau\neq0\}$ égale à la fonction
$\zeta(1+i\tau)\wh f(1+i\tau)$.

Pour la preuve de $(3)$ on introduit la notation:
\[ A_f^*(x) = \sum_{n\geq1} f(nx) -
\frac12\frac{\int_0^\infty f(t)\,dt}x\;,\]
puis on reprend à l'identique la démonstration de
\ref{theo:muntzdistrib}. On obtient:
\[ \int_0^\infty A_\epsilon^*(t)\phi(t)\,dt = 
\int_{2-i\infty}^{2+i\infty} \alpha(s)\zeta(s)\int_0^\infty
f_\epsilon(t) t^{s-1}\,dt\; \frac{ds}{2\pi i} -
\frac12\alpha(1) \int_0^\infty f_\epsilon(t)\,dt\;.\]
On déforme le contour d'intégration vers le contour allant
verticalement de $1-i\infty$ à $1-i\delta$ puis le long d'un
demi-cercle dans $\Reel(s)\geq1$ de $1-i\delta$ à
$1+i\delta$, puis verticalement de $1+i\delta$ à
$1+i\infty$. On note que grâce à l'hypothèse $f\in
L^1(0,\infty;(1+|\log x|)dx)$ on peut écrire
\[\alpha(s)\zeta(s)\wh{f_\epsilon}(s) = \frac{\alpha(1)
  \wh{f_\epsilon}(1)}{s-1} + O_\epsilon(1)\]
pour $1\leq \Reel(s)$ , $|s-1|\leq1$. La contribution le
long du demi-cercle de rayon $\delta$ est donc
$\frac12\alpha(1)\wh{f_\epsilon}(1) + O_\epsilon(\delta)$.
Ainsi:
\[ \int_0^\infty A_\epsilon^*(t)\phi(t)\,dt = 
\lim_{\delta\to0}\int_{\RR\setminus[-\delta,+\delta]}
\alpha(1+i\tau)\zeta(1+i\tau)\wh{f_\epsilon}(1+i\tau)\;
\frac{d\tau}{2\pi} \;.\]
On notera que, compte tenu de $f\in L^1(0,\infty;(1+|\log
x|)dx)$, l'on a pour $\tau\in\RR$
\[ \wh{f_\epsilon}(1+i\tau) = \wh{f_\epsilon}(1) + O(\tau)
\;,\]
avec une constante implicite ne dépendant pas de
$\epsilon>0$. Cela signifie que la convergence pour
$\delta\to0$ est uniforme par rapport à $\epsilon$. On peut
donc prendre la limite pour $\epsilon\to0$ à l'intérieur de
l'intégrale. Ainsi:
\[ \int_0^\infty A_f^*(t)\phi(t)\,dt = 
\lim_{\delta\to0}\int_{\RR\setminus[-\delta,+\delta]}
\alpha(1+i\tau)\zeta(1+i\tau)\wh{f}(1+i\tau)\;
\frac{d\tau}{2\pi} \;.\]
Ceci prouve que la valeur principale au sens de Cauchy de la
fonction $\zeta(1+i\tau)\wh{f}(1+i\tau)$ existe et est bien
une distribution, en fait une distribution tempérée égale au
sens des distributions à $\int_0^\infty tA_f^*(t)\,
t^{i\tau}\;\frac{dt}t$.  Ceci complète la preuve du théorème
\ref{theo:muntzdistrib1a}.

\begin{cor}
  Au sens des distributions on a
   \[ \text{pour }\tau\neq0\quad \lim_{N\to\infty}
   \sum_{n=1}^N \frac1{n^{1+i\tau}} = \zeta(1+i\tau)\]
\end{cor}

Pour toute fonction de Schwartz $\beta(u)$ on a $\beta(u) =
O(1/(1+u^2))$ et donc:
 \[ \int_\RR  \sum_{1\leq
   n} e^u|f(ne^u)| \beta(u)\,du <\infty \;,\] ce qui
 implique $\lim_{N\to\infty} \sum_{1\leq n\leq N} xf(nx) =
 \sum_{1\leq n} xf(nx)$, au sens des distributions tempérées
 en $\log(x)$.
 La transformée de Fourier d'une limite au sens des
 distributions est une limite au sens des distributions. Il
 suffit alors de choisir $f\in L^1(0,\infty;dx)$ de manière
 à ce que $\wh f$ soit égale à $1$ sur un intervalle
 $[\tau_1,\tau_2]$ ne rencontrant pas $0$, et de restreindre
 ensuite à l'intervalle ouvert $]\tau_1,\tau_2[$, pour
 obtenir la conclusion.
 
 Une démonstration directe et plus simple s'obtient à partir
 de l'expression
\[ \zeta(s) = 1 + \frac1{2^s} + \dots
+\frac1{N^s} +\frac{N^{1-s}}{s-1} - s\int_N^\infty
\frac{\{t\}}{t^{s+1}}\,dt
\]
qui est valable pour $N\geq1$, $\Reel(s)>0$, $s\neq1$. En
effet, au sens des distributions en $\tau$:
$\lim_{N\to\infty} N^{-i\tau} = 0$.
     
\subsection{La transformation de Fourier de la fonction dzêta} 

Revenons à la formule
   \[ \int_0^\infty \left(\sum_{n\geq1} xf(nx) -
     \frac12\int_0^\infty f(t)\,dt\right)
   x^{i\tau}\;\frac{dx}x = \mbox{\rm
     V.P.}\;\zeta(1+i\tau)\wh f(1+i\tau)\] pour $f\in
   L^1(0,\infty;\,(1+|\log x|)dx)$. On peut aussi écrire
   sous cette hypothèse:
\[ \mbox{\rm V.P.}\;\zeta(1+i\tau)\wh f(1+i\tau) = \wh f(1) \mbox{\rm
  V.P.}\;\zeta(1+i\tau) + \zeta(1+i\tau)(\wh f(1+i\tau) -
\wh f(1))\]
Choisissons $f(x)$ de la forme $\theta(u)e^{-u}$, $u=\log
x$, avec $\theta$ de classe $C^\infty$ à support compact.
Ainsi:
\[ \sum_{n\geq1} xf(nx) -
\frac12\int_0^\infty f(t)\,dt = \sum_{n\geq1} \frac1n
\theta(\log n+u) - \frac12 \int_\RR \theta(v)\,dv\]
On a par ailleurs
\[ \wh f(1+i\tau) = \int_0^\infty f(x)x^{i\tau}\,dx =
\int_\RR \theta(u) e^{i\tau u} \,du \]
Prenons $\theta$ avec $\int_\RR \theta(u)\,du = 1$ et
remplaçons $\theta(u)$ par $\theta_A(u) = A\theta(A u)$,
avec $A>0$, $A\to\infty$. On obtient l'identité de
distributions tempérées:
\[\begin{split}\int_0^\infty \left(
    \sum_{n\geq1} \frac1n \theta_A(\log n+u) - \frac12
  \right)e^{i\tau u}\,du = \\
  \zeta(1+i\tau)(\wh f(1+i\frac\tau A) - 1) + \mbox{\rm
    V.P.}\;\zeta(1+i\tau)\end{split} \]
On passe à la limite au sens des distributions pour
$A\to\infty$, et on obtient:
\[\int_0^\infty \left(
  \sum_{n\geq1} \frac1n \delta(\log n+u) - \frac12
\right)e^{i\tau u}\,du = \mbox{\rm V.P.}\;\zeta(1+i\tau) \]
Nous avons démontré:

\begin{prop}
  La transformée de Fourier au sens des distributions de
  $\mbox{\rm V.P.}\;\zeta(1+i\tau)$ est donnée par la
  formule:
\[ \int_\RR \left(\mbox{\rm
    V.P.}\;\zeta(1+i\tau)\right) e^{-i\tau
  u}\,\frac{d\tau}{2\pi} = \sum_{n\geq1} \frac1n \delta(\log
n+u) - \frac12 \]
\end{prop}

On peut aussi le voir de manière plus directe en partant de
l'égalité:
\[ \zeta(s) = \frac{s}{s-1} - s\int_1^\infty
\frac{\{t\}}{t^{s+1}}\,dt\;.
\]
Avec $s=1+i\tau$, $t=e^{-u}$, et en menant les calculs au
sens des distributions:
\[\zeta(1+i\tau) = 1 + \frac{1}{i\tau} -
(1+i\tau)\int_{-\infty}^0 \{e^{-u}\} e^{1+i\tau u}\,du\]
\[ = \frac{1}{i\tau} + 1 + \int_{]-\infty,0[} e^{(1+i\tau) u}
\frac{d}{du}\{e^{-u}\}\,du \;.\] On remarque:
\[ \frac{d}{du}\{e^{-u}\} =  \sum_{n\geq1}
\delta(u+\log n)-e^{-u} \;.\]
Donc
\[ \mbox{\rm V.P.}\; \zeta(1+i\tau) =  \mbox{\rm V.P.}\;
\frac{1}{i\tau} + 1 + \int_{]-\infty,0[} \left(\sum_{n\geq1}
  \frac1n \delta(u+\log n) - 1\right) e^{i\tau u}\,du\;,\]
\[ \mbox{\rm V.P.}\; \zeta(1+i\tau) =  \mbox{\rm V.P.}\;
\frac{1}{i\tau} + \int_\RR \left(\sum_{n\geq1} \frac1n
  \delta(u+\log n) - \Un_{u<0}(u)\right) e^{i\tau u}\,du
\;.\]
En utilisant alors la formule connue:
\[ \int_\RR \left(\mbox{\rm V.P.}\;\frac1{i\tau}\right)
e^{-i\tau u}\frac{d\tau}{2\pi} = \frac12\Un_{u<0}(u) -
\frac12 \Un_{u>0}(u)\;,\]
sous la forme:
\[ \mbox{\rm V.P.}\;\frac1{i\tau} = \int_\RR \left(\frac12\Un_{u<0}(u) -
  \frac12 \Un_{u>0}(u)\right) e^{i\tau u}\,du\;,\]
on obtient finalement
\[ \mbox{\rm V.P.}\; \zeta(1+i\tau) = \int_\RR
\left(\sum_{n\geq1} \frac1n \delta(u+\log n) -\frac12
\right) e^{i\tau u}\,du \;.\]

Dans le cas $\sigma>1$ on a bien évidemment:
\[ \sigma>1\Rightarrow\quad \int_\RR \zeta(\sigma+i\tau)e^{-i\tau
  u}\frac{d\tau}{2\pi} = \sum_{n\geq1} \frac1{n^\sigma}
\delta(u+\log n)\;.\]
Et pour $\sigma<1$ on a:

\begin{prop}\label{prop:fourierzeta}
  Soit $-\infty<\sigma<1$. La distribution tempérée
  $\zeta(\sigma+i\tau)$ a une transformée de Fourier qui est
  donnée par la formule:
\[ \int_\RR \zeta(\sigma+i\tau)e^{-i\tau
  u}\frac{d\tau}{2\pi} = \sum_{n\geq1} \frac1{n^\sigma}
\delta(u+\log n) - e^{(\sigma-1)u}\;.
\]
\end{prop}

On remarquera que la somme est convergente au sens des
distributions, mais pas au sens des distributions tempérées;
ce n'est qu'après avoir soustrait l'exponentielle que l'on
retrouve une distribution tempérée. Ce sont des exercices
intéressants de montrer directement que cette différence est
bien une distribution tempérée pour tout $\sigma<1$, ou de
montrer la proposition par des calculs analogues à ceux de
la preuve précédente. Il est aussi possible d'établir la
proposition \ref{prop:fourierzeta} de la manière suivante:
soit $f\in C^\infty([a,b])$, $0<a<b$. Soit (pour $x>0$)
$A_f(x) = \sum_{n\geq1} f(nx) - (\int_0^\infty f(t)\,dt)/x$.
Alors $A_f(x)$ est $O(x^N)$ lorsque $x\to0^+$ pour tout
$N\in\NN$ (en fait $A_f(x) = \phi(1/x)/x$ avec
$\phi\in\cS$). Donc $\wh{A_f}(s) = \int_0^\infty A_f(x)
x^{s-1}\,dx$ est absolument convergent dès que $\Reel(s)<1$,
et définit donc une fonction analytique dans ce demi-plan.
Par la formule de Müntz \eqref{eq:muntz}, on a $\wh{A_f}(s)
= \zeta(s)\wh f(s)$ dans la bande critique donc dans tout ce
demi-plan.  Cette identité, restreinte à la droite
$\Reel(s)=\sigma$ détermine la distribution (tempérée)
$\int_\RR \zeta(\sigma+i\tau)e^{-i\tau u}\frac{d\tau}{2\pi}
$ comme étant donnée par la formule de la proposition
\ref{prop:fourierzeta}.

On peut donc considérer pour $0<\sigma<1$ que la formule de
la proposition \ref{prop:fourierzeta} est structurellement
équivalente à la formule de Müntz \eqref{eq:muntz}, à
condition toutefois que $f$ soit de support compact en $\log
x$, ou d'autres conditions s'en rapprochant.

\subsection{Fonctions de carrés intégrables}

Jusqu'à présent, nous avons pris $f$ dans
$L^1(0,\infty;\,dx)$. Supposons que l'hypothèse soit $f\in
L^2(0,\infty;\,dx)$. On ne peut plus alors définir de
fonction $A_f$ par la formule \ref{eq:A}.

Pour $f$ de carré intégrable la transformée de Mellin
\[ \wh f(s) = \int_0^\infty f(x) x^{s-1}\,dx \]
n'existe, en général, que sur la droite critique $\Reel(s) =
\frac12$, et seulement au sens $L^2$ (comme limite $L^2$ des
intégrales prises sur $[a,A]$, $a\to0$, $A\to\infty$). On
utilisera parfois l'appellation \og transformation de
Mellin-Plancherel\fg\ dans ce contexte. La transformation
inverse, au sens $L^2$, est
\[ f(x) = \int_{\Reel(s)=\frac12} x^{-s} \wh
f(s)\,\frac{|ds|}{2\pi} \;.\]

\begin{theo}
  \label{theo:muntzL2}
  Soit $f\in L^2(0,\infty;dx)$, et soit pour $\Lambda>0$:
  \[ A_\Lambda(x) = \sum_{nx\leq \Lambda} f(nx) -
  \frac{\int_0^\Lambda f(t)\,dt}x \;.\]
  Les fonctions localement intégrables $\sqrt x\,
  A_\Lambda(x)$ sont des distributions tempérées en
  $\log(x)$ et convergent au sens des distributions
  tempérées en $\log(x)$ lorsque $\Lambda\to\infty$ vers une
  distribution $\sqrt xD_f(x)$ sur $]0,\infty[$. On a au
  sens des distributions tempérées en $\log x$ et en $\tau$:
  \begin{equation}\label{eq:Da}
    \int_0^\infty \sqrt x\,D_f(x)x^{i\tau}\,\frac{dx}x =
  \zeta(\frac12+i\tau)\wh f(\frac12 + i\tau)\;.
  \end{equation}
  La distribution $D_f$ sur $]0,\infty[$ est caractérisée
  par la formule suivante:
  \begin{equation}\label{eq:Db}
    \int_0^\infty D_f(x)\phi(x)\,dx = \int_0^\infty
  f(x)\left( \sum_{n\geq1} \frac{\phi(x/n)}n - \int_0^\infty
    \frac{\phi(1/t)}t\,dt\right)\,dx \;,
  \end{equation}
  pour toute fonction $\phi$ de classe $C^\infty$ à support
  compact éloigné de $0$. Elle vérifie également, au sens
  des distributions sur $]0,\infty[$:
  \begin{equation}\label{eq:Dc}
    D_f(x) = x\frac{d}{dx} 
    \int_0^\infty f(xu)\frac{\{u\}}u\,du\;.
  \end{equation}
\end{theo}

\begin{remarque}
  Dans \eqref{eq:Db}, le terme entre parenthèses est dans la
  classe de Schwartz et il s'agit donc d'une intégrale
  absolument convergente. Dans \eqref{eq:Dc} l'intégrale est
  une fonction continue de $x>0$.
\end{remarque}

Soit $f_\Lambda(x) = \Un_{0<x\leq\Lambda}(x)f(x)$. On a
$f_\Lambda\in L^1(0,\infty;\,dx)$ et $A_\Lambda =
A_{f_\Lambda}$.  Les fonctions $A_\Lambda$ sont dans
$L^2(\epsilon,\infty;dx)$ pour tout $\epsilon>0$ et aussi
dans $L^1(0,\infty;dx/(1+\log^2x))$ par le lemme
\ref{lem:g}.

Soit $\sigma$ vérifiant $\frac12<\sigma<1$. Comme
$\int_0^\infty |f_\Lambda(x)| x^{\sigma-1}\,dx < \infty$, on
peut appliquer le théorème \ref{theo:muntzdistrib}, en fait
sous la forme obtenue dans sa démonstration, à savoir
l'identité d'intégrales absolument convergentes:
\begin{equation}\label{eq:truc}
 \int_0^\infty
A_\Lambda(t)\phi(t)\,dt = \int_{s=\sigma+i\tau}
\left(\int_0^\infty \phi(t)t^{-s}\,dt\right)\zeta(s)\wh{f_\Lambda}(s)
\frac{d\tau}{2\pi}\;,
\end{equation}
pour toute fonction $\phi(t)$ de classe $C^\infty$ et à
support dans $[a,b]$, $0<a<b<\infty$. La fonction
$\wh{f_\Lambda}(s) = \int_0^\Lambda f(x) x^{s-1}dx$ est
$\Lambda^{s-\frac12}$ fois une fonction de l'espace de Hardy
du demi-plan $\Reel(s)>\frac12$. Ses restrictions aux
droites $\Reel(s) = \sigma$ convergent donc au sens $L^2$
vers $\wh{f_\Lambda}(\frac12 +i\tau)$ (elle-même n'est
définie qu'au sens $L^2$). On peut donc prendre la limite
pour $\sigma\to\frac12$ dans \eqref{eq:truc}, ce qui donne:
\[
\int_0^\infty A_\Lambda(t)\phi(t)\,dt = \int_\RR
\left(\int_0^\infty
  \phi(t)t^{-\frac12-i\tau}\,dt\right)\zeta(\frac12+i\tau)
\wh{f_\Lambda}(\frac12+i\tau) \frac{d\tau}{2\pi}\;.
\]
Écrivons $\phi(t) = \sqrt t\,\psi(t)$. Ainsi pour toute
fonction $\psi\in C^\infty$, à support dans $[a,b]$,
$0<a<b<\infty$:
\begin{multline}\label{eq:truc2}
  \int_0^\infty \sqrt t\,A_\Lambda(t)\psi(t)\,dt =\\
  \int_\RR \left(\int_0^\infty
    \psi(t)t^{-i\tau}\,dt\right)\zeta(\frac12 +
  i\tau)\wh{f_\Lambda}(\frac12+i\tau) \frac{d\tau}{2\pi}\;.
\end{multline}
Ceci prouve que la distribution tempérée
$\zeta(\frac12+i\tau)\wh{f_\Lambda}(\frac12+i\tau)$ a comme
transformée de Mellin inverse $\int_\RR
t^{-i\tau}\zeta(\frac12 +
i\tau)\wh{f_\Lambda}(\frac12+i\tau) \frac{d\tau}{2\pi}$ au
sens des distributions en $\tau$ et en $\log t$ la fonction
localement intégrable $\sqrt t\, A_\Lambda(t)$.  Ainsi,
cette fonction est une distribution tempérée en $\log t$.
Comme on a convergence au sens $L^2$ de $\wh{f_\Lambda}$
vers $\wh f$, on peut prendre la limite pour
$\Lambda\to\infty$ dans \eqref{eq:truc2}. On obtient:
\begin{multline*}
  \lim_{\Lambda\to\infty}\int_0^\infty \sqrt t\,
  A_\Lambda(t)\psi(t)\,dt =\\
  \int_\RR \left(\int_0^\infty
    \psi(t)t^{-i\tau}\,dt\right)\zeta(\frac12 +
  i\tau)\wh{f}(\frac12+i\tau) \frac{d\tau}{2\pi}\;.
\end{multline*}
Ceci identifie la transformée de Mellin inverse de
$\zeta(\frac12+i\tau)\wh{f}(\frac12+i\tau)$ comme étant la
limite au sens des distributions en $\log x$, lorsque
$\Lambda$ tend vers $\infty$, des fonctions $\sqrt x\,
A_\Lambda(x)$.

Revenons à l'intégrale absolument convergente $\int_0^\infty
A_\Lambda(t)\phi(t)\,dt$. On a:
\[ \int_0^\infty A_\Lambda(t)\phi(t)\,dt =
\int_0^\infty f_\Lambda(t)\left(\sum_{n\geq1}
  \frac{\phi(t/n)}n - \int_0^\infty
  \frac{\phi(1/u)}u\,du\right)\,dt\; ,\]
et on passe à la limite lorsque $\Lambda\to\infty$. Ceci
prouve \eqref{eq:Db}.

Posons $\phi_1(t)=\phi(1/t)/t$, de sorte que le terme entre
parenthèses dans \eqref{eq:Db} s'écrive
\[ \sum_{n\geq1} \frac{\phi_1(n/t)}t - \int_0^\infty
\phi_1(u)\,du = \int_0^\infty \frac{\{tx\}}t
\phi_1'(x)\,dx\; .\]
L'intégrale double
\[\int_0^\infty \int_0^\infty f(t)\frac{\{tx\}}t
\phi_1'(x)\,dt\,dx \]
est absolument convergente, puisque
\[ \int_0^\infty |f(t)|\frac{\{tx\}}t\,dt \leq \sqrt x\, {\|
  f\|}_2\; {\|\frac{\{t\}}t\|}_2 \; .\]
On peut donc écrire:
\[\begin{split}
  \int_0^\infty D_f(x)\phi(x)\,dx &= \int_0^\infty
  \left(\int_0^\infty f(t)\frac{\{tx\}}t\,dt\right)
  x\frac{d}{dx} \phi_1(x)\,\frac{dx}x\\
  &= - \int_0^\infty \left(\int_0^\infty
    f(t)\frac{\{t/x\}}t\,dt\right)
  x\frac{d}{dx}(x\phi(x))\,\frac{dx}{x}\\
  &= \int_0^\infty x\frac{d}{dx} \left(\int_0^\infty
    f(t)\frac{\{t/x\}}t\,dt\right) \phi(x)\,dx\;,
\end{split}
\]
la dernière ligne étant écrite au sens des distributions sur
$]0,\infty[$. On a donc en ce sens
\[ D_f(x) = x\frac{d}{dx} 
\int_0^\infty f(t)\frac{\{t/x\}}t\,dt = x\frac{d}{dx}
\int_0^\infty f(xu)\frac{\{u\}}u\,du\;.\]

La preuve du théorème \ref{theo:muntzL2} est complète.
 
\begin{remarque}
  La fonction $x\mapsto G(x) = \int_0^\infty
  f(xu)\{u\}u^{-1}\,du$ est une fonction continue de $x>0$,
  que l'on peut aussi exprimer par l'identité de Parseval,
  compte tenu de $\zeta(s)/s = -\int_0^\infty \{u\}
  u^{-s-1}\,du$ ($0<\Reel(s)<1$), sous la forme:
\[ x\mapsto -\int_{s = \frac12+i\tau} x^{-s}\wh f(s)
\frac{\overline{\zeta(1-s)}}{\;\overline{1-s}\;}\,\frac{d\tau}{2\pi}
= - x^{-\frac12}\int_ {s = \frac12+i\tau} x^{-i\tau}\wh f(s)
\frac{\zeta(s)}{s}\,\frac{d\tau}{2\pi}
\]
Cela montre que la fonction $\sqrt x\,G(x)$, comme fonction
en $\log x$ est la transformée de Fourier inverse de
$-\wh{f}(\frac12+i\tau)\zeta(\frac12+i\tau)/(\frac12+i\tau)$.
Au sens des distributions en $\log x$ et en $\tau$ on a
$i\tau = -xd/dx$, donc par \eqref{eq:Da}:
\[ \sqrt x\,D_f(x) = -(\frac12 - x\frac{d}{dx})(\sqrt x\,
G(x))\;.\]
L'on retrouve (et l'on pourrait la démontrer ainsi) la
formule:
\[ D_f(x) = x\frac{d}{dx} G(x) = x\frac{d}{dx}\int_0^\infty
f(xu)\frac{\{u\}}u\,du\;.\]
\end{remarque}

\begin{prop}
  \label{prop:muntzL2b}
  Soit $f\in L^2(0,\infty;dx)$, et soit pour $\epsilon>0$:
\[ A_\epsilon(x) = \sum_{n\geq1} f(nx)e^{-\epsilon n x} -
\frac{\int_0^\infty f(t)e^{-\epsilon t}\,dt}x \;.\]
Les fonctions localement intégrables $\sqrt x\,
A_\epsilon(x)$ sont des distributions tempérées en $\log(x)$
et convergent au sens des distributions tempérées en
$\log(x)$ lorsque $\epsilon\to0$ vers la distribution $\sqrt
x D_f(x)$ sur $]0,\infty[$.
\end{prop}

On a $A_\epsilon = A_{f_\epsilon}$ avec $f_\epsilon(x) =
f(x)e^{-\epsilon x}$. La fonction $f_\epsilon$ est dans
$L^1$ donc la somme dans $A_\epsilon$ est presque partout
absolument convergente et $A_\epsilon$ est localement (en
$\log x$) intégrable par le lemme \ref{lem:g}. De plus
$f_\epsilon$ est dans $L^1(0,\infty; x^{\sigma-1}\,dx)$ pour
tout $\sigma>\frac12$; les fonctions
$\wh{f_\epsilon}(\sigma+i\tau)$ convergent au sens $L^2$
vers $\wh{f_\epsilon}(\frac12 + i \tau)$ pour
$\sigma\to\frac12$ (puisque ce sont les transformées de
Mellin-Plancherel des fonctions
$f_\epsilon(x)x^{\sigma-\frac12}$); les fonctions
$\wh{f_\epsilon}(\frac12 + i\tau)$ convergent au sens $L^2$
vers $\wh{f}(\frac12 + i\tau)$. Compte tenu de ces éléments
la proposition est démontrée par une preuve exactement
semblable à celle de \ref{theo:muntzL2}.

\begin{prop}
  \label{prop:poissonL2}
  Soit $f\in L^2(0,\infty;\,dx)$ et soit
  \[ D_f(x) =
  x\frac{d}{dx}\int_0^\infty f(xu)\frac{\{u\}}u\,du\]
  la distribution sur $]0,\infty[$ définie dans
  \ref{theo:muntzL2}.  On a
  \[ D_f(x) = \frac1x D_{\wt f}(\frac 1x) \; , \]
  avec $\wt f$ la transformée en cosinus de $f$.
\end{prop}

Il suffira pour la preuve de rappeler par exemple la formule
\eqref{eq:eqfct2} (ici, on a Mellin gauche!):
\[
\zeta(\frac12-i\tau)\wh{\wt f\;}(\frac12 -i\tau)
=\zeta(\frac12 + i\tau)\wh{f}(\frac12 + i\tau)\;,
\]
et d'invoquer l'équation \eqref{eq:Da}. Ou encore on
utilisera \eqref{eq:Db}:
  \[ \int_0^\infty D_f(x)\phi(x)\,dx = \int_0^\infty
  f(x)\left( \sum_{n\geq1} \frac{\phi(x/n)}n - \int_0^\infty
    \frac{\phi(1/t)}t\,dt\right)\,dx \;,\]
  et la formule de co-Poisson pour $\phi$.

\begin{prop}\label{prop:fourierfract}
  La transformée en cosinus $\int_0^\infty 2\cos(2\pi
  uv)g(u)\,du$ de la fonction de $L^2(0,\infty;\,du)$
  \[ u\mapsto g(u) = \frac{\{u\}}u \]
  est la fonction de $L^2(0,\infty;\,dv)$:
  \[ v\mapsto -\frac{\{v\}}v + \int_v^\infty \{u\}\frac{du}{u^2}\;. \]
\end{prop}

Compte tenu de la formule connue $\zeta(s)/s =
-\int_0^\infty \{u\} u^{-s-1}\,du$ pour $0<\Reel(s)<1$, on a
\[ \wh g(\frac12 + i\tau) =
\frac{\zeta(\frac12-i\tau)}{\frac12-i\tau}\;.
\]
Or, par \eqref{eq:eqfct2}
\[
\zeta(\frac12-i\tau)\wh{\wt g\;}(\frac12 -i\tau)
=\zeta(\frac12 + i\tau)\wh{g}(\frac12 + i\tau)\;,
\]
donc:
\[ \wh{\wt g\;}(\frac12 +i\tau) =
\frac{\zeta(\frac12-i\tau)}{\frac12+i\tau} =
\frac{\frac12-i\tau}{\frac12+i\tau}\; \wh
g(\frac12+i\tau)\;.
\]
On vérifie aisément que l'opérateur de Hardy $f(x)\mapsto
\frac1x \int_0^x f(t)\,dt$ correspond à la multiplication
par $1/(1-s)$ sur la droite critique (pour la transformée de
Mellin $\wh f(s) = \int_0^\infty f(x)x^{s-1}\,dx$), et donc
que l'opérateur $f(x)\mapsto \int_x^\infty \frac{f(t)}t\,dt$
correspond à la multiplication par $1/s$. Comme
$(\frac12-i\tau)/(\frac12+i\tau) = -1 +1/s$ cela donne la
formule de la proposition.

\begin{remarque}
  On peut considérer que la proposition précédente
  \ref{prop:poissonL2} est un corollaire. On a en effet par
  Parseval:
  \[ 
\begin{split}
  G_f(x) &= \int_0^\infty f(xu)\frac{\{u\}}u\,du =
  \int_0^\infty \frac{\wt f(v/x)}x \left( -\frac{\{v\}}v +
    \int_v^\infty
    \{u\}\frac{du}{u^2}\right)\,dv \\
  &= - \int_0^\infty \frac{\wt f(v/x)}x \frac{\{v\}}v\,dv +
  \int_0^\infty
  \left(\int_0^{u} \wt f(v/x)\,\frac{dv}{ux}\right)  \frac{\{u\}}u\,du \\
  &= - \int_0^\infty \frac{\wt f(u/x)}x \frac{\{u\}}u\,du +
  \int_0^\infty
  \left(\int_0^{1/x}\wt f(uw)\,dw\right)  \frac{\{u\}}u\,du \\
  &= -\frac1x \int_0^\infty \wt f(\frac ux)\frac{\{u\}}u\,du
  + \int_0^{1/x} \left(\int_0^\infty \wt
    f(uw)\frac{\{u\}}u\,du\right)\,dw
 \end{split}
\]
Ainsi (les dérivées sont au sens des distributions
sur $]0,\infty[$):
\[ x\frac{d}{dx} G_f(x) = 
- \frac{d}{dx}\int_0^\infty \wt f(\frac ux)\frac{\{u\}}u\,du
= - \frac{d}{dx} G_{\wt f}(\frac1x) = \frac1{x^2}G_{\wt
  f}'(\frac1x)\]
ce qui donne exactement par \eqref{eq:Dc}:
\[ D_f(x) = \frac1x D_{\wt f}(\frac1x) \]
et donc prouve \ref{prop:poissonL2}.
\end{remarque}

\begin{remarque}
  La preuve de \ref{prop:fourierfract} montre plus
  généralement que pour une fonction paire $f\in
  L^2(\RR,dx)$, l'identité pour $y>0$:
  \[ \wt f(y) = -f(y) + \int_y^\infty f(x)\frac1x\,dx \; , \]
  équivaut pour sa transformée de Mellin gauche $\wh f(s) =
  \int_0^\infty f(x) x^{s-1}\,dx$, $\Reel(s)=\frac12$,
  d'être de la forme $Z(1-s)/(1-s)$, où $Z$ vérifie la même
  équation fonctionnelle que la fonction $\zeta(s)$. Cela
  équivaut aussi, comme on le voit facilement, à ce que la
  distribution tempérée $D(x) = \frac{d}{dx} xf(x)$ soit
  invariante sous Fourier (ici $D(x) = 1 - \sum_{n\neq0}
  \delta_n(x)$).
\end{remarque}

\section{Entrelacement et fonctions méromorphes}

\begin{remarque}
  Dans tout ce chapitre la transformée de Mellin est la
  transformée de Mellin droite:
  \[ \wh f(s) = \int_0^\infty f(t)t^{-s}\,dt\]
\end{remarque}

La fonction paire de carré intégrable $u\mapsto
g(u)=\frac{\{|u|\}}{|u|}$ a les deux propriétés suivantes:
\begin{enumerate}
\item elle est constante (égale à $1$) sur $]-a,a[$ avec
  $a=1$, et sa transformée de Fourier (par la proposition
  \ref{prop:fourierfract}) est de la forme $A+B\log|u|$ sur
  le même intervalle ($A=-\gamma, B=-1$).
\item sa transformée de Mellin complète
\[\pi^{-s/2}\Gamma(\frac s2)\int_0^\infty g(t)t^{-s}\,dt =
-\pi^{-s/2}\Gamma(\frac s2)\frac{\zeta(s)}{s}\] est une
fonction méromorphe dans tout le plan complexe, avec un pôle
double en $s=0$ et un pôle simple en $s=1$.
\end{enumerate}

Nous allons voir dans ce chapitre, non seulement pour les
fonctions, mais aussi pour les distributions, que la
deuxième propriété est une conséquence de la première, et
aussi comment l'on peut renforcer la deuxième propriété pour
la rendre équivalente à la première. Nous étudierons tout
d'abord les distributions (tempérées), paires, nulles et de
Fourier nulles dans $]-a,a[$, puis dans un deuxième temps
des situations plus générales. Nous aurons besoin de notions
sur la convolution multiplicative, la transformation de
Mellin pour les distributions, et les distributions
quasi-homogènes. Ne connaissant pas de référence commode
pour les résultats qui nous seront nécessaires, nous y
consacrons quelques développements rapides.

Nous emploierons souvent la notation $\b D,\phi\k$ pour
représenter l'appariement d'une distribution $D$ et d'une
fonction test $\phi$.

\subsection{Convolution multiplicative}

Nous notons $\cL_c(\RR^\times)$ l'espace des fonctions
intégrables sur $\RR$ qui ont un support compact éloigné de
l'origine. Soit $g\in \cL_c(\RR^\times)$. Comme on le
vérifie aisément l'application linéaire:
\begin{equation}
\phi(x)\mapsto  \int_\RR g(t)\phi(tx)\,dt
\end{equation}
va de $\cD(\RR) = \cC_c^\infty(\RR)$ vers lui-même et est
continue pour sa topologie.  Elle va aussi de $\cS(\RR)$
vers lui-même et est continue pour sa topologie.

Nous obtenons par dualité des applications linéaires
continues,
\begin{equation}
D(x)\mapsto (g*D)(x)
\end{equation}
de l'espace des distributions, resp. des distributions
tempérées, vers lui-même. Si $D$ est tempérée les deux
définitions donnent la même convolution multiplicative
$g*D$. Aussi, si $D$ est en fait elle-même une fonction test
$\phi(x)$ alors \[(g*\phi)(x) = \int_\RR
g(t)\phi(x/t)\,\frac{dt}{|t|}\;.\] Nous voyons donc que $
\int_\RR g(t)\phi(tx)\,dt$ peut aussi être écrite sous la
forme $(I(g)*\phi)(x)$ avec $I(g)(t) = g(1/t)/|t|$. La
définition de $g*D$ est donc
\[\b g*D, \phi\k  := \b D,I(g)*\phi\k \;.\]

Supposons que la fonction $g(x)$ soit en fait supportée dans
$]0,\infty[$. La restriction de $g*D$ à l'intervalle
$]0,\infty[$ ne dépend alors que de la restriction de $D$ à
$]0,\infty[$. Un changement de variable $x = \exp(u)$,
$u\in\RR$, ramène l'étude de $g*D$ sur $]0,\infty[$ au cas
de la convolution additive usuelle d'une distribution avec
une fonction à support compact, pour laquelle nous disposons
des résultats inclus dans les traités classiques
\cite{schwartz, hormander}. Par exemple, le théorème de
Titchmarsh-Lions \cite[IV]{hormander} donne ici:

\begin{prop}\label{prop:supportconvol}
  Soit $g(x)$, non-nulle, supportée dans $[a,A]$, avec $a>0$
  le plus petit point du support fermé essentiel de $g$.
  Supposons que la restriction de $D$ à $]0,\infty[$ ait le
  plus petit point de son support fermé en $b>0$. Alors le
  plus petit point du support fermé de la restriction de
  $g*D$ à $]0,\infty[$ est $ab$. En particulier $g*D\neq 0$.
\end{prop}

La propriété de régularisation (lissage) associée à la
convolution ne vaut ici qu'en dehors de l'origine. Par
exemple, la convolution multiplicative du Dirac à l'origine
$\delta(x)$ est un multiple de $\delta(x)$. Plus
généralement on a:

\begin{lem}
  Soit $P(T)$ un polynôme, et soit $g\in \cL_c(\RR^\times)$.
  On a:
  \[ g* \big(\,P(\frac{d}{dx})\delta\,\big) =  \Big(\int_\RR
  g(t)P(t\frac{d}{dx})\,dt\Big)\delta \]
\end{lem}

Il suffit de le prouver pour $P(T) = T^N$, $N\in\NN$. On
calcule (rappelons $\b \delta^{(k)}, \phi\k = (-1)^k
\phi^{(k)}(0)$):
\[ \begin{split}
  \b g * \delta^{(N)},\phi\k &=
  \left(-\frac{d}{dx}\right)^N\Big(\int_\RR
  g(t)\phi(tx)\,dt\Big)(0)\\
  &= \Big(\int_\RR g(t) t^N\,dt\Big)\b \delta^{(N)},\phi\k
\end{split}\]
et donc
\[ g * \delta^{(N)} = \Big(\int_\RR g(t)(t\frac{d}{dx})^N\,dt\Big)\delta\;.\]
En dehors de l'origine, la convolution régularise:

\begin{lem}\label{lem:conv1}
  Soit $g\in\cL_c(\RR^\times)$ (essentiellement) bornée.
  Soit $F$ une fonction localement intégrable. Alors $g*F$
  est une fonction localement intégrable qui est donnée par:
  \[\forall t\in\RR\quad (g*F)(t) =  \int_\RR g(x)F\left(\frac
    tx\right)\frac{dx}{|x|} = \int_\RR
  \frac{g(1/x)}{|x|}F(tx)dx\]
  Pour $t\neq 0$ on a:
  \[(g*F)(t)  =  \int_\RR\frac{g(t/x)}{|x|}F(x)dx\]
  Si $F$ est continue, resp. de classe $C^N$, alors $g*F$
  est continue, resp. de classe $C^N$, sur $\RR$ tout
  entier. Si $g$ est continue, resp. $C^N$, alors $g*F$ est
  continue, resp. $C^N$, sur $\RR^\times$.
\end{lem}

Pour la preuve on applique le théorème de Fubini à
l'intégrale \hbox{$\b g*F,\phi\k$}. Les propriétés de
continuité et de dérivabilité découlent des formules
intégrales.

\begin{lem}\label{lem:conv2}
  On a
  \[(g * D)^\prime = \frac{g(x)}x*D^\prime\]
\end{lem}

En effet
\begin{eqnarray*}
\b \frac{g(x)}x*\frac d{dx}D, \phi\k  &=& \b \frac d{dx}D, \int
\frac{g(t)}t\phi(tx)dt\k \\
&=& -\b D , \int g(t)\phi^\prime(tx)dt\k \\
&=&-\b g*D , \frac d{dx}\phi\k  = \b \frac d{dx}(g*D), \phi\!>\;,
\end{eqnarray*}
ce qui complète la preuve.

\begin{prop}\label{prop:conva}
  On suppose que la fonction $g$ est dans
  $\cC_c^\infty(\RR^\times)$. La restriction de la
  convolution multiplicative $g*D$ à $\RR^\times$ est alors
  une fonction, qui est donnée par la formule:
  \[ t\neq 0\Rightarrow\quad (g*D)(t) = \b  D,
  \frac{g(t/x)}{|x|}\k \]
  Cette fonction est de classe $C^\infty$ sur $\RR^\times$.
\end{prop}

L'application $x\mapsto g(t/x)/|x|$, pour $t\neq0$ donné,
est de classe $C^\infty$ et à support compact, donc la
formule a au moins un sens. On peut tout aussi bien pour la
preuve supposer que $g$ est supportée dans $x>0$. Un
changement de variable réduit l'énoncé au cas additif. Pour
une autre démonstration, on peut commencer par dire que le
terme de droite est une fonction $C^\infty$ de $t\neq0$, par
\cite[IV.1]{schwartz}. Nous savons par le lemme
\ref{lem:conv1} que la formule est valable lorsque $D$ est
une fonction continue. Et pour étudier $g*D$ dans un
voisinage d'un certain $t\neq0$ on peut tout aussi bien
supposer $D$ de support compact (si $\theta\in\cD(\RR)$ est
$1$ sur un intervalle suffisamment grand alors $g*D =
g*\theta D$ dans un voisinage de $t$). Toute distribution à
support compact est la dérivée à un certain ordre d'une
fonction continue, donc il suffira de montrer que si la
formule vaut pour $D$ elle vaut pour $D^\prime$. De
\ref{lem:conv2} on a:
  \[(g*D^\prime)(t) = \frac d{dt}(xg(x) * D)(t) =  \frac d{dt}\b D(x),
  \frac{tg(t/x)}{x|x|}\k \]
  Par \cite[IV.1]{schwartz} on peut mettre la dérivée à
  l'intérieur de l'appariement:
  \[ (g*D^\prime)(t) =  \b D, \frac\partial{\partial
  t}\frac{tg(t/x)}{x|x|}\k \] 
  Par homogénéité:
  \[ t\frac\partial{\partial
    t}\frac{tg(t/x)}{|x|} = - x\frac\partial{\partial x}
  \frac{tg(t/x)}{|x|}\quad\Rightarrow\quad
  \frac\partial{\partial t} \frac{tg(t/x)}{x|x|} =
  -\frac\partial{\partial x} \frac{g(t/x)}{|x|}\] Donc:
  \[ (g*D^\prime)(t)  = - \b D(x), \frac\partial{\partial x}
  \frac{g(t/x)}{|x|}\k =
  \b D^\prime , \frac{g(t/x)}{|x|}\k \;,\]%
  ce qui complète la preuve.

\begin{cor}
  Soit $a\neq0$ et soit $P(T)$ un polynôme. Soit
  $g\in\cC_c^\infty(\RR^\times)$. La convolution
  multiplicative $g*P(\frac d{dx})\delta(x-a)$ est la
  fonction de classe $C^\infty$, à support compact éloigné
  de l'origine, donnée par la formule:
  \[t\mapsto
  \left.P\left(-\frac{\partial}{\partial
        x}\right)\frac{g(t/x)}{|x|}\right|_{x=a}\]
\end{cor}

\begin{remarque}
  Par le lemme \ref{lem:conv2} on a
  \[(g*\delta^{(j)}_a)(t) = \left(\frac d{dt}\right)^j(x^jg(x)*\delta_a(x))(t)
  = \left(\frac \partial{\partial t}\right)^j \frac{t^j
    g(t/a)}{a^j |a|}\]
  et la compatibilité avec la formule du corollaire équivaut
  à l'identité:
  \[\left(\frac \partial {\partial t}\right)^j t^j\;
  \frac{g(t/a)}{|a|} = (-1)^j a^j \left(\frac \partial
    {\partial a}\right)^j \frac{g(t/a)}{|a|}\;,\]
  que l'on peut bien sûr vérifier directement. Nous avons
  employé la notation $\delta^{(j)}_a(x) = \left(\frac
    d{dx}\right)^j \delta(x-a)$. Ainsi $\b \delta^{(j)}_a,
  \phi\k = (-1)^j \phi^{(j)}(a)$.
\end{remarque}

\subsection{Le théorème d'entrelacement}

La transformation de Fourier $\cF$ entrelace les dilatations
et les contractions. Ceci donne un théorème simple et
important:

\begin{theo}\label{theo:entrelacement}
  Soit $D$ une distribution tempérée et soit
  $g\in\cL_c(\RR^\times)$ (c'est-à-dire $g$ est intégrable
  et de support essentiel compact et éloigné de l'origine).
  Soit $I(g)$ la fonction $g(1/t)/|t|$. On a alors la
  formule d'entrelacement de la convolution multiplicative
  et de la transformation de Fourier:
  \[\cF(g*D) = I(g)*\cF(D)\]
\end{theo}

Soit $\phi$ une fonction de Schwartz. Soit $k=I(g)$. On a
(les intégrales étant toutes étendues sur $\RR$ et
absolument convergentes):
\begin{eqnarray*}
(k*\cF(\phi))(x) &=& \int g(t)\cF(\phi)(tx)dt = \int g(t)\left(\int
e^{2\pi i txz}\phi(z)dz\right)dt\\
&=& \int\!\!\int e^{2\pi i tzx}\phi(z)g(t)\,dt\,dz
= \int\!\!\int e^{2\pi i wx}\phi(\frac wt)g(t)\,\frac{dt}{|t|}\,dw\\
&=& \int\!\!\int e^{2\pi i wx}\phi(uw)\frac{g(1/u)}{|u|}\,du\,dw
= \cF(g*\phi)(x)\\
\end{eqnarray*}
Ainsi $k*\cF(\phi) = \cF(g*\phi)$ et donc en dualisant:
  \begin{multline*}
    \b g*D, \cF(\phi)\k = \b D, k*\cF(\phi)\k = \b D,
    \cF(g*\phi)\k \\= \b\cF(D), g*\phi\k = \b k*\cF(D),
    \phi\k\;,\end{multline*} ce qui complète la preuve du
  théorème.

\subsection{Transformation de Mellin}

Nous développons ici la notion de transformation de Mellin
pour une distribution. Nous supposerons dans cette section
que la distribution tempérée $D$ a son support dans
$[a,\infty[$ pour un certain $a>0$. Les définitions faites
ne dépendront pas du choix de $a$.

\medskip

Considérons tout d'abord le cas où $D(x)$ est une fonction
continue $C(x)$ avec croissance au plus polynomiale. Alors
\[ \wh C(s) = \int_0^\infty C(x)x^{-s}dx\]
existe comme fonction analytique de $s$ pour $\Reel(s)\gg0$.
Pour $\sigma\gg0$ la fonction analytique $\wh C(s)$ est
$O(a^{-s})$ uniformément dans le demi-plan
$\Reel(s)>\sigma$.

\begin{remarque}
  On notera bien que dans tout ce chapitre la transformée de
  Mellin est définie suivant la formule \og droite\fg
  $\int_0^\infty C(x)x^{-s}dx$.
\end{remarque}

Dans le cas général, on peut écrire $D(x)$ comme la
$N$\ieme\ dérivée d'une fonction continue avec croissance
au plus polynomiale, pour $N$ suffisamment grand. Il y a un
unique choix d'une telle fonction continue $C_N(x)$ qui soit
nulle sur $]-\infty,a]$ (et $C_N$ ne dépend pas du choix de
$a>0$). De $C_{N+1}^\prime(x) = C_N(x)$ on déduit
\[\wh{C_N}(s) = s\,\wh{C_{N+1}}(s+1)\]

\begin{definition}\label{def:mellin}
  Nous définissons la transformée de Mellin de la
  distribution tempérée $D$, dont le support est inclus dans
  $[a,\infty[$, $a>0$, comme étant la fonction analytique
  \[ \wh D(s) = s(s+1)\dots(s+N-1) \wh{C_N}(s+N)\]
  avec comme domaine de définition le plus grand demi-plan
  $\Reel(s)>\sigma$ sur lequel un prolongement analytique
  soit possible (à partir de $\Reel(s)\gg0$), ou $\CC$ tout
  entier éventuellement. Pour des raisons typographiques, il
  sera parfois plus aisé d'écrire $M(D)(s)$ en lieu et place
  de $\wh D(s)$.
\end{definition}

Prenons note des propriétés suivantes:

\begin{fact}\label{fact:unicite}
  Si $\wh D(s) \equiv 0$ alors $D \equiv 0$.
\end{fact}

\begin{fact}
  On a $\wh{\frac d{dx}D}(s) = s\wh D(s+1)$.
\end{fact}

\begin{fact}
  On a $\wh{xD}(s) = \wh D(s-1)$.
\end{fact}

En effet c'est certainement vrai si $D$ est une fonction
continue de croissance polynomiale. Et si c'est vrai pour
$D$ c'est aussi vrai pour $D^\prime$: $\wh{xD^\prime}(s) =
\wh{(xD)^\prime}(s) - \wh{D}(s) = s\wh{xD}(s+1) - \wh D(s) =
(s-1)\wh D(s) = \wh{D^\prime}(s-1)$.

\begin{fact}\label{fact:multmellinpars}
  On a $\wh{x\frac d{dx} D}(s) = (s-1)\wh D(s)$ et $\wh
  {\frac d{dx}x D}(s) = s \wh D(s)$.
\end{fact}

Le lemme technique suivant est utile:

\begin{lem}
  Soit $\theta(x)\in\cD(\RR)$ ($\cD(\RR) =
  \cC_c^\infty(\RR)$) avec support dans $[\frac12, 3]$ et
  valant identiquement $1$ sur $[1,2]$. Soit pour $M\geq2$:
  $\theta_M(x)$ définie comme valant $\theta(x)$ sur
  $]-\infty, 1]$, $1$ sur $[1, M]$, $\theta(x - M + 2)$ sur
  $[M,\infty[$. Soit $D$ une distribution tempérée avec
  support inclus dans $[a,\infty[$, $a>0$.  On a alors
  \[\Reel(s)\gg 0\Rightarrow\quad \wh D(s) = \lim_{M\to\infty} \b D,
  x^{-s}\theta_M(\frac xa)\k\;.\] Et l'on peut choisir
  $\sigma_1\gg0$ de sorte que la convergence soit uniforme
  dans toute bande $\sigma_1\leq \Reel(s)\leq \sigma_2$,
  $\sigma_2>\sigma_1$.
\end{lem}

Pour la preuve nous choisissons $N\gg0$ tel que $D$ soit la
$N$\ieme\ dérivée d'une fonction continue $C_N(x)$ de
croissance polynomiale, nulle pour $x\leq a$. Considérons
les fonctions entières:
\[A_M(s) = \,<\! D, x^{-s}\theta_M(\frac xa)\!>\, = (-1)^N
\int_a^\infty C_N(x)\left(\frac d{dx}\right)^N
(x^{-s}\theta_M(\frac xa))\,dx \; .\]
Pour $s$ fixé avec $\Reel(s)\gg0$ on a convergence dominée
pour $M\to\infty$ vers
\[(-1)^N \int_a^\infty C_N(x)\left(\frac d{dx}\right)^N (x^{-s})\,dx
= s(s+1)\dots (s+N-1)\int_a^\infty C_N(x) x^{-s-N}\,dx\; ,
\]
qui est en fait par définition $\wh D(s)$. Et la convergence
dominée sera uniforme en $s$ si
$\sigma_1\leq\Reel(s)\leq\sigma_2$ et $\sigma_1\gg0$.

Énonçons un certain nombre de corollaires, dont les preuves
(simples) seront pour la plupart omises:

\begin{eqnarray}
  M(x^{-w}D(x))(s) &= \wh D(s+w)\\
  M(\log(x)D(x))(s) &= - \frac d{ds} \wh D(s)
\end{eqnarray}

\begin{fact}
  Si $f$, avec son support dans $[a,\infty[$, $a>0$, vérifie
  $\int_a^\infty |f(x)| x^{-\sigma}dx < \infty$ alors $\wh
  f(s)$ définie selon \ref{def:mellin} coïncide sur le
  demi-plan $\Reel(s)>\sigma$ avec $\int_a^\infty
  f(x)x^{-s}dx$ en tant qu'intégrale de Lebesgue. En
  particulier si $D(x)$ appartient à $L^1(a,\infty; \,dx)$
  alors $\wh D(s) = \int_a^\infty D(x)x^{-s}dx$ pour
  $\Reel(s)>0$.
\end{fact}

\begin{fact}
  Si $D(x)$ appartient à $L^2(a,\infty;\,dx)$ alors $\wh
  D(s) = \int_a^\infty D(x)x^{-s}dx$ pour
  $\Reel(s)>\frac12$.
\end{fact}

\begin{fact}
  Si $D$ a son support dans $[a,A]$, $0<a<A<\infty$, alors
  sa transformée de Mellin est une fonction entière, qui
  peut être évaluée par la formule
\[ s\mapsto \b D,
x^{-s}\theta(x)\k\; , \] où $\theta(x)\in\cD(\RR)$ vaut
identiquement $1$ sur $[a,A]$.
\end{fact}

\begin{fact}
  Si $D = \sum_{n} a_n\delta_{\lambda_n}$ avec une suite
  strictement croissante (finie ou infinie) $\lambda_n$
  telle que $S(X) = \sum_{\lambda_n\leq X} |a_n|$ ait une
  croissance (au plus) polynomiale, alors
\[\Reel(s)\gg0\Rightarrow\quad \wh D(s) = \sum_n \frac{a_n}{\lambda_n^s}\]
\end{fact}

\begin{prop}\label{prop:mellin}
  Soit $D$ une distribution tempérée supportée dans
  $[a,\infty[$, $a>0$, et soit $g$ une fonction intégrable
  supportée dans $[b,B]$, $b>0$. Alors, dans tout demi-plan
  $\Reel(s)>\sigma$ où $\wh D$ existe:
\[\wh{g*D}(s) = \wh g(s)\wh D(s) \; .\]
\end{prop}

Si $D(x) = C(x)$ est continue de croissance polynomiale et
$\Reel(s)\gg0$ ceci résulte d'une application directe de
Fubini. Dans le cas général, on a par le lemme
\ref{lem:conv2} $(g*D)(t) = (x^Ng(x) * C_N)^{(N)}(t)$ donc
\[\begin{split}
  M(g*D)(s) = s\cdots(s+N-1)M(x^N g*C_N)(s+N)\\ =
  s\cdots(s+N-1)\wh g(s)\wh{C_N}(s+N) = \wh g(s) \wh D(s) \;
  , \end{split}\]
ce qui prouve la proposition.

\subsection{Propriété \textit{S} et transformées de Mellin entières}

\begin{definition}
  On dira qu'une distribution $D$ sur $\RR$ a la \emph{propriété
  $S$} si elle est tempérée, et si il existe un intervalle
  $]-a,a[$, $a>0$, sur lequel à la fois $D$ et $\cF(D)$ sont
  nulles.
\end{definition}

Si $D$ a la propriété $S$ il en est de même de sa partie
paire et de sa partie impaire.

Si $D$ est paire et $g\in \cL_c(\RR^\times)$ est impaire
alors $g*D = 0$. Si $D$ est impaire et $g$ est paire alors
$g*D = 0$. Pour éviter ces problèmes on prendra toujours
dans la suite $g$ avec son support dans $]0,+\infty[$. Si
$D$ est paire, $g*D$ l'est aussi. Si $D$ est impaire $g*D$
l'est aussi. Pour une telle fonction $g\in
\cL_c(\RR^{+\times})$ on pose $\wh g(s) = \int_0^\infty
g(t)t^{-s}\,dt$. Pour une distribution $D$, soit paire, soit
impaire, et nulle dans un intervalle $]-a,a[$, la
transformée de Mellin $\wh D(s)$ est définie selon
\ref{def:mellin} appliqué à $D|_{]0,\infty[}$. Par la
proposition \ref{prop:mellin} on a:
\[\wh{g*D}(s) = \wh g(s)\wh D(s) \; .\]

\begin{definition}
  On dira parfois que $g*D$ est la $D$-co-somme de $g$, et
  on utilisera parfois la notation $P'_D(g) = g*D$.
\end{definition}

\begin{prop}
  Si $D$ a la propriété $S$ et si $g\in
  \cL_c(\RR^{+\times})$ alors $g*D$ a la propriété $S$.
\end{prop}

En effet, par la proposition \ref{prop:supportconvol} la
convolution $g*D$ est nulle dans l'intervalle $]-ab,ab[$ si
$g$ est nulle dans $]-b,b[$. Par le théorème d'entrelacement
\ref{theo:entrelacement} $\cF(g*D) = I(g)*\cF(D)$ et elle
est donc nulle dans l'intervalle $]-a/B,a/B[$ ($g(t)=0$ pour
$|t|\geq B$). Donc $g*D$ a la propriété $S$.

\begin{theo}\label{theo:dansschwartz}
  Soit $D$ une distribution, qui est soit paire, soit
  impaire, et qui possède la propriété $S$. Soit $g$ une
  fonction à support compact $[b,B]$, $0<b$, et infiniment
  dérivable. Alors la $D$-co-somme de $g$ est une fonction
  dans la classe de Schwartz, qui possède aussi la propriété
  $S$.
\end{theo}

Il reste à montrer que $g*D$ est dans $\cS$. Rappelons qu'en
tout cas par la proposition \ref{prop:conva} elle est de
classe $C^\infty$ sur $\RR$ et qu'elle est donnée pour
$t\neq0$ par la formule:
  \[ t\neq 0\Rightarrow\quad (g*D)(t) = \b  D,
  \frac{g(t/x)}{|x|}\k \; . \]
  Soit $k = I(g)$ et soit $\gamma(y) = \int_\RR e^{- 2\pi i
    xy}k(x)\,dx$ de sorte que $k = \cF(\gamma)$, et que
  $\gamma\in\cS$. Notons que pour $t\neq 0$ donné on a:
  \[ \frac{g(t/x)}{|x|} = \frac{k(x/t)}{|t|} =
  \cF(\gamma(ty))(x) \; ,\]
  et donc pour $t\neq0$:
\[ (g*D)(t) = \b D, \cF(\gamma(t\,\cdot))\k  = \b \cF(D),
\gamma(t\,\cdot)\k \; .\]
Supposons par exemple que $D$ soit paire; c'est alors aussi
le cas de $E = \cF(D)$. Comme $E(y)$ est tempérée elle est
la dérivée $C^{(2N)}(y)$ d'une fonction paire continue de
croissance (au plus) polynomiale.  La restriction de cette
fonction à l'intervalle $]-a,a[$ est un polynôme de degré au
plus $2N-1$. Après avoir soustrait ce polynôme on pourra
supposer que $C$ elle-même est identiquement nulle dans
l'intervalle $]-a,a[$. On obtient donc la formule suivante,
pour $t\neq0$:
\begin{multline*}
  (g*D)(t) = \b C^{(2N)}(y), \gamma(t y)\k = \int_\RR
  C(y)t^{2N}\gamma^{(2N)}(ty)\,dy\\
  = t^{2N}\int_a^\infty C(y)( \gamma^{(2N)}(ty) +
  \gamma^{(2N)}(-ty))\,dy \; .
\end{multline*}
Par cette formule et $\gamma\in\cS$ le fait que $(g*D)(t)$
et toutes ses dérivées décroissent plus vite que tout
inverse de polynôme lorsque $t\to\pm \infty$ est apparent
(le cas $D$ impaire est traité de manière similaire). Donc
$g*D\in\cS$, ce qu'il fallait démontrer.
 
Nous en arrivons au théorème principal:

\begin{theo}\label{theo:mellindistribsonine}
  Soit $D$ une distribution, paire ou impaire, avec la
  propriété $S$. Sa transformée de Mellin $\wh D(s)$ est une
  fonction entière. Si $D$ est paire, alors $\wh D$ a des
  zéros triviaux en $0$, $-2$, $-4$, \dots\ et vérifie
  l'équation fonctionnelle
\[ \pi^{-\frac s2}\Gamma(\frac s2)\wh D(s) =
\pi^{-\frac{1-s}2}\Gamma(\frac{1-s}2)\wh{\cF(D)}(1-s) \; .\]
Si $D$ est impaire, alors $\wh D$ a des zéros triviaux en
$-1$, $-3$, $-5$, \dots\ et vérifie l'équation fonctionnelle
\[ \pi^{-\frac{1+s}2}\Gamma(\frac{1+s}2)\wh D(s) =
-i\pi^{-\frac{2-s}2}\Gamma(\frac{2-s}2)\wh{\cF(D)}(1-s) \]
\end{theo}

Traitons le cas $D$ paire. Prenons n'importe quelle fonction
$g$ non nulle, infiniment dérivable, supportée dans $[b,B]$,
$0<b<B<\infty$. Nous savons par le théorème
\ref{theo:dansschwartz} que $g*D$ est une fonction paire
dans $\cS$, identiquement nulle dans un voisinage de
l'origine. Sa transformée de Mellin $\int_0^\infty
(g*D)(t)t^{-s}\,dt$ est donc une fonction entière de $s$.
Par la proposition \ref{prop:mellin} on a
\[ \wh{g*D}(s) = \wh g(s)\wh D(s)\;,\]
ce qui fournit pour $\wh D(s)$ un prolongement méromorphe à
tout le plan complexe. En fait les pôles ne peuvent être
qu'aux zéros de $g$, et $g$ peut être choisie
arbitrairement, donc $\wh D(s)$ est bien une fonction
entière. En ce qui concerne la fonction de Schwartz $\phi =
g*D$ paire elle est certainement $L^2$. Par \eqref{eq:eqfct}
elle vérifie presque partout sur la droite critique
\[ \pi^{-\frac s2}\Gamma(\frac s2)\wh{\;\cF(\phi)\;}(s)
=\pi^{-\frac{1-s}2}\Gamma(\frac{1-s}2)\wh{\phi}(1-s) \;.\]
Par le théorème d'entrelacement \ref{theo:entrelacement} on
a $\cF(\phi) = \cF(g*D) = I(g)*\cF(D)$, qui est aussi dans
$\cS$ et est nulle dans un voisinage de l'origine. Donc
$\wh{\;\cF(\phi)\;}$ est aussi une fonction entière. Ainsi
l'identité vaut dans tout le plan complexe, comme une
identité de fonctions méromorphes, en fait entières car les
pôles éventuels du terme de droite ne peuvent en être du
terme de gauche et vice versa. Par ailleurs $\wh{I(g)}(s) =
\wh g(1-s)$, d'où $\wh{\;\cF(\phi)\;}(s) = \wh
g(1-s)\wh{\cF(D)}(s)$, tandis que $\wh{\phi}(s) = \wh
g(s)\wh D(s)$.  Ainsi:
\[ \pi^{-\frac s2}\Gamma(\frac s2)\wh g(1-s)\wh{\cF(D)}(s)
=\pi^{-\frac{1-s}2}\Gamma(\frac{1-s}2)\wh g(1-s)\wh D(1-s)
\;,\]
puis
\[ \pi^{-\frac s2}\Gamma(\frac s2)\wh{\cF(D)}(s)
=\pi^{-\frac{1-s}2}\Gamma(\frac{1-s}2)\wh D(1-s)
\; ,\]
ce qui complète la preuve. Les zéros triviaux sont dus à
l'équation fonctionnelle. En effet la fonction méromorphe
${\mathcal D}(s)= \pi^{-\frac s2}\Gamma(\frac s2)\wh D(s)$
ne peut avoir de pôles qu'en $s = -2n$, $n\in \NN$. Mais
${\mathcal D}(1-s)$ est par l'équation fonctionnelle une
fonction du même type, et donc ${\mathcal D}(s)$ ne peut
avoir de pôles qu'en $s= 1+2n$, $n\in\NN$. La fonction
${\mathcal D}(s)$ est donc en fait une fonction entière,
d'où les zéros triviaux.

Le cas impair est traité quasiment à l'identique.

\subsection{Fonctions modérées et propriété \textit{S}}

\begin{definition}\label{def:moderee}
  Nous dirons qu'une fonction entière $F$ est \emph{modérée}
  si elle vérifie les deux conditions suivantes:
\begin{enumerate}
\item sur toute droite $\Reel(s) = \sigma$, $F$ est de
  croissance (au plus) polynomiale.
\item dans tout demi-plan $\Reel(s)>\sigma$, $F$ est dans la
  classe de Nevanlinna $\cN(\Reel(s)>\sigma)$.
\end{enumerate}
\end{definition}

\begin{remarque}
  On a renoncé à la terminologie plus précise \og $\theta$-modérée
  \fg, où $\theta$ serait un angle et ici $\theta = 0$.
\end{remarque}

La classe de Nevanlinna $\cN$ d'un demi-plan ouvert est
  composé des fonctions analytiques que l'on peut écrire
  dans ce demi-plan sous la forme d'un quotient de deux
  fonctions analytiques bornées. Les fonctions méromorphes
  quotients de deux fonctions analytiques bornées sont dites
  de \og caractéristiques bornées\fg. Pour ces notions et un
  certain nombre de résultats que nous supposerons connus et
  utiliserons sans autre référence, nous renvoyons au livre
  de Rosenblum et Rovnyak \cite{rosrov} qui contient tout ce
  qui nous est nécessaire.

\begin{remarque}\label{rem:nevanlinna}
  Faisons tout de suite la remarque qu'une fonction modérée
  $F(s)$ est, dans tout demi-plan $\Reel(s)\geq\sigma_0$
  donné, $O(A^{\Reel(s)}(1+|s|)^N)$ pour un certain $A>0$ et
  un certain $N\in\NN$. En effet par le théorème de
  factorisation de Smirnov-Nevanlinna, et la terminologie
  introduite par Beurling, on peut factoriser $F$ (supposée
  non-nulle\dots) dans le demi-plan ouvert
  $\Reel(s)>\sigma_0$ en un produit de Blaschke (qui est
  borné par $1$), un facteur intérieur spécial
  $A^{s-\sigma_0}$, un facteur intérieur associé à une
  mesure singulière signée sur la droite
  $\Reel(s)=\sigma_0$, qui est par force \emph{inexistant}
  ici à cause du prolongement analytique de $F$ à travers
  cette droite, et un facteur extérieur dont le module a
  comme logarithme l'intégrale de Poisson de
  $\log|F(\sigma_0 + it)|$. Par hypothèse ceci est $\leq
  \log C + N\log\surd(1+t^2)$ sur la droite, pour un certain
  $N\in\NN$ et une constante $C>0$, et par ailleurs
  l'intégrale de Poisson de $\log\surd(1+t^2)$ est
  $\log|s-\sigma_0+1|$, donc $F(s) =
  O(A^{\Reel(s)}(1+|s|)^N)$ dans le demi-plan
  $\Reel(s)>\sigma_0$ (ou $\geq \sigma_0$).  Signalons de
  plus que le $A$ du facteur intérieur spécial est donné par
  la formule:
  \[ \log(A) = \limsup_{\sigma\to+\infty}
  \frac{\log|F(\sigma)|}\sigma \; ,\]
  et ne dépend donc pas du demi-plan considéré. Suivant la
  terminologie de de~Branges \cite{bra}, on dit que
  $\log(A)$ est le \og type moyen\fg\ de $F$ dans le
  demi-plan considéré. Notons de plus le fait utile suivant:
  si $a>0$ est dans le support fermé essentiel de $f\in
  L^1(a,\infty;\,x^{-\sigma_0}dx)$ alors la fonction
  analytique $F(s) = \int_a^\infty f(x)x^{-s}\,dx$ est de
  Nevanlinna dans le demi-plan $\Reel(s)>\sigma_0$ et son type
  moyen est $\log\frac1a$. De même si $a>0$ est la borne
  inférieure du support fermé essentiel de $f\in
  L^2(0,\infty;\,dx)$ alors la fonction $\wh f(s) =
  \int_0^\infty f(x)x^{-s}\,dx$ est de Nevanlinna dans le
  demi-plan $\Reel(s)>\frac12$ et son type moyen est
  $\log\frac1a$.
\end{remarque}

\begin{remarque}
  La propriété pour
  une fonction entière $F$ d'être modérée équivaut à l'une ou
  l'autre des assertions:
  \[ \forall\sigma_0\in\RR\ \exists A>0\ \exists N\in\NN \quad
  F(s) = O(A^{\Reel(s)}(1+|s|)^N)\ \text{dans}\ 
  \{\Reel(s)>\sigma_0\}\]
  \[ \exists A>0\ \forall\sigma_0\in\RR\ \exists N\in\NN \quad
  F(s) = O(A^{\Reel(s)}(1+|s|)^N)\ \text{dans}\ 
  \{\Reel(s)>\sigma_0\}\]
  car elles impliquent $(1)$ et $(2)$ de \ref{def:moderee},
  et sont impliquées par \ref{def:moderee}, par la remarque
  qui précède.
\end{remarque}

Nous rappelons la notation
\[ \chi(s) = \frac{\zeta(s)}{\zeta(1-s)} =
\pi^{s-\frac12}\frac{\Gamma(\frac{1-s}2)}{\Gamma(\frac s2)}
= \frac{(2\pi)^{s}}{2\cos(\frac{\pi s}2)\Gamma(s)}\]
et le fait que cette fonction méromorphe et son inverse
$\chi(1-s)$ sont de croissance au plus polynomiale lorsque
$|\Imag(s)|\to\infty$ dans toute bande
$-\infty<\sigma_1\leq\Reel(s)\leq \sigma_2<\infty$.

\begin{theo}
  \label{theo:conversepair}
  Soit $D$ une distribution tempérée paire avec la
  propriété $S$. La fonction entière $F(s) =\wh D(s)$ est
  modérée et $\chi(s)F(1-s)$ est aussi une fonction entière
  et modérée. Réciproquement si $F$ est une fonction entière
  modérée et si $\chi(s)F(1-s)$ est aussi entière et modérée
  alors il existe une (unique) distribution tempérée paire
  $D$ avec la propriété $S$ telle que $F$ en soit la
  transformée de Mellin.
\end{theo}

Supposons que la distribution paire $D$ ait la propriété
$S$, pour l'intervalle $]-a,a[$, $a>0$. 
Soit $g$ une fonction non nulle quelconque, choisie
infiniment dérivable et de support dans $[b,B]$,
$0<b<B<\infty$.
Nous savons alors par \ref{theo:dansschwartz} que $\phi:=
g*D$ est une fonction paire, nulle dans $]-c,c[$, $c=ab$, et
est dans la classe de Schwartz. Soit $\sigma\in\RR$. Pour
$\Reel(s)\geq \sigma$, on a $|\wh{g*D}(s)| \leq
\int_c^\infty |\phi(x)|x^{-\Reel(s)}\,dx \leq
c^{-\Reel(s)}\int_1^\infty
|\phi(cx)|x^{-\sigma}x^{-\Reel(s)+\sigma}c\,dx =
O_\sigma(c^{-\Reel(s)})$. La fonction $c^s\wh\phi(s)$ est
ainsi bornée dans $\Reel(s)\geq\sigma$ et il en est de même
de $b^s\wh g(s)$. Comme $\wh D(s) = \wh\phi(s)/\wh g(s)$ on
en déduit que $F(s)=\wh D(s)$ est dans la classe de
Nevanlinna du demi-plan $\Reel(s)>\sigma$. Il en est de même
de $G(s) =\chi(s)F(1-s)$ puisque par le théorème
\ref{theo:mellindistribsonine} on a $G(s) = \wh E(s)$, $E =
\cF(D)$. 

Remarquons que par la définition \ref{def:mellin},
on a pour un certain $N\in \NN$:
\[ \Reel(s)\gg 0\Rightarrow\quad \wh D(s) =
s(s+1)\dots(s+N-1) \wh{C_N}(s+N) \; , \]
où $C_N$ est une fonction continue de croissance
polynomiale, nulle sur $]0,a[$.
La fonction $\int_a^\infty C_N(x)x^{-s-N}\,dx$ est
$O_{\sigma_0}(a^{-\Reel(s)})$ dans tout demi-plan
$\Reel(s)\geq\sigma\geq\sigma_0$, $\sigma_0\gg0$. Donc $F(s)
= \wh D(s)$ est $O(A^{\Reel(s)}(1 + |s|)^N)$ dans un tel
demi-plan ($A = a^{-1}$). Appliquons ce résultat à $E =
\cF(D)$. Nous obtenons en particulier que sur toute droite
$\Reel(s)=\sigma_1$, $\sigma_1\ll0$, la fonction $\wh
E(1-s)$ donc aussi $\chi(s)\wh E(1-s) = \wh D(s)$ est de
croissance polynomiale. Or nous savons déjà que $\wh D$ est dans
$\cN(\Reel(s)>\sigma_0)$, pour tout $\sigma_0\in\RR$, donc par
la remarque \ref{rem:nevanlinna}, le fait qu'elle soit de
croissance polynomiale sur la droite $\Reel(s) = \sigma_1$
implique qu'elle l'est aussi sur toute droite $\Reel(s) =
\sigma$, $\sigma>\sigma_1$. Comme $\sigma_1\ll0$ est
arbitraire, nous avons montré que $F(s) = \wh D(s)$ est de
croissance polynomiale sur toute droite verticale. Nous
avons donc établi que $F$ vérifie $(1)$ et $(2)$ de la
définition \ref{def:moderee}. 
Comme de plus $\chi(s)F(1-s) = \wh E(s)$, avec $E = \cF(D)$,
il est aussi établi que $\chi(s)F(1-s)$ est une fonction
entière modérée.

Il nous faut montrer la réciproque. Soit donc $F$ une
fonction entière modérée, telle que $\chi(s)F(1-s)$
soit aussi une fonction entière modérée. Nous savons par la
remarque \ref{rem:nevanlinna} qu'il existe un certain $a>0$
et un $N\in\NN$ tel que $a^s F(s)$ est $O(|s|^N)$ dans
$\Reel(s)\geq\frac12$. La fonction $a^s F(s)/s^{N+1}$ étant
$O(1/|s|)$ dans ce demi-plan, les normes $L^2$ sur les
droites verticales de ce demi-plan sont uniformément
bornées. Par un théorème de Paley-Wiener, cela implique
qu'il existe une fonction $f\in L^2(1,\infty;\,dx)$ avec
\[ a^s \frac{F(s)}{s^{N+1}} = \int_1^\infty f(x)x^{-s}\,dx =
\int_a^\infty f(x/a)a^{s-1} x^{-s}\,dx \]
dans ce demi-plan. Notons $g(x) = f(x/a)/a$, pour $x>a$, $0$
pour $|x|\leq a$, $g(-x)$ pour $x<-a$. Ainsi
\[ F(s) = s^{N+1} \wh g(s)\; , \]
et donc (par \ref{fact:multmellinpars}), $F = \wh D$, avec
$D$ la distribution tempérée paire qui est égale à $(\frac
d{dx}\,x)^{N+1} g(x)$, et est donc aussi nulle dans l'intervalle
$]-a,a[$. 

En répétant l'argument pour $G(s) = \chi(s)F(1-s)$ en lieu
et place de $F$ on voit qu'il existe aussi $a'>0$,
$N'\in\NN$ et une fonction $k$ paire, de carré intégrable,
nulle sur $]-a',a'[$, telle que
\[ G(s) = s^{N'+1} \wh k(s)\; , \]
pour $\Reel(s)\geq\frac12$ (presque partout sur la droite
critique). Quitte à remplacer $N$ et $N'$ par $\max(N,N')$
et à adapter $g$ ou $k$ l'on peut supposer $N=N'$.  
Par \eqref{eq:eqfct} la transformée de Fourier $\wt g$ de la
fonction de carré intégrable $g$ a une transformée de Mellin
qui est donnée par la formule $\wh{\wt g}(s) = \chi(s)\wh
g(1-s)$ sur la droite critique, c'est-à-dire ici:
\[ \wh{\wt g}(s) = \chi(s)F(1-s)(1-s)^{-N-1}  =  G(s)
(1-s)^{-N-1} = \left(\frac s{1-s}\right)^{N+1} \wh k(s)\; .
\]
Rappelons que (via la transformée de Mellin droite)
l'opérateur $L$ de multiplication par $s/(1-s) = -1
+1/(1-s)$ agit sur $L^2(0,\infty;\,dx)$ selon $k(x)\mapsto
L(k)(x) = -k(x)+\int_x^\infty k(y)\frac1y\,dy$. La fonction
$L(k)$ est constante pour $0<x<a'$. La fonction $L(L(k))$
est affine en $\log x$ pour $0<x<a'$, etc\dots . La fonction
paire $\wt g = L^{N+1}(k)$ est donc polynomiale en $\log x$ de
degré au plus $N$ pour $0<x<a'$. 

On note par ailleurs qu'il est correct
d'écrire au sens des distributions (paires) sur $\RR$:
\[ x \, \frac d{dx} (\log|x|)^{k} = k (\log|x|)^{k-1}\;.\]
Donc la distribution paire
\[ (x \, \frac d{dx})^{N+1} \cF(g)(x) \]
est identiquement nulle dans $]-a',a'[$. Il ne reste plus
maintenant qu'à utiliser
\[ \cF \frac d{dx}\,x = - x\,\frac d{dx} \cF\,\]
pour obtenir
\[ D(x) = (\frac d{dx}\,x)^{N+1} g(x)\Rightarrow\quad \cF(D)(x) =
(-1)^{N+1} (x \, \frac d{dx})^{N+1} \cF(g)(x) \; , \]
et donc pouvoir conclure que la distribution paire tempérée
$D$, qui est nulle dans un intervalle $]-a,a[$, $a>0$, a sa
transformée de Fourier $\cF(D)$ identiquement nulle dans un
intervalle $]-a',a'[$, $a'>0$. La distribution $D$ a donc la
propriété $S$ et le théorème est complètement
démontré.

\medskip

Définissons:
\[ \chi^{\mathrm{sin}}(s) =
i\,\pi^{s-\frac12}\frac{\Gamma(\frac{2-s}2)}{\Gamma(\frac{s+1}2)}
= i\,\frac{(2\pi)^s}{2\sin(\frac{\pi s}2)\Gamma(s)}\;.
\]
%
Si $g$ est, par exemple, dans $L^2(0,\infty;dx)$, et est vue
comme une fonction \emph{impaire} alors sa transformée de
Fourier $k = \wt g$, $k(x) = i\,\int_0^\infty 2\sin(2\pi
xy)g(y)\,dy$ a une transformée de Mellin droite qui est
reliée à celle de $g$ par la formule (presque partout sur la
droite critique):
\[ \wh k(s) = \chi^{\mathrm{sin}}(s)\wh g(1-s) \;.\]
Si l'on préfèrait travailler avec la transformée en sinus
$\int_0^\infty 2\sin(2\pi xy)g(y)\,dy$ il vaudrait donc
mieux supprimer le \og $i$\fg\ dans la définition de
$\chi^{\mathrm{sin}}$. Quoi qu'il en soit, le théorème suivant
a une démonstration tout-à-fait semblable à celle de
\ref{theo:conversepair}, et elle est laissée au lecteur.

\begin{theo}
  \label{theo:converseimpair}
  Soit $D$ une distribution tempérée impaire avec la
  propriété $S$. La fonction entière $F(s) =\wh D(s)$ est
  modérée et $\chi^{\mathrm{sin}}(s)F(1-s)$ est aussi une
  fonction entière et modérée. Réciproquement si $F$ est une
  fonction entière modérée et si
  $\chi^{\mathrm{sin}}(s)F(1-s)$ est aussi entière et modérée alors
  il existe une (unique) distribution tempérée impaire $D$
  avec la propriété $S$ telle que $F$ en soit la transformée
  de Mellin.
\end{theo}

Nous établissons quelques résultats concernant les zéros de
$\wh D(s)$ lorsque $D$ a la propriété $S$. Par l'équation
fonctionnelle du théorème \ref{theo:mellindistribsonine}
nous savons déjà que $\wh D(s)$ a des zéros triviaux. Si $D$
est paire la fonction $\pi^{-\frac s2}\Gamma(\frac s2)\wh
D(s)$ est entière, et si $D$ est impaire la fonction
$\pi^{-\frac{1+s}2}\Gamma(\frac{1+s}2)\wh D(s)$ est entière:
un \og zéro non-trivial\fg\ est par définition un zéro de
cette fonction entière.

\begin{lem}
  Soit $D$ une distribution, paire ou impaire, avec la
  propriété $S$. Alors $\wh D(s)$ a au moins un zéro non trivial.
\end{lem}

Nous supposerons que $D$ est paire (et non identiquement
nulle), la preuve étant similaire dans le cas impaire. Soit
$f(s) = \pi^{-\frac s2}\Gamma(\frac s2)\wh D(s)$. On établit
aisément que $\pi^{-\frac s2}\Gamma(\frac s2)$ est
$O((1+|s-\frac12|)^{\frac12|s-\frac12|})$ dans le demi-plan
$\Reel(s)\geq\frac12$.  Comme $\wh D(s)$ est modérée, on en
déduit que $|f(s)|$ est $O((1+|s-\frac12|)^{|s-\frac12|})$
dans ce demi-plan. Mais $f(1-s)$ est une fonction du même
type, donc $|f(s)|$ est $O((1+|s-\frac12|)^{|s-\frac12|})$
dans tout le plan complexe. Si elle n'avait pas de zéro, on
en déduirait que la fonction harmonique $\Reel(\log f(s))$
serait $\leq O(1)+|s-\frac12|\log(1+|s-\frac12|)$ dans
$\CC$. Cela n'est possible que si $\log f(s)$ est de la
forme $\alpha s + \beta$. On obtiendrait dans ce cas:
\[ \wh D(s) = \frac{e^{\alpha s + \beta}}{\pi^{-\frac
  s2}\Gamma(\frac s2)} \; , \]
puis (par exemple), en invoquant à nouveau la formule de Stirling:
\[ \lim_{\sigma\to+\infty} \frac{\log|\wh
  D(\sigma)|}{\sigma} = -\infty \; .\] 
Mais ceci est impossible, puisque contradictoire avec la
formule pour le type moyen de $\wh D(s)$ dans les demi-plans
$\Reel(s)>\sigma_0$ (\textsl{cf.} remarque
\ref{rem:nevanlinna}).

\begin{lem}
  Soit $D$ une distribution paire (resp. impaire) avec la
  propriété $S$. Soit $\rho$ un zéro non-trivial de $\wh
  D(s)$. Alors il existe une distribution paire
  (resp. impaire) $D_1$ avec la propriété $S$ telle que 
  \[ D(x) = \left( \frac{d}{dx} x - \rho \right) D_1(x) \; .
  \]
  On a $\wh D(s) = (s - \rho)\wh D_1(s)$.
\end{lem}

En effet, soit $F_1(s)$ la fonction entière $F(s)/(s-\rho)$.
Il s'agit certainement d'une fonction modérée, et comme
$\rho$ n'est pas un zéro trivial, la fonction
$\chi(s)F_1(1-s)$ (resp. $\chi^{\textrm{sin}}(s)F_1(1-s)$)
est aussi entière (et modérée). Par le théorème
\ref{theo:conversepair} (resp.  \ref{theo:converseimpair})
$F_1(s) = \wh D_1(s)$ pour une certaine distribution $D_1$
avec la propriété $S$. On a $\wh D(s) = (s - \rho)\wh
D_1(s)$ et donc $D(x) = \left( \frac{d}{dx} x - \rho \right)
D_1(x)$ par la propriété \ref{fact:multmellinpars}.

\begin{prop}
  Soit $D$ une distribution paire (resp. impaire) avec la
  propriété $S$. Alors $\wh D(s)$ a une infinité de zéros
  non-triviaux. Il existe un entier $N$ tel que pour tout
  choix de $N$ zéros non-triviaux $\rho_1$, \dots, $\rho_N$,
  il existe une fonction de carré intégrable $C$ paire (resp.  impaire)
  avec la propriété $S$, telle que
  l'identité de distributions
  \[ D(x) = P(\frac{d}{dx} x) C(x) \]
  ait lieu avec $P(T) = \prod_{1\leq j\leq N} (T - \rho_j)$.
\end{prop}

Par les deux lemmes précédents et une récurrence, pour tout
choix de $N$ zéros non triviaux on aura $D(x) = \prod_{1\leq
  j\leq N} (\frac{d}{dx} x - \rho_j) D_N(x)$ avec une
distribution (dépendant des $N$ zéros) possédant la
propriété $S$. De plus $\wh D(s) = \prod_{1\leq j\leq N} (s
- \rho_j) \wh D_N(s)$; comme $\wh D(s)$ est une fonction
modérée, pour $N$ suffisamment grand, et indépendamment du
choix des zéros, on aura $\wh D_N(s) = O(1/|s|)$ sur la
droite critique.  La distribution $D_N$ est donc en fait une
fonction de carré intégrable puisque sa transformée de
Mellin est de carré intégrable sur la droite critique.

\begin{remarque}
  Si l'on prend un zéro de plus les fonctions $C$ seront
  continues et $o(1/\sqrt x)$ pour $x\to+\infty$ à cause de
  $C(x) = \int_\RR \wh D_N(\frac12 + i\tau) |x|^{-\frac12 +
    i\tau}\frac{d\tau}{2\pi}$.
\end{remarque}


\begin{theo}
  Soit $D$ une distribution paire (resp. impaire) avec la
  propriété $S$, non identiquement nulle. Soit
  $0<\epsilon<\frac\pi2$. Le nombre des zéros de modules au
  plus $T$ de la fonction entière $\wh D(s)$ dans le secteur
  angulaire $|\arg(s-\frac12) - \frac\pi2|<\epsilon$ est
  asymptotiquement équivalent à $\frac{T}{2\pi}\log(T)$, et
  de même pour $|\arg(s-\frac12) + \frac\pi2|<\epsilon$. Le
  nombre des zéros non-triviaux de modules au plus $T$ dans
  les secteur angulaires $|\arg(s)|<\frac\pi2 -\epsilon$ et
  $|\arg(s)-\pi|<\frac\pi2 -\epsilon$ est $o(T)$.
\end{theo}

Nous avons établi ce théorème pour les fonctions paires de
carrés intégrables avec la propriété $S$ dans \cite{twosys}
(Theorem $7.7$). Il vaut aussi pour les fonctions impaires
de carré intégrables par une preuve quasi-identique. Le
théorème dans le cas général peut alors être établi en
invoquant la proposition précédente. En fait on peut aussi
simplement vérifier que la preuve donnée pour \cite[Theorem
7.7]{twosys} s'adapte aux cas envisagés ici.

\subsection{Distributions homogènes et quasi-homogènes}

Rappelons tout d'abord quelques notions sur les
distributions homogènes (voir par exemple \cite{gelfand}).
On dira que la distribution $D$ sur $\RR$ est homogène si il
existe une fonction $\lambda:\RR^\times\to\CC$ telle que
$D(t x) = \lambda(t) D(x)$ pour tout $t \neq 0$. Si $D$
n'est pas la distribution nulle alors $\lambda$ est un
morphisme continu de $\RR^\times$ vers $\CC^\times$, qui est
donc soit $|t|^{s-1}$ pour un certain $s\in\CC$,
soit $\frac t{|t|} |t|^{s-1}$ pour un certain
$s\in\CC$. À une constante multiplicative près il existe une
unique distribution (non-nulle) d'une homogénéité donnée.

La fonction
\[ s\mapsto D^s(x) = \frac{|x|^{s-1}}{\pi^{-\frac s2}\Gamma(\frac
  s2)}
\]
est analytique sur le demi-plan $\Reel(s)>0$ à valeurs dans
les distributions, et elle admet un prolongement analytique
en une fonction entière et sans zéro sur $\CC$ à valeurs
dans les distributions paires sur $\RR$. Pour tout $s\in\CC$
la distribution $D^s$ est homogène d'homogénéité
$|t|^{s-1}$.  La restriction de $D^s$ à l'ouvert
$\RR\setminus\{0\}$ est la fonction $|x|^{s-1}/\pi^{-\frac
  s2}\Gamma(\frac s2)$ (donc nulle aux pôles du
dénominateur).

La fonction
\[ s\mapsto E^s(x) =
\frac{x}{|x|}\frac{|x|^{s-1}}{\pi^{-\frac{s+1}2}\Gamma(\frac{s+1}2)} 
\]
est analytique sur le demi-plan $\Reel(s)>0$ à valeurs dans
les distributions, et elle admet un prolongement analytique
en une fonction entière et sans zéro sur $\CC$ à valeurs
dans les distributions impaires sur $\RR$. Pour tout
$s\in\CC$ la distribution $E^s$ est homogène d'homogénéité
$\frac t{|t|} |t|^{s-1}$.

\begin{fact} Les identités suivantes sont vérifiées:
\[ \frac{d}{dx} x D^s(x) = s D^s(x)\qquad \frac{d}{dx} x
E^s(x) = s E^s(s)\]
\end{fact}

On fait le calcul pour $\Reel(s)>0$ et ensuite on invoque un
prolongement analytique en $s$.

\begin{fact}
  Soit $g\in \cL_c(\RR^\times)$ une fonction intégrable à
  support compact éloigné de l'origine. On a:
  \begin{eqnarray*}
    g*D^s = \Big(\int_\RR g(t)|t|^{-s}\,dt\Big)D^s\\
    g*E^s = \Big(\int_\RR \frac t{|t|} g(t)|t|^{-s}\,dt\Big)E^s
  \end{eqnarray*}
\end{fact}

On fait le calcul pour $\Reel(s)>0$ grâce au lemme
\ref{lem:conv2} et ensuite on invoque un prolongement
analytique en $s$.

\begin{fact}
Si $D$ est une distribution non-nulle d'homogénéité
$\lambda$ alors elle est tempérée et $\cF(D)$ est une
distribution non-nulle d'homogénéité $\lambda(1/t)/|t|$. On
a en fait très exactement:
\begin{eqnarray*}\label{eqn:tate}
\cF(D^s) = D^{1-s}\\
\cF(E^s) = i\;E^{1-s}
\end{eqnarray*}
\end{fact}

Il suffit d'évaluer sur les fonctions test $\exp(-\pi x^2)$
ou $x\exp(-\pi x^2)$.

Comme les fonctions $s\mapsto D^s$ et $s\mapsto E^s$ sont
des fonctions entières on peut les dériver par rapport à
$s$. Nous adopterons la définition suivante:

\begin{definition}
  On dira qu'une distribution $D$ sur $\RR$ est
  \emph{quasi-homogène} si elle est la combinaison linéaire
  complexe d'un nombre fini de distributions du type
  $D^{s,N} = \left.\left(\frac\partial{\partial z}\right)^N
    D^z\right|_{z=s}$ ou $E^{s,M} =
  \left.\left(\frac\partial{\partial z}\right)^M
    E^z\right|_{z=s}$.
\end{definition}

Une distribution quasi-homogène est déterminée par sa
restriction à n'importe quel intervalle $]-a,a[$, $a>0$. Sa
transformée de Fourier est aussi quasi-homogène.

\begin{fact}\label{fact:ddxdses}
Pour tout $s\in \CC$ et tout $N\in\NN$ on a 
\[\begin{matrix}
\left( \frac{d}{dx} x - s \right)  D^{s,N+1} = (N+1) D^{s,N}\hfill\\
\left( \frac{d}{dx} x - s \right)^{N} D^{s,N} = N! D^s\hfill\\
\left( \frac{d}{dx} x - s \right)^{N+1} D^{s,N} =
0\hfill
\end{matrix}
\qquad\begin{matrix}
\left( \frac{d}{dx} x - s \right) E^{s,N+1}= (N+1) E^{s,N}
\hfill \\
\left( \frac{d}{dx} x - s \right)^{N} E^{s,N} = N! E^s \hfill \\
\left( \frac{d}{dx} x - s \right)^{N+1} E^{s,N} = 0 \hfill
\end{matrix}\]
\end{fact}

On dérive $\frac {d}{dx} x D^s = s D^s$ par
rapport à $s$ et on obtient $\frac{d}{dx} x D^{s,N+1}(x) = s
D^{s,N+1}(x) + (N+1) D^{s,N}(x)$ (de même pour $E^s$). Les
assertions en découlent. Remarquons que $D^{s,N}$ et 
$E^{s,N}$ ne sont jamais la distribution nulle.

\begin{fact}
  Ce n'est pas vraiment une définition, ni une remarque,
  mais plutôt un commentaire heuristique: la transformée de
  Mellin d'une distribution quasi-homogène (paire ou
  impaire) est la
  fonction complexe identiquement nulle.
\end{fact}

\begin{definition}
  On dira qu'une distribution $D$ est $a$-quasi-homogène si
  sa restriction à un intervalle $]-a,a[$, $a>0$, coïncide
  avec la restriction à cet intervalle d'une distribution
  quasi-homogène $D^{qh}$. Celle-ci est uniquement
  déterminée et ne dépend pas de $a$.
\end{definition}

\begin{definition}\label{def:qhsurintervalle}
  Soit $D$ une distribution tempérée paire ou impaire
  $a$-quasi-homogène. La transformée de Mellin de $D$ est
  définie comme étant la transformée de Mellin au sens de
  \ref{def:mellin} de la distribution (paire ou impaire)
  $D-D^{qh}$.
\end{definition}

Les deux propositions suivantes sont vérifiées aisément:

\begin{prop}
  Soit $D$ une distribution paire (resp. impaire)
  $a$-quasi-homogène. Les distributions impaires
  (resp. paires) $xD$ et $D'$ sont aussi
  $a$-quasi-homogènes. Les formules
  \[ \wh{\frac d{dx} D}(s) = s\wh D(s+1)\qquad
  \wh{xD}(s) = \wh D(s-1) \]
  sont valables et donc aussi:
  \[ \wh{x \frac d{dx} D}(s) = (s-1)\wh D(s)\qquad 
  \wh{\frac d{dx} x D}(s) = s\wh D(s) \; .\]
\end{prop}

\begin{prop}
  Soit $D$ une distribution $a$-quasi-homogène, paire ou
  impaire. Soit $g\in \cL_c(\RR^{+\times})$. La convolution
  multiplicative $g*D$ est $a'$-quasi-homogène pour un
  certain $a'>0$ et la formule
  \[ \wh{g*D}(s) = \wh g(s)\wh D(s) \]
  est valable.
\end{prop}

Nous aurons aussi particulèrement besoin des lemmes suivant:

\begin{lem}\label{lem:inverseqh}
  Soit $D$ une distribution quasi-homogène et soit
  $s\in\CC$. Il existe une distribution quasi-homogène $D_1$
  avec $(\frac{d}{dx}x - s) D_1 = D$.
\end{lem}

Cela découle directement de \ref{fact:ddxdses} pour $D =
D^{s,N}$ ou $E^{s,N}$; et aussi pour $D = D^{w,N}$, $w\neq
s$, puisque l'opérateur $L = \frac{d}{dx}x - s$ sur l'espace
vectoriel des $D^{w,k}$ (resp. $E^{w,k}$), $0\leq k\leq M$
  est de la forme $w-s+N$ avec $N$ nilpotent.

\begin{lem}
  \label{lem:dsescaract}
  Soit $D$ une distribution vérifiant
  \[ \left( \frac{d}{dx} x  - s \right ) D = 0 \; .\]
  Il existe alors deux constantes (uniques) telles que $D =
  \alpha D^s + \beta E^s$.
\end{lem}

Soit $t>0$ et soit $\phi\in\cD(\RR)$. Soit $f(t) = \int_\RR
D(x)\phi(tx)\,dx$. Par \cite[IV.1]{schwartz} $f$ est une
fonction de classe $C^\infty$ de $t>0$ et l'on peut dériver
sous le signe \og intégrale\fg.
\[ f'(t) = \int_\RR
D(x)x\phi'(tx)\,dx = -\b \frac{d}{dx} x D, \frac1t
\phi(tx)\k = -s \frac1t f(t) \; .\] 
Donc $f(t) = C t^{-s}$ et 
\[ \b D(x/t)/t, \phi(x)\k = t^{-s} \b D(x), \phi(x)\k\; .\]
Autrement dit, $D(x/t)/t = t^{-s}D(x)$, $D(tx) =
t^{s-1}D(x)$ pour tout $t>0$. Si $D$ est paire, on obtient
que $D$ est homogène pour le quasi-caractère $|t|^{s-1}$, et
donc $D$ est un multiple de $D^s$. Si $D$ est impaire on
obtient qu'elle est un multiple de $E^s$. Dans le cas
général on applique le résultat séparément aux parties
paires et impaires de $D$ car elles vérifient la même
équation.



\begin{lem}
  \label{lem:inverse}
  Soit $D$ une distribution nulle
  dans un intervalle $]-a,a[$, $a>0$. Pour tout $s\in\CC$ il
  existe une distribution (unique) 
  $D_1$ qui est nulle dans $]-a,a[$ et qui vérifie:
  \[ \left( \frac{d}{dx} x  - s \right ) D_1(x) = D(x)
  \;. \]
   Si $D$ est tempérée, $D_1$ l'est aussi. Si $D$ est paire
   (resp. impaire) $D_1$ l'est aussi.
\end{lem}

Pour la preuve de cet autre énoncé standard, indiquons
simplement que comme $D$ est nulle dans $]-a,a[$ on peut
considérer la distribution identiquement nulle dans $]-a,a[$
et qui est donnée formellement pour $x\neq 0$ par la formule
\[ D_1(x) = |x|^{s}\frac1x \int_0^x |y|^{-s} D(y)\,dy \; . \]
Cette distribution est l'unique solution. 

\begin{prop}
  \label{lem:qhcaract}
  Soit $D$ une distribution telle qu'il existe un $a>0$ et
  un polynôme unitaire $P(T)$ tels que la distribution
  $P(\frac{d}{dx} x)D$ soit nulle dans l'intervalle
  $]-a,a[$. Alors $D$ est $a$-quasi-homogène. La réciproque
  vaut aussi.
\end{prop}

Par récurrence sur le degré de $P$, et il suffira de montrer
que si $E = (\frac{d}{dx} x - s)D$ est $a$-quasi-homogène
alors $D$ l'est aussi. Par le lemme précédent appliqué à $E
- E^{qh}$, il existe $D_1$ nulle dans $]-a,a[$ telle que
$(\frac{d}{dx} x - s)D_1 = E - E^{qh}$.  De plus par le
lemme \ref{lem:inverseqh} il existe une distribution
quasi-homogène $D_2$ vérifiant $(\frac{d}{dx} x - s)D_2 =
E^{qh}$. La distribution $D_3 = D_1 + D_2$ est
$a$-quasi-homogène et $(\frac{d}{dx} x - s)D_3 =
(\frac{d}{dx} x - s)D$. Il ne reste plus qu'à faire appel au
lemme \ref{lem:dsescaract}. La réciproque vaut par
\ref{fact:ddxdses}.

\subsection{Propriété \textit{S}-étendue et fonctions méromorphes} 

\begin{definition}
  \label{def:Setendu}
  On dira qu'une distribution tempérée $D$ a la
  \emph{propriété $S$-étendue} si il existe $a>0$ tel que à
  la fois $D$ et $\cF(D)$ sont $a$-quasi-homogènes.
\end{definition}

Si $D$ a la propriété $S$-étendue, il en est de même de $xD$
et de $D'$, ainsi que de $g*D$ pour toute fonction $g\in
\cL_c(\RR^\times)$, par le théorème d'entrelacement
\ref{theo:entrelacement}. Si $D$ a la propriété $S$-étendue,
il en est de même de $D + D_1$ pour toute distribution
quasi-homogène $D_1$. Les distributions quasi-homogènes ont
la propriété $S$-étendue, et ce sont les seules à l'avoir
pour tous les $a>0$.

\begin{definition}
  On dira qu'une fonction méromorphe $F(s)$ est
  \emph{modérée} s'il existe un polynôme non-nul $P(s)$ tel
  que la fonction $P(s)F(s)$ est entière et modérée au sens
  de la definition \ref{def:moderee}.
\end{definition}

\begin{theo}\label{theo:S-etendue}
  Soit $D$ une distribution paire avec la propriété
  $S$-étendue. Sa transformée de Mellin $F(s) = \wh D(s)$
  est alors une fonction méromorphe dans le plan complexe,
  modérée, et telle que la fonction $\chi(s)F(1-s)$ est
  aussi une fonction méromorphe modérée. Réciproquement si
  $F$ est une fonction méromorphe modérée telle que
  $\chi(s)F(1-s)$ est aussi méromorphe et modérée alors il
  existe une distribution paire avec la propriété
  $S$-étendue telle que $F(s) = \wh D(s)$; on peut fixer $D$
  de manière unique en lui imposant d'être nulle dans un
  voisinage de l'origine. Les seules autres distributions
  paires avec la propriété $S$-étendue et la même
  transformée de Mellin sont celles qui diffèrent de $D$ par une
  distribution quasi-homogène paire. On a $\wh{\cF(D)}(s) =
  \chi(s)\wh D(1-s)$. Le même énoncé vaut dans le cas impair
  en remplaçant $\chi(s)$ par $\chi^{\mathrm{sin}}(s)$.
\end{theo}

Soit $D$ une distribution paire avec la propriété
$S$-étendue, $D^{qh}$ sa partie quasi-homogène, $E =
\cF(D)$, $E^{qh}$ sa partie quasi-homogène. On a
\[ \cF(D - D^{qh}) = E - E^{qh} + E^{qh} - \cF(D^{qh})\; .\]
Il existe un polynôme unitaire $Q(T)$ tel que
$Q(\frac{d}{dx} x)$ annule la distribution quasi-homogène
$E^{qh}(x) - \cF(D^{qh})(x)$. Soit $D_1(x) =
Q(-x\frac{d}{dx})(D(x) - D^{qh}(x))$.  Comme $\cF x\frac{d}{dx}
= - \frac{d}{dx} x \cF$ on a
\[ \cF(D_1)(x) = Q(\frac{d}{dx} x)(E(x) - E^{qh}(x)) \; . \]
La distribution paire $D_1$ a donc la propriété $S$ stricte.
Sa transformée de Mellin $F_1(s) = \wh D_1(s)$ est ainsi une
fonction entière modérée telle que $\chi(s)F_1(1-s)$ est
aussi entière et modérée. Or 
\[ F_1(s) = Q(1-s)\wh D(s) \; , \]
donc $F(s) = \wh D(s)$ est une fonction méromorphe modérée,
et $\chi(s)F(1-s) = \chi(s)F_1(1-s)/Q(s)$ également.

Pour la réciproque, si $F$ est méromorphe et modérée, et
aussi $\chi(s)F(1-s)$ alors la fonction méromorphe
${\mathcal F}(s) = \pi^{-\frac s2}\Gamma(\frac s2)F(s)$ ne
peut avoir qu'un nombre fini de pôles. Soit $P(s)$ le
polynôme unitaire dont les racines sont ces pôles, avec
multiplicités éventuelles. La fonction $F_1(s) = P(s)F(s)$
est entière et modérée et la fonction $\chi(s)F_1(1-s) =
\chi(s)P(1-s)F(1-s) = P(1-s){\mathcal F}(1-s)/(\pi^{-\frac
  s2}\Gamma(\frac s2))$ l'est également. Donc, par le
théorème \ref{theo:conversepair} c'est qu'il existe une
distribution $D_1$ paire avec la propriété $S$ stricte pour
un certain intervalle $]-a,a[$, $a>0$, telle que $F_1(s) =
\wh D_1(s)$. Par le lemme \ref{lem:inverse} il existe une
distribution paire tempérée $D(x)$ (unique), identiquement
nulle dans $]-a,a[$, et vérifiant: $P(\frac{d}{dx} x)D(x) =
D_1(x)$. C'est donc que $P(s)\wh D(s) = \wh D_1(s)$ ou
encore $\wh D(s) = F(s)$. On a $P(-x\frac{d}{dx})\cF(D)(x) =
\cF(D_1)(x)$, donc par la proposition \ref{lem:qhcaract}, la
distribution $\cF(D)$ est $a$-quasi-homogène. Donc $D$ a la
propriété $S$-étendue. De plus $P(1-s)\wh{\cF(D)}(s) =
\wh{\cF(D_1)}(s) = \chi(s)\wh{D_1}(1-s) = \chi(s)P(1-s)\wh
D(1-s)$ donc $\wh{\cF(D)}(s) = \chi(s)\wh D(1-s)$. Si $D_2$
a la même transformation de Mellin, alors par la définition
\ref{def:qhsurintervalle} et par \ref{fact:unicite}, $D_2 -
D_2^{qh} = D - D^{qh}$, donc $D_2$ diffère de $D$ par une
distribution quasi-homogène. Réciproquement, si c'est le cas
alors $\wh{D_2} = \wh D$. Si de plus à la fois $D_2$ et $D$
sont nulles dans un voisinage de l'origine alors $D_2 = D$.
La démonstration du théorème dans le cas impair est
quasi-identique.

Dans ce théorème, on  a normalisé $D$ on lui imposant d'être
nulle dans un voisinage de l'origine. On peut tout aussi
bien imposer cette condition à $\cF(D)$, quitte à remplacer
$D$ par $D - \cF^{-1}(\cF(D)^{qh})$. C'est cette
normalisation que nous utiliserons pour décrire la partie
polaire de $\pi^{-\frac s2}\Gamma(\frac s2)\wh D(s)$.

\begin{theo}\label{theo:polaire}
  Soit $D$ une distribution paire avec la propriété
  $S$-étendue, normalisée de sorte que $\cF(D)$ ait une
  partie quasi-homogène nulle. Soit
  $\sum_{w\in\CC,\;k\in\NN} c_{w,k}D^{w,k}$ la partie
  quasi-homogène $D^{qh}$ de $D$ (avec un nombre fini de
  coefficients non nuls). La partie polaire de la
  transformée de Mellin complète $\pi^{-\frac
    s2}\Gamma(\frac s2)\wh D(s)$ est:
  \[ \sum_{w\in\CC,k\in\NN}  c_{w,k}\;\frac{(-1) k!}{(s-w)^{k+1}}
  \;.\]
  Dans le cas où $D$ est une distribution impaire avec la
  propriété $S$-étendue, normalisée de sorte que
  $\cF(D)^{qh} \equiv 0$, la partie polaire de
  $\pi^{-\frac{s+1}2}\Gamma(\frac{s+1}2)\wh D(s)$ est donnée
  par la même formule avec $D^{qh} = \sum_{w\in\CC,\;
    k\in\NN} c_{w,k}E^{w,k}$.
\end{theo}

La partie polaire ne dépend (linéairement, et indépendamment
de $a>0$) que de $D^{qh}$, puisque si $D^{qh} = D_1^{qh}$
alors $D-D_1$ a la propriété $S$-stricte et donc les
transformées de Mellin complètes de $D$ et de $D_1$
diffèrent par une fonction entière et ont la même partie
polaire. De plus si la partie polaire s'annule alors $D$ par
le théorème \ref{theo:S-etendue} diffère d'une distribution
avec la propriété $S$-stricte par une distribution
quasi-homogène. Comme $\cF(D)\equiv0$ dans un voisinage de
l'origine cela n'est possible que si cette distribution
quasi-homogène est nulle. Donc $D^{qh}$ est nulle. Il existe
ainsi une correspondance linéaire bijective entre les
distributions quasi-homogènes réalisables et les parties
polaires réalisables, qu'il s'agit de préciser. Notons que
toutes les parties polaires sont effectivement réalisées,
comme on le voit en partant d'une quelconque distribution
non-nulle $E$ avec la propriété $S$ stricte en formant des
combinaisons linéaires avec les fonctions
$(s-w)^{-k-1}\pi^{-\frac s2}\Gamma(\frac s2)\wh E(s)$,
$w\in\CC$, $k\in\NN$ ($k$ au moins égal à la multiplicité
éventuelle de $w$ comme zéro non-trivial).  Supposons alors
que la distribution $D$ avec $\cF(D)^{qh}\equiv0$ soit telle
que seul un certain $w$ soit un pôle, d'ordre au plus $k+1$,
$k\in\NN$.  Alors $(\frac d{dx} x - w)^{k+1}(D)$ a une
transformée de Mellin complète entière; elle diffère donc
d'une distribution avec la propriété $S$-stricte par une
distribution quasi-homogène, mais celle-ci doit être nulle
car sa transformée de Fourier doit elle aussi être
identiquement nulle dans un voisinage de l'origine. Donc
$(\frac d{dx} x - w)^{k+1}$ annule $D^{qh}$.  Les
distributions quasi-homogènes compatibles avec cette
condition sont les $D^{w,j}$, $0\leq j\leq k$, et leurs
combinaisons linéaires qui forment un espace de dimension
$k+1$. Toutes sont donc réalisées puisque toutes les parties
polaires sont réalisées et que les dimensions correspondent.

Soit donc $D$ telle que $D^{qh} = D^{w,0}$ pour un certain
$w\in\CC$ (et $\cF(D)^{qh}\equiv0$). Supposons tout d'abord
que $w$ n'est pas $0$, $-2$, $-4$, \dots. Soit $g$
infiniment dérivable, non nulle, à support dans
$]0,+\infty[$, compact et éloigné de l'origine. Par la
méthode de la démonstration de \ref{theo:dansschwartz}
$\cF(D)^{qh} \equiv0$ implique que la fonction $(g*D)(x)$
(de classe $C^\infty$ pour $x\neq0$) est à décroissance
rapide pour $x\to+\infty$. Sa partie quasi-homogène est $\wh
g(w)D^{w,0}$. La transformée de Mellin de $g*D$ est donc
donnée par définition, pour $\Reel(s)\gg0$, par
\[\int_a^\infty \left((g*D)(x)+(-1)\wh
  g(w)\frac{x^{w-1}}{\pi^{-\frac w2}\Gamma(\frac
    w2)}\right)x^{-s}\,dx\; , \]
avec $a>0$ suffisamment petit. Le premier terme donne une
fonction entière et le deuxième terme a comme partie polaire
de son prolongement méromorphe $-\wh g(w)(\pi^{-\frac
  w2}\Gamma(\frac w2))^{-1}\frac1{s-w}$. La partie polaire
de la transformée de Mellin complète de $g*D$ est donc $-\wh
g(w)\frac1{s-w}$ (on savait déjà a priori par le paragraphe
précédent que cette partie polaire était un multiple de
$\frac1{s-w}$, donc malgré la multiplication par
$\pi^{-\frac s2}\Gamma(\frac s2)$, il n'y a pas création de
pôles additionnels en $s=0,-2,-4,\dots$). On peut toujours
choisir $g$ avec $\wh g(w)\neq0$, et on conclut que la
partie polaire de $\pi^{-\frac s2}\Gamma(\frac s2)\wh D(s)$
est exactement $\frac{-1}{s-w}$, lorsque $D^{qh} = D^{w,0}$
et $\cF(D)^{qh} \equiv0$. Ceci a été démontré pour $w\neq
0,-2,-4,\dots$, mais en remplaçant $D$ par $E =
\cF(D-D^{w,0}) = \cF(D) - D^{1-w,0}$, on remplace $D^{w,0}$
par $-D^{1-w,0}$ et $s$ par $1-s$ (par l'équation
fonctionnelle du théorème \ref{theo:S-etendue}), donc la
conclusion vaut aussi sous l'autre condition $w\neq
1,3,5,\dots$. Elle vaut donc pour tout $w\in\CC$.

Soit maintenant $\phi(x)$ une fonction paire non nulle dans
$\cS$ et avec la propriété
$S$-stricte.  Posons, pour $w\in\CC$: $\phi_w(x) =
|x|^{w}x^{-1}\int_0^x |y|^{-w}\phi(y)dy$, de sorte que
$\phi_w$ est à nouveau nulle dans $]-a,a[$ et vérifie
$(\frac d{dx} x - w)\phi_w = \phi$, et donc
$(s-w)\wh{\phi_w}(s) = \wh\phi(s)$. On notera
$\Phi(s):=\pi^{-\frac s2}\Gamma(\frac s2)\wh{\phi}(s)$, qui
est une fonction entière, et $\Phi_w(s) := \pi^{-\frac
  s2}\Gamma(\frac s2)\wh{\phi_w}(s) = \frac{\Phi(s)}{s-w}$.
Considérons, pour $w$ restreint à l'ouvert $U$
complémentaire dans $\CC$ des zéros non triviaux de $\wh
\phi(s)$, la distribution paire
\[ \Delta_w := \frac{\phi_w -\cF^{-1}(\cF(\phi_w)^{qh})}{\Phi(w)}
\;,\]
qui est normalisée de sorte que $\cF(\Delta_w)^{qh}\equiv0$,
et dont la transformée de Mellin est
$\frac1{s-w}\frac{\Phi(s)}{\Phi(w)}$. Par ce qui précède on
sait sans calcul que la partie quasi-homogène de $\Delta_w$
est exactement $-D^{w,0}$, et donc 
l'on peut aussi écrire $\Delta_w = \frac1{\Phi(w)}\phi_w -
D^{w,0}$. Par ailleurs la formule qui définit $\phi_w$
prouve qu'en tant que distribution elle dépend
analytiquement de $w$.  La distribution $\Delta_w$ dépend
donc analytiquement de $w$. Il faut aussi se convaincre que
le passage à la transformée de Mellin est une opération qui
commute à la dérivation un nombre quelconque de fois par
rapport au paramètre $w$. Au facteur
$\frac1{\Phi(w)}\pi^{-\frac s2}\Gamma(\frac s2)$
près il s'agit de regarder, au moins pour $\Reel(s)\gg0$,
$\int_a^\infty \phi_w(x)x^{-s}\,dx$. Certainement on peut
permuter $(\frac\partial{\partial w})^k$ et l'intégrale,
pour tout $k\in\NN$, pour $w$ dans un ouvert d'adhérence
compacte dans $U$ et pour $\Reel(s)\ge\sigma$ avec
$\sigma\gg0$. Par le principe du prolongement analytique,
$(\frac\partial{\partial w})^k
\frac1{s-w}\frac{\Phi(s)}{\Phi(w)}$ est la transformée de
Mellin complète de $(\frac\partial{\partial w})^k \Delta_w$
pour tout $k\in\NN$ (et tout $w\in U$). Or la partie
quasi-homogène de $(\frac\partial{\partial w})^k \Delta_w$
est $-D^{w,k}$, tandis que
\[\frac1{s-w}\frac{\Phi(s)}{\Phi(w)} = \frac1{s-w} +
\frac1{\Phi(w)}\frac{\Phi(s)-\Phi(w)}{s-w}\]
diffère de $\frac1{s-w}$ par une fonction $F(s,w)$
conjointement analytique en $s\in\CC$, $w\in U$. Donc la
partie polaire de $(\frac\partial{\partial w})^k
\frac1{s-w}\frac{\Phi(s)}{\Phi(w)}$, comme fonction de
$s\in\CC$ est exactement $\frac{k!}{(s-w)^{k+1}}$. Ainsi la
partie polaire correspondant à la distribution $D^{w,k}$ est
exactement $-\frac{k!}{(s-w)^{k+1}}$. Ce résultat est établi
pour $w$ distinct des zéros triviaux de $\wh\phi(s)$, et en
changeant $\phi$ il est donc établi pour tout $w\in\CC$.  Le
théorème \ref{theo:polaire} est donc établi dans le cas
pair, et le cas impair est traité quasiment à l'identique.
\medskip

Il peut parfois être plus commode de disposer de la partie
polaire de $\wh D(s)$:

\begin{cor}
  Soit $D$ une distribution non nulle, paire ou impaire,
  avec la propriété $S$-étendue, et normalisée par
  $\cF(D)^{qh}\equiv0$. Soit $w$ distinct de $0$, $-2$,
  $-4$, \dots\ si $D$ est paire ou distinct de $-1$, $-3$,
  $-5$, \dots\ si $D$ est impaire. Soit $(c_0+ c_1 \log(x) +
  c_2 \log^2(x) + \dots + c_N \log^N(x)) x^{w-1} $ la
  $w$-composante de la restriction (quasi-homogène) de $D$ à
  $]0,a[$, $a>0$ suffisamment petit. Alors la $w$-composante
  de la partie polaire de la transformée de Mellin simple
  $\wh D(s)$ est
  \[ \sum_{0\leq k\leq N}  c_{k}\;\frac{(-1) k!}{(s-w)^{k+1}}
  \;.\]
 \end{cor}

Nous traitons le cas pair. Notons $\gamma(s) = \pi^{-\frac
  s2}\Gamma(\frac s2)$. Il suffit de considérer un cas où la
  restriction de $D$ à $]-a,a[$ est $\log^N(|x|) |x|^{w-1}$. La
  composante quasi-homogène est donc
  $(\frac\partial{\partial w})^N \gamma(w)D^{w,0}$, que l'on
  peut écrire
\[ \sum_{0\leq j\leq N}
N!\frac{\gamma^{(N-j)}(w)}{(N-j)!}\frac{D^{w,j}}{j!}\; .\]
Par le théorème précédent la partie polaire de
$\gamma(s)\wh D(s)$ est
\[ \sum_{0\leq j\leq N}
N!\frac{\gamma^{(N-j)}(w)}{(N-j)!}\frac{-1}{(s-w)^{j+1}}\;
,\] 
qui est aussi la partie $w$-polaire du produit 
\[ N! \gamma(s)\frac{-1}{(s-w)^{N+1}} \; . \]
Ainsi $\gamma(s)\Big(\wh D(s) + \frac{N!}{(s-w)^{N+1}}\Big)$ est
holomorphe en $w$ et donc la partie polaire de $\wh D(s)$ en
$w$ est exactement $\frac{-N!}{(s-w)^{N+1}}$, ce qu'il
fallait démontrer.

\subsection{Exemples}

Dans cette dernière section nous évoquerons brièvement
quelques exemples, notables ou anodins, de distributions et
de fonctions de carrés intégrables avec la propriété $S$ ou
la propriété $S$-étendue. À tout seigneur tout honneur notre
premier exemple est probablement le plus important: il
s'agit bien sûr de la distribution de Poisson
\[ D(x) = \sum_{n\in\ZZ} \delta_n(x) \; .\]
Sa transformée de Mellin est $\zeta(s)$. Si l'on remplace
$D$ par $D-1$ de sorte que $\cF(D)^{qh} \equiv 0$, on a
$D^{qh} = \delta_0 - 1 = D^{0,0} - D^{1,0}$. La partie
polaire de $\pi^{-\frac s2}\Gamma(\frac s2)\zeta(s)$ est
donc par le théorème \ref{theo:polaire}: $\frac{-1}s +
\frac1{s-1}$, ce qui n'est certes pas un résultat d'une
fracassante nouveauté. Notons néanmoins que la preuve du
prolongement analytique et de l'équation fonctionnelle de la
fonction $\zeta(s)$ implicite dans cette approche est
nouvelle, et que ses mots-clés sont: régularisation
multiplicative et transformée de Mellin, formule
d'entrelacement (de co-Poisson), opérateur invariant
d'échelle $\cF I$, espaces de Hardy, espaces de de~Branges,
propriété $S$-étendue, fonctions de Nevanlinna.

Citons ensuite la fonction paire de carré intégrable:
\[ f(x) = \frac{\{|x|\}}{|x|} \; , \]
qui vérifie
\[ - \frac{d}{dx}x\, f(x) = \sum_{n\neq 0}\delta_n(x) - 1 =
D_1(x)\;,\] et dont la transformée de Mellin est donc
$-\frac{\zeta(s)}s$ (on laisse au lecteur la vérification de
la conformité de la partie polaire avec le théorème
\ref{theo:polaire}).
L'avantage de la distribution $D_1$ c'est que les
convolutions $g*D_1$ avec des fonctions de classe $C^\infty$
à support compact éloigné de l'origine sont des fonctions de
Schwartz (dont la transformée de Fourier est calculée par la
formule de co-Poisson \ref{eq:copoisson}, c'est-à-dire par
la formule d'entrelacement \ref{theo:entrelacement}).

Des exemples d'un tout autre type sont fournis par les
distributions paires $A_a$ et $B_a$, et leurs analogues
impaires, qui ont été obtenues dans
\cite{crassonine}. L'opérateur sur $L^2(0,a;\,dx)$ donné par
\[ \phi(x)\mapsto (x\mapsto \int_0^a 2\cos(2\pi
xy)\phi(y)\,dy) \]
est un opérateur auto-adjoint compact, et comme il n'existe
pas de fonction non nulle à support compact dont la
transformée de Fourier soit à support compact, la norme de
cet opérateur est strictement inférieure à $1$. Il existe
donc une unique fonction $\phi_a\in L^2(0,a; \,dx)$
vérifiant
\[ \phi_a(x) + \int_0^a 2\cos(2\pi xy)\phi_a(y)\,dy =
2\cos(2\pi ax) \]
Remarquons que l'on peut alors donner un sens par cette
équation à $\phi_a(x)$ pour tout $x\in\RR$ et même pour tout
$x\in\CC$ puisque l'intégrale et le cosinus sont des
fonctions entières de $x$. Donc $\phi_a(x)$ est une fonction
entière, paire, de $x$. Elle est $2\cos(2\pi ax) + O(1/x)$
pour $x\to+\infty$, donc tempérée comme distribution de
$x$. Par l'équation qui la définit sa transformée de Fourier
vérifie l'équation:
\[ \cF(\phi_a)(y) + \Un_{]-a,a[}(y)\phi_a(y) = \delta_a(y) +
\delta_{-a}(y)\; ,  \]
donc:
\[ \cF(\phi_a)(y) + \phi_a(y) = \delta_a(y) +
\delta_{-a}(y) + \Un_{|y|\geq a}(y)\phi_a(y) \;. \]
et ainsi la distribution tempérée paire
\[ A_a = \frac{\sqrt{a}}2\, (\cF(\phi_a) + \phi_a ) \]
est identiquement nulle dans $]-a,a[$. Comme elle est sa
propre transformée de Fourier elle a la propriété $S$
stricte (le facteur $\sqrt a$ est introduit pour d'autres
raisons; \textsl{cf.} \cite{crassonine, crasdirac}).
De même en partant de 
\[ \phi_a^-(x) - \int_0^a 2\cos(2\pi xy)\phi_a^-(y)\,dy =
2\cos(2\pi ax) \]
et en définissant
\[ B_a = \frac{i\sqrt{a}}2\, (\cF(\phi_a^-) - \phi_a^-) \]
on obtient une autre distribution tempérée paire,
anti-invariante sous Fourier, avec la propriété $S$ stricte
pour l'intervalle $]-a,a[$.

Les fonctions 
\[ \cA_a(s) = \pi^{-\frac s2}\Gamma(\frac
s2)\wh{A_a}(s)\qquad \cB_a(s) = \pi^{-\frac s2}\Gamma(\frac
s2)\wh{B_a}(s) \]
sont donc par le théorème \ref{theo:mellindistribsonine} des
fonctions entières vérifiant les équations fonctionnelles
\[ \cA_a(s) = \cA_a(1-s)\;,\qquad \cB_a(s) = - \cB_a(1-s)\; ,\]
et on notera aussi  que
\[ \cA_a(s) = \overline{\cA_a(\overline s)} \;,\qquad 
\cB_a(s) = - \overline{\cB_a(\overline s)} \; . \]

Nous avons le théorème suivant (et des définitions analogues
avec la fonction sinus mènent à des distributions impaires
avec la propriété $S$ pour l'intervalle $]-a,a[$ qui
vérifient aussi ce théorème):

\begin{theo}[{\cite{crassonine}}]\label{theo:riemannAB}
  Les fonctions entières $\wh{A_a}(s)$ et $\wh{B_a}(s)$
  ont tous leurs zéros non-triviaux sur la droite critique.  
\end{theo}

Cela résulte des équations fonctionnelles et du fait
suivant: la fonction entière $\cE_a(s) = \cA_a(s) - i
\cB_a(s)$ vérifie $|\cE_a(s)|>|\cE_a(\overline{1-s})|$ pour
$\Reel(s)>\frac12$, car $(|\cE_a(s)|^2 -
|\cE_a(\overline{1-s})|^2)/(2\Reel(s)-1)$ est le carré de la
norme de l'évaluation en $s$ de $\pi^{-\frac s2}\Gamma(\frac
s2)\wh f(s)$ pour $f$ paire, de carré intégrable, nulle et
de Fourier nulle dans $]-a,a[$.  Que de telles transformées
de Mellin (à un changement de variable près) forment un
espace de de~Branges au sens de \cite{bra} avait été
démontré par de~Branges dans \cite{bra64}, et de cela
découlait déjà l'existence d'une fonction entière $\cE_a$
donnant la norme des évaluateurs par la formule ci-dessus,
sans toutefois que l'on puisse en donner de formule
explicite.  Les distributions $A_a$ et $B_a$ ont été
obtenues par l'auteur dans son article \cite{crassonine} qui
résolvait ce problème de donner une formule \og
explicite\fg\ pour $\cE_a$.

On peut aussi prouver \ref{theo:riemannAB}, en donnant de
plus une interprétation opératorielle aux zéros, en
utilisant les équations différentielles par rapport à $a$
qui ont été établies par l'auteur dans \cite{crasdirac} et
qui donnent la déformation des espaces de Hilbert de
de~Branges-Rovnyak-Fourier par rapport à $a>0$.

Rappelons de plus que comme conséquence du théorème
\ref{theo:mellindistribsonine} (dont la preuve a fait appel
à \cite[Theorem 7.7]{twosys}) le nombre des zéros de parties
imaginaires entre $0$ et $T$ est équivalent à
$T\log(T)/2\pi$.

Revenons pour conclure à des exemples déduits de la
distribution de Poisson (ou de la fonction dzêta de
Riemann):

\begin{fact}
Soit $0<a<1$, $A=1/a$. La fonction paire donnée pour 
$x>0$ selon :
  \[f_a(x) = \sum_{ax\leq n\leq Ax} \frac{3(\frac xn + \frac
  nx) - (A+a+4)}{\sqrt{nx}}\]
est de carré intégrable sur $\RR$, vérifie $\cF(f_a) = +
f_a$, et est identiquement nulle sur l'intervalle $]-a,a[$.
\end{fact}

\begin{fact}
Soit $0<a<1$, $A=1/a$. Soient $f_a(x)$,
$g_a(x)$, $k_a(x)$ les fonctions impaires données pour $x>0$
selon:
\begin{eqnarray*}
f_a(x) &=& \sum_{a x\leq n +\frac12 \leq
A x}\frac{(-1)^n}{x}\\ g_a(x) &=& \sum_{a x\leq n
+\frac12 \leq A x}\frac{(-1)^n}{n+\frac12}\\ k_a(x) &=&
\sum_{a x\leq n +\frac12 \leq
A x}\frac{(-1)^n}{\sqrt{(n+\frac12)x}}\\
\end{eqnarray*}
Elles sont de carrés intégrables sur $\RR$ et identiquement
nulles sur $]-\frac a2,+\frac a2[$. La transformée en sinus de $f_a$ est
$g_a$. La fonction $k_a$ est sa propre transformée en sinus.
\end{fact}

\begin{fact}
Soit $N\in\NN$. Soit 
\[ Q_N(n) = \prod_{0\leq j < N} \left(\strut n(n+1)-
j(j+1)\right)\]
La fonction paire (non-identiquement nulle!!) qui est donnée
pour $x>0$ par la formule:
  \[\sum_{\stackrel{n\ \text{tel que}}{x\leq \frac{n +
        \frac12}{\sqrt{N+\frac12}} \leq 2x}} 
(-1)^n Q_N(n)\left(\Bigg(\frac{x}{n+\frac12} -
\frac1{\sqrt{4N+2}}\Bigg)\Bigg(\frac{x}{n+\frac12} -
\frac1{\sqrt{N+\frac12}}\Bigg)\right)^{2N+1}\]
est de carré intégrable sur $\RR$. Elle-même et sa
transformée de Fourier sont identiquement nulles dans
l'intervalle $]-\frac12(N+\frac12)^{1/2},
+\frac12(N+\frac12)^{1/2}\strut[$.
\end{fact}

Nous avons donné ces dernières assertions purement à titre
d'exercices.

\bigskip
\bigskip

\begin{small}

\noindent\textbf{Remerciements.} Je remercie {Jean-Pierre
  Kahane}: dans sa lettre (22 mars 2002) il m'indiquait une
construction explicite de fonctions dans $\cS$ avec la
propriété de support pour un intervalle $]-a,a[$ quelconque,
par modification puis régularisation par convolution
additive de la distribution de Poisson. En étudiant la
régularisation multiplicative, et les propriétés des
fonctions analytiques obtenues par transformation de Mellin,
j'avais alors établi plusieurs des résultats qui ont été
présentés ici et qui furent déjà employés dans
\cite{crassonine}. Je remercie {Luis B\'aez-Duarte}:
plusieurs échanges depuis décembre 2001 sur la formule de
co-Poisson, en particulier sur l'intérêt de résultats dans
un style \og classique\fg\ (tel que celui de l'ouvrage de
Titchmarsh \cite{titchfourier}), m'ont incité à obtenir
plusieurs parmi les propositions qui ont été présentées ici.

\medskip Le présent manuscrit est celui qui a été annoncé
par l'auteur dans de précédentes publications sous le titre
\og Co-Poisson intertwining: distribution and function
theoretic aspects\fg .

\end{small}

\bigskip

\bibliographystyle{amsplain}

\vfill

\begin{small}
\textsc{Université Lille 1, UFR de Mathématiques, Cité
  Scientifique M2, F-59655 Villeneuve d'Ascq cedex, France}

\texttt{burnol@math.univ-lille1.fr}
\end{small}

\clearpage





31 août 2004. Modifications par rapport à la version du
1\ier\ août 2004:

\begin{itemize}
\item page 8: on a rajouté \og(plus précisément, des sommes
  $\sum_n \frac{c_n}n F(x^{\frac1n})$)\fg.
\item page 9: on rajoute à l'introduction l'énoncé du
  théorème inverse qui caractérise les transformées de
  Mellin de distributions avec la propriété $S$ ou
  $S$-étendue.
\item page 12: notation: $\mu([0,t[)$ au lieu de $\mu([0,t))$.
\item page 48: notations: $\frac\partial{\partial x}$ et
  $\frac\partial{\partial t}$ dans la preuve de 4.8.
\item page 56: dans la preuve de 4.37 on a utilisé
  $O((1+|s-\frac12|)^{|s-\frac12|})$ pour un passage à
  $1-s$ plus correct.
\item pages 61-62: on rajoute au
  théorème \ref{theo:S-etendue}: \og Les seules
  autres distributions \dots \dots On a 
  $\wh{\cF(D)}(s) = \chi(s)\wh D(1-s)$.\fg
\item pages 62 et suiv.: on a ajouté à la section 4.G le
  théorème \ref{theo:polaire} et son corollaire décrivant
  les parties polaires.
\item On ajoute quelques lignes au début de la section
  4.H. illustrant le théorème sur les parties polaires.
\item On donne plus d'informations relatives au
  théorème \ref{theo:riemannAB} de la section 4.H.
\end{itemize}

\end{document}